\documentclass[]{article}
\usepackage[utf8]{inputenc} 
\usepackage[T1]{fontenc}    
\usepackage{hyperref}       
\usepackage{url}            
\usepackage{graphics}
\usepackage{graphicx}
\usepackage{epsfig}
\usepackage{array}
\usepackage{color} 
\usepackage{tabularx}
\usepackage{amsmath}
\usepackage{amssymb}
\usepackage{amsfonts}
\usepackage{txfonts}
\usepackage{arxiv}
\usepackage{xcolor}
\usepackage{bm}
\usepackage{soul}
\usepackage{float}
\usepackage{lineno}
\usepackage{cite}
\usepackage{booktabs, multirow}
\usepackage{textcomp}
\usepackage{epstopdf}
\usepackage{hhline}

\hypersetup{
	colorlinks=true,                          
	linkcolor=blue, 
	citecolor=blue, 
	urlcolor=black  
} 
\usepackage[colorinlistoftodos]{todonotes}
\allowdisplaybreaks

\title{Semi-implicit hybrid finite volume/finite element method for the GPR model of continuum mechanics}

\author{
	\href{https://orcid.org/0000-0002-6509-4269}{\includegraphics[scale=0.06]{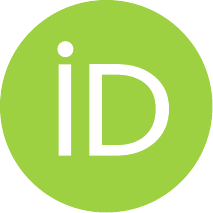}\hspace{1mm}Saray Busto}
	\thanks{corresponding author.\newline Email addresses: \texttt{saray.busto.ulloa@usc.es} (S. Busto) 
	\texttt{laura.delrio@unitn.it} (L. R\'io-Mart\'in)}$^{\;\;a}$ 
	\hspace{3mm}
	\href{https://orcid.org/0000-0002-7812-3161}{\includegraphics[scale=0.06]{orcid.pdf}\hspace{1mm}Laura R\'io-Mart\'in$^{b}$}
}

\affiliation{$^a$  Department of Applied Mathematics, Universidade de Santiago de Compostela, 15782 Santiago de Compostela, Spain; \\Galician Center for Mathematical Research and Technology, CITMAga, 15782 Santiago de Compostela, Spain\\
	$^b$Department of Information Engineering and Computer Science, University of Trento, via Sommarive 9, Povo, 38123 Trento, Italy; \\Laboratory of Applied Mathematics, DICAM, University of Trento, Via Mesiano 77, 38123 Trento, Italy.}
\date{}

\renewcommand{\headeright}{S. Busto, L. R\'io-Mart\'in}
\renewcommand{\undertitle}{}
\renewcommand{\shorttitle}{}

\hypersetup{
	pdftitle={},
	pdfsubject={},
	pdfauthor={Saray Busto, Laura del Rio Martin},
	pdfkeywords={First keyword, Second keyword, More},
}
\DeclareMathOperator{\dive}{\nabla \cdot}
\DeclareMathOperator{\gra}{\nabla}

\DeclareMathOperator{\lapla}{\Delta}

\DeclareMathOperator{\dS}{\mathrm{dS}}
\DeclareMathOperator{\dV}{\mathrm{dV}}

\newcommand{\halb}{\frac{1}{2}}

\newcommand{\Dt}{\Delta t\,}
\newcommand{\x}{\mathbf{x}}
\newcommand{\ci}{{\mathfrak{i}}}
\newcommand{\cell}[1]{C_{#1}}
\newcommand{\vcell}[1]{\left|\cell{#1}\right|}
\newcommand{\cj}{{\mathfrak{j}}}
\newcommand{\cl}{{\mathfrak{l}}}
\newcommand{\ck}{\boldsymbol{\kappa}}
\newcommand{\cm}{\mathfrak{m}}
\newcommand{\Pcell}[1]{T_{#1}}
\newcommand{\nn}{\boldsymbol{\eta}}
\newcommand{\un}{\mathbf{n}}
\newcommand{\vel}{u}
\newcommand{\mom}{\rho\vel}
\newcommand{\tmom}{\rho\vel^{\ast}}
\newcommand{\bvel}{\mathbf{\vel}}
\newcommand{\bmom}{\rho\bvel}
\newcommand{\tbmom}{\rho\bvel^{\ast}}
\newcommand{\press}{p}
\newcommand{\tpress}{\widetilde{\press}}
\newcommand{\A}{A}
\newcommand{\bA}{\mathbf{\A}}
\newcommand{\bG}{\mathbf{G}}
\newcommand{\devG}{\mathring{\bG}}
\newcommand{\J}{J}
\newcommand{\bJ}{\mathbf{\J}}
\newcommand{\bY}{\mathbf{Y}}

\newcommand{\hq}{q}
\newcommand{\g}{\mathbf{g}}
\newcommand{\q}{\mathbf{q}}
\newcommand{\Q}{\mathbf{Q}}
\newcommand{\eQ}{\overline{\Q}}
\newcommand{\tQ}{\Q^{\ast}}
\newcommand{\Flux}{\boldsymbol{\mathcal{F}}}
\newcommand{\source}{\boldsymbol{\mathcal{S}}}
\newcommand{\BNCP}[1]{\boldsymbol{\mathcal{B}}(#1)}
\newcommand{\NCP}[1]{\BNCP{#1}\cdot \gra #1}
\newcommand{\E}{E}
\newcommand{\bE}{\mathbf{E}}

\newcommand{\bsigma}{\boldsymbol{\sigma}}
\newcommand{\bomega}{\boldsymbol{\omega}}

\newcommand{\pxk}{\frac{\partial}{\partial {x_{k}}}}
\newcommand{\pxj}{\frac{\partial}{\partial {x_{j}}}}
\newcommand{\pxi}{\frac{\partial}{\partial {x_{i}}}}
\newcommand{\pxm}{\frac{\partial}{\partial {x_{m}}}}
\newcommand{\pt}{\frac{\partial}{\partial t}}

\newtheorem{weakproblem}{Weak problem}

\begin{document}
\maketitle
\begin{abstract}
We present a new hybrid semi-implicit finite volume / finite element numerical scheme for the solution of incompressible and weakly compressible media. From the continuum mechanics model proposed by Godunov, Peshkov and Romenski (GPR), we derive the incompressible GPR formulation as well as a weakly compressible GPR system. As for the original GPR model, the new formulations are able to describe different media, from elastoplastic solids to viscous fluids, depending on the values set for the model's relaxation parameters. Then, we propose a new numerical method for the solution of both models based on the splitting of the original systems into three subsystems: one containing the convective part and non-conservative products, a second subsystem for the source terms of the distortion tensor and heat flux equations and, finally, a pressure subsystem. In the first stage of the algorithm, the transport subsystem is solved by employing an explicit finite volume method, while the source terms are solved implicitly. Next, the pressure subsystem is implicitly discretised using finite elements. Within this methodology, unstructured grids are employed, with the pressure defined in the primal grid and the rest of the variables computed in the dual grid. 
To evaluate the performance of the proposed scheme, a numerical convergence analysis is carried out, which confirms the second order of accuracy in space. A wide range of benchmarks is reproduced for the incompressible and weakly compressible cases, considering both solid and fluid media. These results demonstrate the good behaviour and robustness of the proposed scheme in a variety of scenarios and conditions. 
\end{abstract}



\section{Introduction} \label{sec:intro}
Incompressible and weakly compressible flows appear in numerous industrial and biological applications, ranging from fluid dynamics in hydraulic systems to the study of airflow in aeronautical engineering. The study of such flows is crucial because of their prevalence in situations where the flow velocity is low compared to the speed of sound, making compressibility effects small but generally not negligible. To accurately model these flows, a classical approach is the use of Navier-Stokes equations, whose numerical solution has been extensively studied. A recently developed alternative is the use of the Godunov-Peshkov-Romenski (GPR) model for continuum mechanics, which offers a unified formulation for modelling different media, from solids with large deformations to compressible viscous fluids. The use of this model allows capturing complex phenomena of continuum mechanics, providing greater versatility in the simulations than former classical approaches.

The main objective of this work is to present a novel methodology for addressing problems in the weakly compressible regime by solving incompressible and weakly compressible formulations of the GPR first-order hyperbolic model of continuum mechanics \cite{DPRZ16,PeshRom2016,Rom1998}. The compressible GPR model, which employs a single system of hyperbolic equations, describes both solids and fluids in a unified manner by an appropriate choice of the model relaxation parameters. This model originates from the elastoplastic deformation model of Godunov and Romenski \cite{Godunov:2003a}, and was first introduced in \cite{PeshRom2016}. Meanwhile, the heat equation is derived from the model proposed by Malyshev and Romenski in \cite{MalyshevRom1986} and implemented in \cite{Rom2007TwoPhase,Rom2010TwoPhase}. 

Since its introduction, the compressible GPR model has been successfully applied to model a variety of problems, including non-Newtonian and viscous flows, and elastic, elastoplastic and porous solids \cite{DPRZ16,BDL16,PDBRCI2021GPRNN,Dumbser2018viscous,Peshkov2019hyperhypo,Romenski2021threephase,Romenski2021porous}. In addition, it has been extended to include the effects of electrodynamics \cite{DPRZ2017MHD}, surface tension \cite{Chiocchetti2022} or general relativity \cite{Romenski2020GR}, as well as to model nonlinear dispersive systems \cite{Dhaouadi2022NSKGPR} and to include solids or viscous fluids in multiphase systems \cite{Ferrari2024unified}.
Different types of schemes have been successfully developed to solve these equations as explicit finite volume (FV) methods \cite{Dhaouadi2023NSKSP}, semi-implicit second order finite volume schemes \cite{Boscheri2021SIGPR,PDBRCI2021GPRNN}, high order IMEX methods \cite{Boscheri2022IMEX}, high order ADER-FV and ADER-DG approaches \cite{DPRZ16,BDL16,DPRZ2017MHD,BCDGP20}, or smooth particle methods (SPH) \cite{Kincl2023}.

It is noteworthy to remark that the GPR model presents at the continuous level, a wide set of properties that should be conveyed to the discrete level. 
Firstly, the model falls within the framework of Symmetric Hyperbolic Thermodynamically Compatible (SHTC) systems so that thermodynamically compatibility at the discrete level may also be pursued, as done in the modern family of HTC-FV and HTC-DG schemes \cite{Busto2022HTCGPR,Busto2022HTCDG,HTCAbgrall,HTCTotalAbgrall}. 
Moreover, the model presents natural involution constraints on the curl-free property of the distortion and heat flux fields that have been addressed, e.g., in \cite{Boscheri2021SIGPR}. We furthermore highlight the asymptotic preserving property of the model in the fluid relaxation limit, i.e., as the corresponding relaxation times go to zero, the Navier-Stokes-Fourier limit is obtained, \cite{DPRZ16}. 
To the best of the authors knowledge, simultaneous preservation of all these properties at the discrete level has not been achieved yet, and the effort is first placed on developing efficient schemes preserving at least some of these properties. In particular, the methodology proposed in this paper falls in the family of asymptotic preserving schemes in the fluid limit of the model. Furthermore, analogously to the asymptotic preserving properties of the Navier-Stokes equations in the incompressible limit, considering the low Mach limit we recover an incompressible formulation of the GPR model.
%

Let us note that addressing low Mach number flows using the Navier-Stokes equations is a wide field of research where two main approaches have been initially followed: the development of explicit density-based and pressure-based solvers. The first family was first proposed in the framework of high Mach number flows and has as its main shortcoming the highly restrictive time step condition related to the pressure wave velocity in the low Mach number limit. Moreover, unless properly corrected, those schemes may present excessive numerical diffusivity due to an incorrect scaling with respect to the Mach number, \cite{GM04,Dellacherie1,LG08}.
On the other hand, the second family has been initially designed in the context of incompressible flows and is usually based on a non-conservative formulation of the equations leading to important errors in the presence of strong discontinuities, \cite{CG84,MDR07}. Aiming at profiting from the main advantages of both families, the first semi-implicit pressure-based solvers were proposed in \cite{PM05,HVW03}. Since then, they have been successfully extended to different families of numerical methods including, e.g., finite volume approaches (FV), \cite{CordierDegond,DC16,RussoAllMach,BDLTV2020,Thomann2020,BP2021}, discontinuous Galerkin schemes (DG), \cite{TD17,BTBD20,BT22,TB24}, or hybrid finite volume/finite element methods (hybrid FV/FE), \cite{BFSV14,BFTVC17,Hybrid1,BRMVCD21} and hybrid finite volume/virtual element schemes (hybrid FV/VEM), \cite{HybridFVVEMinc,BBD24}. 

The same advantages and drawbacks arising in the discretisation of low Mach number flows when using the Navier-Stokes equations could also be expected for discretising the weakly compressible GPR model. Therefore, accounting for the promising results obtained in the framework of the Navier-Stokes equations and also in previous semi-implicit schemes for continuum mechanics models, \cite{AbateIolloPuppo,Boscheri2021SIGPR}, we focus on the development of a semi-implicit scheme for the incompressible and weakly compressible GPR models. In particular, the methodology presented in this paper follows the seminal ideas in \cite{BFSV14} and is based on a semi-implicit hybrid scheme that combines finite volume and finite element methods, taking advantage of the benefits of both numerical approaches. This family of semi-implicit hybrid FV/FE schemes has been developed in the last decade aiming at solving the Navier-Stokes equations, both for incompressible and all Mach number flows as well as the shallow water and the MHD equations, on fixed and moving staggered unstructured grids in two and three spatial dimensions \cite{BFTVC17,BRMVCD21,HybridMPI,BD22,HybridALE,HybridHexa1,HybridMHD,HybridImplicit}. The procedure behind this methodology consists in splitting the system into two parts which decouples the pressure field from the convective system. This division of the original system into subsystems is performed following the Toro-V\'azquez splitting technique \cite{TV12}. Then, a second-order explicit finite volume discretisation is used for the convective terms, while a second-order continuous Lagrangian finite element scheme is employed to solve the pressure subsystem. Moreover, within this methodology, unstructured staggered grids are considered, in a similar way to that described in \cite{TD16,TD17,BTBD20}. Since we focus on the weakly compressible regime, the semi-implicit approach further improves the stability properties and the efficiency of the scheme with respect to classical fully explicit approaches.

This paper is organized as follows. Section~\ref{sec:goveq} presents a brief review of the original GPR model, followed by the novel derivation of the related incompressible and weakly compressible GPR formulations. Section~\ref{sec:numdisc} details the proposed numerical discretisation for both systems. In particular, the splitting performed to get the transport and the Poisson-type pressure subsystem is described. Then, the algorithms used to solve both subsystems are detailed, as well as the interpolation approach performed between the dual and the primal grids. 
Numerical validation of the proposed methodology is presented in Section \ref{sec:numericalresults} for two and three dimensions. The wide set of benchmarks studied include both incompressible and weakly compressible cases, considering both solid and fluid media. The paper closes with some remarks and an outlook for future work in Section \ref{sec:conclusions}.

\section{Governing equations} \label{sec:goveq}
Before introducing the incompressible and the weakly compressible GPR models, we start recalling the original compressible Godunov-Peshkov-Romenski  model for continuum mechanics, \cite{DPRZ16,PeshRom2016,Dumbser2018viscous}, that reads
\begin{subequations}\label{eqn.GPR}
	\begin{eqnarray}
		\pt \rho + \pxk \left(\rho\vel_{k}\right) &=& 0, \label{eqn.GPR_rho}\\
		\pt \left(\rho\vel_i\right) + \pxk \left(\rho \vel_{i} \vel_{k} \right) + \pxi \press +\pxk\left( \sigma_{ik}+\omega_{ik} \right) &=& \rho g_{i}, \label{eqn.GPR_mom}\\
		\pt \A_{ik} + \pxk \left(\vel_{m}\A_{im}\right) + \vel_{j}\left(\pxj \A_{ik} -\pxk \A_{ij}\right) &=& -\frac{1}{\theta_{1} \left(\tau_{1}\right)} \E_{\A_{ik}},\qquad  \label{eqn.GPR_A}\\
		\pt \J_{k} + \pxk \left(\J_{m}\vel_{m} \right) +\pxk T + \vel_{j}\left(\pxj J_{k} - \pxk J_{j}\right) &=&  -\frac{1}{\theta_{2} \left(\tau_{2}\right)} \E_{\J_{k}},\qquad  \label{eqn.GPR_J}\\
		\pt \left( \rho S\right) + \pxk \left(\rho S\vel_{k}\right) + \pxk\E_{\J_{k}} = \frac{\rho}{T}\left(\frac{1}{\theta_{1}\left(\tau_{1}\right)}\E_{\A_{ik}}\E_{\A_{ik}} +  \frac{1}{\theta_{2}\left(\tau_{2}\right)}\E_{\J_{k}}\E_{\J_{k}}\right) &\geq& 0, \label{eqn.GPR_S}\\
		\pt \left( \rho \E\right)  +\pxk \left(\rho\E\vel_{k}\right) + \pxk \left(\press\vel_{k}\right) + \pxk \left(\vel_{i}\sigma_{ik} \right)+ \pxk \left(\vel_{i}\omega_{ik} \right) +\pxk \hq_{k} &=& \rho g_{i} \vel_{i}. \label{eqn.GPR_E}
	\end{eqnarray}
\end{subequations}
In this paper, we apply the Einstein summation notation for repeated indexes. Moreover, we denote $\rho$ the density, $\bvel=\left(\vel_{1},\vel_{2},\vel_{3}\right)$ the velocity vector and its components, $\press$ the pressure, $\bA=\left( A_{ik}\right) $ the distortion field, given by a $3\times 3$ tensor, $\bJ=\left(\J_{1},\J_{2},\J_{3}\right)$ the heat flux vector, $S$ the entropy, $E$ the total energy, to be further described later,  and $\g=\left(g_{1},g_{2},g_{3}\right)$ the gravity vector. 
The non-isotropic part of the stress tensor, containing the shear and thermal stresses, is given by:
\begin{equation}
	\sigma_{ik} = \A_{ji}\, \partial_{A_{jk}}\!\! \left( \rho E\right)   = \rho c_{s}^2 G_{ij} \mathring{G}_{jk}, \qquad 
	\omega_{ik} = \J_{i}\, \partial_{J_k}\!\! \left( \rho E  \right)     = \rho c_{h}^{2} \J_i\J_{k},
	\label{eqn.def_sigma_omega}
\end{equation}
where $c_{s}^{2}$ and $c_{h}^{2}$ are the characteristic velocities for propagation of shear and thermal perturbations while $\mathring{G}_{ik}$ denotes the trace-free part of the metric tensor $G_{ik}=\A_{ji}\A_{jk}$,
\begin{equation*}
	\mathring{G}_{ik} = G_{ik}- \frac{1}{3}G_{mm}\delta_{ik}. 
\end{equation*}
The heat flux is
\begin{equation}
\hq_{k} = \partial_{\rho S} \E \partial_{J_k} \E = \rho c_{h}^2 T \J_k, \label{eqn.def_q}
\end{equation}
and $T$ corresponds to the temperature
\begin{equation*}
	T = \partial_{S}\E.
\end{equation*}
Furthermore, the shear and thermal stress relaxation functions read
\begin{equation}
	\theta_{1}\left(\tau_{1}\right) = \frac{1}{3} \rho_{0}\tau_{1} c_{s}^2 \left|\bA\right|^{-\frac{5}{3}},\qquad
	\theta_{2}\left(\tau_{2}\right) = \frac{\rho_{0} T_{0}}{T} \tau_{2} c_{h}^{2},
	\label{eqn.relaxation_fun}
\end{equation}
with $\tau_{1}$ and $\tau_{2}$ the corresponding relaxation times. 
The total energy $\E$ can be divided into four terms as
\begin{equation}
	\E\left(\rho,\bvel,\bA,\bJ,S\right) = \E_{1}\left(\bvel\right)
	  +\E_{2}\left(\bA\right)  
	  +\E_{3}\left(\bJ\right)  
	  +\E_{4}\left(\rho,S\right),\label{eqn.energydecomp}
\end{equation}
while in case the heat flux contributions are neglected, the term $\E_3$ is no more taken into account.
The first contribution to the energy, $\E_{1}$, corresponds to the specific kinetic energy per unit mass,
\begin{equation}
	\E_{1}\left(\bvel\right) = \halb \left|\bvel\right|^{2}.\label{eqn.E1}
\end{equation}
The second and third terms provide the contribution of the mesoscopic, non-equilibrium, part of the total energy related to the material deformations and the thermal impulse,
\begin{equation}
	\E_{2}\left(\bA\right) = \frac{1}{4} c_{s}^2 	\mathring{G}_{ij}	\mathring{G}_{ij}, \qquad \E_{3}\left(\bJ\right) = \frac{1}{2} c_{h}^2 	\J_{i}	\J_{i}.\label{eqn.E23}
\end{equation}
The last term in \eqref{eqn.energydecomp} is the internal energy, related to the kinetic energy of the molecular motion, that we assume given by the ideal gas equation of state
\begin{equation*}
	\E_{4}\left(\rho,S\right) = \frac{\rho^{\gamma-1}}{\left( \gamma-1\right) } e^{\frac{S}{c_{v}}} , \label{eqn.E4s}
\end{equation*}
which is equivalent to consider
\begin{equation}
	\E_{4}\left(\rho,\press\right) = \frac{\press}{\rho\left( \gamma-1\right) }, \label{eqn.E4}
\end{equation}
with $\gamma=\frac{c_{\press}}{c_{v}}$ the ratio of specific heat at constant pressure, $c_{\press}$, and at constant volume, $c_{v}$.

\subsection{Incompressible GPR model}
To get the simplified GPR model in the incompressible limit, we proceed as for the incompressible Navier-Stokes equations. Accordingly, we suppose the fluid to be incompressible, homogeneous and non heat-conducting, obtaining the following system of conservation laws:

\begin{subequations}\label{eqn.incGPR}
	\begin{eqnarray}
		\pxk \left( \mom_{k}\right)  &=& 0,\\
		\pt \left(\rho\vel_i\right) + \pxk \left(\rho \vel_{i} \vel_{k} \right) + \pxi \press +\pxk\sigma_{ik} &=& \rho g_{i},\\
		\pt \A_{ik} + \pxk \left(\vel_{m}\A_{im}\right) + \vel_{j} \left(\pxj \A_{ik} -\pxk \A_{ij}\right) &=& -\frac{1}{\theta_{1} \left(\tau_{1}\right)} \E_{\A_{ik}},
	\end{eqnarray}
\end{subequations}
where we have furthermore assumed a constant density.

\subsection{Weakly compressible GPR model}
On the other hand, the pressure-based reformulation of the energy equation would provide a PDE system able to solve weakly compressible flows. We start multiplying the momentum equations by the corresponding dual variable $\vel_{i}$ yielding
\begin{gather*}
	\vel_{i} \pt \left(\rho\vel_i\right) + \vel_{i} \pxk \left(\rho \vel_{i} \vel_{k} \right) + \vel_{i} \pxi \press +\vel_{i} \pxk\left( \sigma_{ik}+\omega_{ik} \right) = \rho g_{i}\vel_{i},\\
	\pt \left(\halb \rho\vel_i^2\right) + \pxk \left(\halb\rho \vel_{i}^2 \vel_{k} \right) + \vel_{i}\pxi \press +\vel_{i}\pxk\left( \sigma_{ik}+\omega_{ik} \right) = \rho g_{i}\vel_{i}.
\end{gather*}
Summing all momentum equations we get
\begin{gather*}
	\pt \left(\halb \rho\left| \bvel \right|^2\right) + \dive \left(\halb\rho \left| \bvel \right|^2 \bvel \right) + \bvel\cdot \gra \press +\bvel \cdot\dive \left( \bsigma+\bomega \right) = \rho \g\cdot\bvel.
\end{gather*}
Subtracting this relation from the total energy equation and taking into account that
\begin{equation*}
	\dive \left(\bsigma \bvel \right) = \bvel\cdot \dive \bsigma + \bsigma\cdot\gra\bvel, \quad
	\dive \left(\bomega \bvel \right) = \bvel\cdot \dive \bomega + \bomega\cdot\gra\bvel, \quad
	\dive \left(\bvel \press\right)   = \press \dive \bvel + \bvel\cdot\gra\press,
\end{equation*}
it results
\begin{gather*}
	\pt \left( \rho \E\right) 	-\pt \left(\halb \rho\left| \bvel \right|^2\right)
	+ \dive \left(\rho\E\bvel\right) -\dive \left(\halb\rho \left| \bvel \right|^2 \bvel \right)
	+ \press \dive \bvel 
	+ \bsigma\cdot\gra\bvel   
	+ \bomega\cdot\gra\bvel  
	+ \dive \q  = 0.
\end{gather*}
Decomposing the total energy \eqref{eqn.energydecomp} into its four components, and taking into account \eqref{eqn.E1} and the equation of state for ideal gasses in $E_4$, \eqref{eqn.E4}, we get
\begin{gather*}
	\pt \left( \frac{\press}{\gamma-1}\right) + \dive \left( \frac{\press}{ \gamma-1} \bvel\right)
	+ \pt \left( \rho \E_2\right)
	+ \dive \left(\rho\E_2\bvel\right) 
	+ \pt \left( \rho \E_3\right)
	+ \dive \left(\rho\E_3\bvel\right) 
	\notag\\	
	+ \press \dive \bvel
	+ \bsigma \cdot \gra \bvel
	+ \bomega \cdot \gra \bvel
	+ \dive \q  = 0.
\end{gather*}
Since $\gamma = \dfrac{c^2\rho}{\press}$ and $\dive \left( \press\bvel\right) = \press\dive \bvel+\bvel\cdot\gra\press$, we obtain
\begin{gather*}
	\frac{\partial \press}{\partial t}  
	+ \bvel\cdot\gra\press
	+ \left( \gamma-1\right) \left( \pt \left( \rho \E_2\right)
	+ \dive \left(\rho\E_2\bvel\right) 
	+ \pt \left( \rho \E_3\right)
	+ \dive \left(\rho\E_3\bvel\right) \right) 
	\notag\\
	+ c^2\rho \dive \bvel
	+ \left( \gamma-1\right) \left( \bsigma \cdot \gra \bvel
	+ \bomega \cdot \gra \bvel
	+ \dive \q\right)   = 0.
\end{gather*}
Hence, taking into account $ \dive \left( \rho\bvel \right) = \rho \dive \bvel + \bvel\cdot \gra \rho$, it yields
\begin{gather}
	\frac{\partial \press}{\partial t}  
	+ \bvel\cdot\gra\press
	+ \left( \gamma-1\right) \left( \pt \left( \rho \E_2\right)
	+ \dive \left(\rho\E_2\bvel\right) 
	+ \pt \left( \rho \E_3\right)
	+ \dive \left(\rho\E_3\bvel\right) \right) 	
	\notag\\
	+ c^2 \dive \left( \rho\bvel \right) - c^2 \bvel \cdot\gra \rho
	+ \left( \gamma-1\right) \left( \bsigma \cdot \gra \bvel
	+ \bomega \cdot \gra \bvel
	+ \dive \q\right)   = 0. \label{eqn.pressder1}
\end{gather}
To get rid of the time derivative term on $\E_2$, we first multiply \eqref{eqn.GPR_A} by $\rho\E_{A_{ik}}$, and sum over all equations of the distortion field components, obtaining
\begin{equation*}
	\rho\E_{\A_{ik}}\pt \A_{ik} + \rho\E_{\A_{ik}}\pxk \left(\vel_{m}\A_{im}\right) + \rho\E_{\A_{ik}}\vel_{j} \pxj \A_{ik} -\rho\E_{\A_{ik}}\vel_{j}\pxk \A_{ij} = -\frac{\rho}{\theta_{1} \left(\tau_{1}\right)} \E_{\A_{ik}}\E_{\A_{ik}},
\end{equation*}
hence 
\begin{equation*}
	\rho\pt \E_2 + \rho\bvel\cdot \gra \E_2 + \rho\E_{\A_{ik}}\left( \pxk \left(\vel_{m}\A_{im}\right)  -\vel_{m}\pxk \A_{im}\right)  = -\frac{\rho}{\theta_{1} \left(\tau_{1}\right)} \E_{\A_{ik}}\E_{\A_{ik}}.
\end{equation*}
Adding \eqref{eqn.GPR_rho} multiplied by $\E_2$, 
\begin{equation*}
	\E_2 \pt \rho + \E_2 \dive \left(\rho\bvel\right) = 0,
\end{equation*}
we get
\begin{equation}
	\pt\left( \rho\E_2\right) + \dive \left(\rho \E_2 \bvel \right) = - \rho\E_{\A_{ik}} \A_{im} \pxk \vel_{m}   -\frac{\rho}{\theta_{1} \left(\tau_{1}\right)} \E_{\A_{ik}}\E_{\A_{ik}}. \label{eqn.E2}
\end{equation}
Substituting \eqref{eqn.E2} in \eqref{eqn.pressder1}, we obtain
\begin{gather}
	\frac{\partial \press}{\partial t}  
	+ \bvel\cdot\gra\press
	+ \left( \gamma-1\right) \left( - \rho\E_{\A_{ik}} \A_{im} \pxk \vel_{m}   -\frac{\rho}{\theta_{1} \left(\tau_{1}\right)} \E_{\A_{ik}}\E_{\A_{ik}}
	+ \pt \left( \rho \E_3\right)
	+ \dive \left(\rho\E_3\bvel\right) \right) 
	\notag\\
	+ c^2 \dive \left( \rho\bvel \right) - c^2 \bvel\cdot \gra \rho
	+ \left( \gamma-1\right) \left( \bsigma \cdot \gra \bvel
	+ \bomega \cdot \gra \bvel
	+ \dive \q\right)   = 0. \label{eqn.pressder2}
\end{gather}
Following an analogous procedure, we next substitute the time derivative term in $E_3$. Multiplying \eqref{eqn.GPR_J} by $\rho \E_{J_k}$, summing up all heat flux equations, and adding \eqref{eqn.GPR_rho} multiplied by $\E_3$, lead to
\begin{equation}
	\pt \left( \rho\E_3 \right)  + \dive \left( \rho \E_3 \bvel\right)  = - \rho \E_{J_k}  \J_{m}\pxk \vel_{m} - \rho \E_{J_k}\pxk T - \frac{\rho}{\theta_{2} \left(\tau_{2}\right)} \E_{J_k} \E_{\J_{k}}. \label{eqn.E3}
\end{equation}
Then substituting \eqref{eqn.E3} in \eqref{eqn.pressder2}, we get
\begin{align*}
	\frac{\partial \press}{\partial t}
	&+ \bvel\cdot\gra\press
	+ c^2 \dive \left( \rho\bvel \right)
	- c^2 \bvel\cdot \gra \rho
	+ \left( \gamma-1\right) \left( \bsigma \cdot \gra \bvel
	+ \bomega \cdot \gra \bvel
	+ \dive \q\right)
	\notag\\&
	- \left( \gamma-1\right) \left( \rho\E_{\A_{ik}} \A_{im} \pxk \vel_{m}+ \rho \E_{J_k}  \J_{m}\pxk \vel_{m}
	+ \rho \E_{J_k}\pxk T
	\right) 
	\notag\\& = \left( \gamma-1\right) \left( 
	 \frac{\rho}{\theta_{1} \left(\tau_{1}\right)} \E_{\A_{ik}}\E_{\A_{ik}}
	+ \frac{\rho}{\theta_{2} \left(\tau_{2}\right)} \E_{J_k} \E_{\J_{k}}
	\right) . 
\end{align*}
Replacing $\E_{\J_{k}} = c_{h}^2 \J_k$, $\E_{\A_{ik}}= c_{s}^2\A_{ij} \mathring{G}_{jk}$ in the left hand side of the former equation leads to
\begin{align}
	\frac{\partial \press}{\partial t}
	&+ \bvel\cdot\gra\press
	+ c^2 \dive \left( \rho\bvel \right)
	- c^2 \bvel \cdot\gra \rho
	+ \left( \gamma-1\right) \left( \bsigma \cdot \gra \bvel
	+ \bomega \cdot \gra \bvel
	+ \dive \q\right)
	\notag\\&
	- \left( \gamma-1\right) \left( \rho c_{s}^2 \A_{ij} \mathring{G}_{jk} \A_{im} \pxk \vel_{m}
	+ \rho c_{h}^2 \J_k  \J_{m}\pxk \vel_{m}
	+ \rho c_{h}^2 \J_k \pxk T
	\right) \notag\\&
	= \left( \gamma-1\right) \left( 
	  \frac{\rho}{\theta_{1} \left(\tau_{1}\right)} \E_{\A_{ik}}\E_{\A_{ik}}
	+ \frac{\rho}{\theta_{2} \left(\tau_{2}\right)} \E_{J_k} \E_{\J_{k}}
	\right) . \label{eqn.pressder4}
\end{align}
%
Finally, from \eqref{eqn.def_sigma_omega} and \eqref{eqn.def_q}, we observe that
\begin{gather*}
	\bsigma \cdot \gra \bvel = \sigma_{ik} \pxk \vel_{i} = \A_{ji} \partial_{\A_{jk}} \left(\rho \E\right) \pxk \vel_{i} = \rho c_s^2 \A_{mi}\A_{mj}\mathring{G}_{jk}\pxk \vel_{i},\\
	\bomega\cdot\gra\bvel = \omega_{ik}\pxk\vel_{i} = \rho c_h^2 J_iJ_k\pxk\vel_{i},\\
	\dive \q = \pxk \q_{k} = c_h^2 \pxk \left(\rho J_{k} T\right) =  c_h^2 T \pxk \left(\rho J_{k} \right) + c_h^2 \rho J_{k} \pxk T.
\end{gather*}
Hence
\begin{gather*}
	 \bsigma \cdot \gra \bvel
	 + \bomega \cdot \gra \bvel
	 + \dive \q
	 -  \rho c_{s}^2 \A_{ij} \mathring{G}_{jk} \A_{im} \pxk \vel_{m}	
	 - \rho c_{h}^2 \J_k  \J_{m}\pxk \vel_{m}
	 \notag \\	
	 - \rho c_{h}^2 \J_k \pxk T 
	 = c_h^2 T \pxk \left(\rho J_{k} \right).
\end{gather*}
Substitution in \eqref{eqn.pressder4} yields
\begin{gather*}
	\frac{\partial \press}{\partial t}
	+ \vel_k\pxk\press
	+ c^2 \pxk \left( \rho\vel_k \right)
	- c^2 \vel_k \pxk \rho
	+ c_h^2 \left( \gamma-1\right) T \pxk \left(\rho J_{k} \right)
	\\= \left( \gamma-1\right) \left( 
	\frac{\rho}{\theta_{1} \left(\tau_{1}\right)}  \E_{\A_{ik}}\E_{\A_{ik}}
	+ \frac{\rho}{\theta_{2} \left(\tau_{2}\right)}  \E_{\J_k} \E_{\J_{k}}
	\right) . \label{eqn.pressder6}
\end{gather*}
Or, equivalently,
\begin{gather*}
	\frac{\partial \press}{\partial t}
	+ \bvel\cdot\gra\press
	+ c^2 \dive \left( \rho\bvel \right)
	- c^2 \bvel \cdot\gra \rho
	+ c_h^2 \left( \gamma-1\right) T \dive\left(\rho \bJ\right)
	\notag\\ = \left( \gamma-1\right) \left( 
	\frac{\rho}{\theta_{1} \left(\tau_{1}\right)} \bE_{\bA}\cdot\bE_{\bA}
	+ \frac{\rho}{\theta_{2} \left(\tau_{2}\right)} \bE_{\bJ}\cdot \bE_{\bJ}
	\right) . \label{eqn.pressder7}
\end{gather*}

Therefore, the weakly compressible GPR model, where the energy conservation equation has been replaced by a non-conservative pressure equation, reads: 
\begin{subequations}\label{eqn.WGPR}
	\begin{eqnarray}
		\pt \rho + \pxk \left(\rho\vel_{k}\right) &=& 0, \label{eqn.WGPR_rho}\\
		\pt \left(\rho\vel_i\right) + \pxk \left(\rho \vel_{i} \vel_{k} \right) + \pxi \press +\pxk\left( \sigma_{ik}+\omega_{ik} \right) &=& \rho g_{i},\\
		\pt \A_{ik} + \pxk \left(\vel_{m}\A_{im}\right) + \vel_{j} \pxj \A_{ik} -\vel_{j}\pxk \A_{ij} &=& -\frac{1}{\theta_{1} \left(\tau_{1}\right)} \E_{\A_{ik}}, \\
		\pt \J_{k} + \pxk \left(\J_{m}\vel_{m} \right) +\pxk T + \vel_{j}\left(\pxj J_{k} - \pxk J_{j}\right)&=& -\frac{1}{\theta_{2} \left(\tau_{2}\right)} \E_{\J_{k}}, \\
		\frac{\partial \press}{\partial t}
		\!+\! \vel_k\pxk\press
		\!+\! c^2 \! \pxk \left( \rho\vel_k \right)
		\!-\! c^2 \vel_k \pxk \rho
		\!+\! c_h^2 \left(\! \gamma-1\right) T\! \pxk \left(\rho J_{k} \right) 
		&=&
		\frac{\left( \!\gamma-1\!\right)\!\rho}{\theta_{1} \left(\tau_{1}\right)}  \E_{\A_{ik}}\E_{\A_{ik}}
		\!+\! \frac{\left(\! \gamma-1\!\right)\!\rho}{\theta_{2} \left(\tau_{2}\right)}  \E_{\J_k} \E_{\J_{k}}.\qquad\quad
	\end{eqnarray}
\end{subequations}

\section{Numerical discretisation} \label{sec:numdisc}
The discretisation of the former GPR systems will be performed in the framework of the hybrid finite volume/finite element approach introduced in \cite{BFSV14,BFTVC17,Hybrid1,BRMVCD21,HybridMPI,HybridNNT} for incompressible, weakly compressible and all Mach number flows and the shallow water equations.
In particular, we are interested in the low Mach number limit, so we first address the incompressible GPR model \eqref{eqn.incGPR}, and we then also propose a hybrid FV/FE approach for the discretisation of a pressure-based formulation of \eqref{eqn.WGPR} able to address weakly compressible flows.

\subsection{Semi-discretisation in time of the incompressible GPR model}
A semi-discretisation in time of system \eqref{eqn.incGPR}, yields
\begin{subequations}\label{eqn.incGPR_semidiscrete}
	\begin{eqnarray}
		\pxk \left(\mom_{k}^{n+1}\right) &=& 0, \label{eqn.incGPR_semidiscrete_rho}\\
		\frac{1}{\Dt} \left(\mom_i^{n+1}-\mom_i^{n}\right) + \pxk \left(\mom_{i}^{n} \vel_{k}^{n} \right) + \pxi \press^{n+1} +\pxk\sigma_{ik}^{n} &=& \rho^{n} g_{i}, \label{eqn.incGPR_semidiscrete_mom}\\
		\frac{1}{\Dt} \left( \A_{ik}^{n+1}-\A_{ik}^{n} \right) + \pxk \left(\vel_{m}^{n}\A_{im}^{n}\right) + \vel_{j}^{n} \pxj \A_{ik}^{n} -\vel_{j}^{n}\pxk \A_{ij}^{n} &=& -\frac{1}{\theta_{1}^{n+1} \left(\tau_{1}\right)} \E_{\A_{ik}}^{n+1}.\qquad \qquad \label{eqn.incGPR_semidiscrete_A}
	\end{eqnarray}
\end{subequations}
Following classical projection methods, \cite{PABCHL95,Guer06,BFTVC17}, we split the momentum equation into two parts and gather the equations into a transport-diffusion and a pressure subsystem:
\paragraph{Transport-diffusion subsystem}
\begin{subequations}\label{eqn.ingpr_disc_transp}
	\begin{eqnarray}
	\tmom_i & = & \mom_i^{n}- \Dt\left(  \pxk \left(\mom_{i}^{n} \vel_{k}^{n} \right) + \pxi \press^{n} +\pxk\sigma_{ik}^{n} - \rho^{n} g_{i}\right) , \\
	 \A_{ik}^{n+1} & = &  \A_{ik}^{n} -\Dt \left( \pxk \left(\vel_{m}^{n}\A_{im}^{n}\right) + \vel_{j}^{n} \pxj \A_{ik}^{n} -\vel_{j}^{n}\pxk \A_{ij}^{n} + \frac{1}{\theta_{1}^{n+1} \left(\tau_{1}\right)} \E_{\A_{ik}}^{n+1}\right).\qquad\qquad\label{eqn.ingpr_disc_transp.A}
\end{eqnarray}
\end{subequations}

\paragraph{Pressure subsystem}
\begin{subequations}\label{eqn.incGPR_disc_press}
	\begin{eqnarray}
		\pxk \left(\mom_{k}^{n+1}\right) &=& 0, \label{eqn.incGPR_disc_press1}\\
		\mom_i^{n+1} &=& \tmom_i - \Dt \pxi \delta \press^{n+1}, \qquad \delta \press^{n+1}= \press^{n+1}-\press^{n} . \label{eqn.incGPR_disc_press2}
	\end{eqnarray}
\end{subequations}

Taking into account the nature of these systems, we will employ an explicit finite volume approach for the spatial discretisation of \eqref{eqn.ingpr_disc_transp} while \eqref{eqn.incGPR_disc_press1}-\eqref{eqn.incGPR_disc_press2} will be solved implicitly using continuous finite elements.

\subsection{Semi-discretisation in time of the weakly compressible GPR model}
Similarly, for the weakly compressible GPR model \eqref{eqn.WGPR}, we apply an splitting procedure following \cite{TV12} and, introducing the semi-discretisation in time, we get
\begin{subequations}\label{eqn.WGPR_semidiscrete}
	\begin{align}
		&\frac{1}{\Dt} \left( \rho^{n+1}-\rho^{n} \right) 
		+ \pxk \left(\mom_{k}^{n+1}\right) 
		= 0, \label{eqn.WGPR_semidiscrete_rho}\\
		&\frac{1}{\Dt} \left(\mom_i^{\ast}-\mom_i^{n}\right) 
		+ \pxk \left(\mom_{i}^{n} \vel_{k}^{n} \right) 
		+\pxk\sigma_{ik}^{n} +\pxk\omega_{ik}^{n} 
		= \rho^{n} g_{i}, \label{eqn.WGPR_semidiscrete_mom}\\
		&\frac{1}{\Delta t}\left( \mom^{n+1}_i-\tmom_i\right) 
		+\pxi \press^{n+1} \label{eqn.WGPR_semidiscrete_mompress}
		= 0,\\
		&\frac{1}{\Delta t}\left( \A_{ik}^{n+1}\!-\!\A_{ik}^{n}\right) 
		+ \pxk \left(\vel_{m}^{n}\A_{im}^{n}\right) 
		+ \vel_{j}^{n} \pxj \A_{ik}^{n} 
		-\vel_{j}^{n}\pxk \A_{ij}^{n} 
		= -\frac{1}{\theta_{1}^{n+1} \left(\tau_{1}\right)} \E_{\A_{ik}}^{n+1}, \label{eqn.WGPR_semidiscrete_A}\\
		&\frac{1}{\Delta t}\!\left( \J_{k}^{n+1}\!-\!\J_k^{n}\right) 
		\!+\! \pxk \left(\J_{m}^{n}\vel_{m}^{n} \right) +\pxk T^{n} 
		\!+\! \vel_{j}^{n}\left(\!\pxj J_{k}^{n} \!-\! \pxk J_{j}^{n}\!\right) 
		=  -\frac{1}{\theta_{2}^{\star} \left(\tau_{2}\right)} \E_{\J_{k}}^{n+1}, \label{eqn.WGPR_semidiscrete_J}\\
		&\frac{1}{\Delta t}\left( \widetilde{\press}_{\press}-\press^{n}\right) 
		+ \vel^{n}_{k} \pxk \press^{n}
		+  c_h^2 \left(\gamma-1\right)  T^{n} \pxk \left( \rho^{n} J_k^{n}\right) 
		= \mathcal{S}^{\press}\left(\Q^{n+1}\right),  \label{eq:disc_pres_tilde} \\
		&\frac{1}{\Delta t} \widetilde{\press}_{\rho} 
		- c^{2}\vel^{n}_{k} \pxk \rho^{n}  
		= 0, \label{eq:disc_pres_tilderho} \\
		&\frac{1}{\Delta t}\left( \press^{n+1}- \widetilde{\press}\,\right) 
		+ c^{2} \pxk \left( \mom^{n+1}_k\right)   
		= 0. \label{eq:WGPR_disc_pres2}
	\end{align}
\end{subequations}
with
\begin{equation*}
	\mathcal{S}^{\press}\left(\Q^{n+1}\right)= \frac{\rho^{n+1}\left(\! \gamma\!-\! 1\!\right)}{\theta_{1}^{n+1} \left(\tau_{1}\right)} \E_{\A_{ik}}^{n+1}\!\E_{\A_{ik}}^{n+1}
	\!+ \!\frac{\rho^{n+1}\left(\! \gamma\! -\! 1\!\right)}{\theta_{2}^{\star} \left(\tau_{2}\right)} \E_{J_k}^{n+1}\! \E_{\J_{k}}^{n+1}
\end{equation*}	
and $\Q=\left(\rho,\bvel,\bA,\bJ,\press\right)$.
As for the incompressible GPR model, gathering \eqref{eqn.WGPR_semidiscrete_mompress} and \eqref{eq:WGPR_disc_pres2}, we get a pressure system of the form
\begin{subequations}\label{eqn.pressuresystem_WGPR}
\begin{eqnarray}
	\frac{1}{\Delta t}\left( \press^{n+1}- \widetilde{\press}\,\right) + c^{2} \pxk \left( \mom^{n+1}_k\right)& = &0, \label{eqn.WGPR_presuresystem_WGPR1}	\\
	\mom^{n+1}_i &=& \tmom_i -\Dt \pxi \press^{n+1}, \label{eqn.WGPR_semidiscrete_mompress2}	
\end{eqnarray}
\end{subequations}
corresponding to a Poisson-type problem, and a set of transport equations containing the conservative fluxes and non-conservative products, \eqref{eqn.WGPR_semidiscrete_rho}-\eqref{eqn.WGPR_semidiscrete_mom}, \eqref{eqn.WGPR_semidiscrete_A}-\eqref{eq:disc_pres_tilderho}.

\subsection{Overall algorithm}
Attending to the nature of the different equations involved in \eqref{eqn.ingpr_disc_transp}-\eqref{eqn.incGPR_disc_press} and \eqref{eqn.WGPR_semidiscrete},
the proposed hybrid FV/FE methodology is divided into the following stages:

\begin{itemize}
	\item Transport stage. The equations containing the convective and non-conservative terms, i.e., system \eqref{eqn.ingpr_disc_transp} for the incompressible GPR model and equations \eqref{eqn.WGPR_semidiscrete_rho}, \eqref{eqn.WGPR_semidiscrete_mom}, \eqref{eqn.WGPR_semidiscrete_A}, \eqref{eqn.WGPR_semidiscrete_J}, and \eqref{eq:disc_pres_tilde} for the weakly compressible GPR model, are discretised explicitly using a finite volume scheme. 
	Let us note that both systems can be recast in the general form
	\begin{equation}
		\partial_{t} \Q + \dive \Flux\left(\Q\right) + \NCP{\Q} = \source\left(\Q\right),\label{eqn.generaltransportncp}
	\end{equation}
	with $\Q$ the vector of conservative variables, $\Flux\left(\Q\right)$ the flux term, $\NCP{\Q}$ the non-conservative products, and $\source\left(\Q\right)$ the source terms.

	\item Interpolation stage. This stage is only necessary for the weakly GPR model where the contribution of the non-conservative product on the density derivative appearing in the pressure equation, \eqref{eq:disc_pres_tilderho}, is approximated by making use of an explicit finite volume approach. Moreover, the pressure intermediate value $\widetilde{\press}_{\press}$ obtained in the dual cells during the convective stage is interpolated to the primal grid.
	
	\item Pressure stage. The pressure subsystems \eqref{eqn.incGPR_disc_press} or \eqref{eqn.pressuresystem_WGPR} are solved using continuous finite elements.
	
	\item Correction stage. In the case of incompressible flows, the use of the pressure gradient at the previous time step is not sufficient to ensure the divergence-free condition of the velocity field, so we must correct $\tbmom$ with the gradient of the pressure variation $\gra \delta \press^{n+1}$. On the other hand, for the weakly compressible GPR model the pressure gradient has been completely neglected in the convective stage, so we need to incorporate its contribution once the new pressure is available, $\gra \press^{n+1}$. Therefore, the intermediate momentum $\tbmom$ is corrected using \eqref{eqn.incGPR_disc_press2} and \eqref{eqn.WGPR_semidiscrete_mompress2} for the incompressible and the weakly compressible GPR models, respectively. 
\end{itemize}

In what follows, we introduce the spatial discretisation and provide a detailed description of each algorithm stage.

\subsection{Spatial discretisation. Unstructured staggered grids}
\begin{figure}[h]
	\begin{center}
		\includegraphics[width=0.35\linewidth]{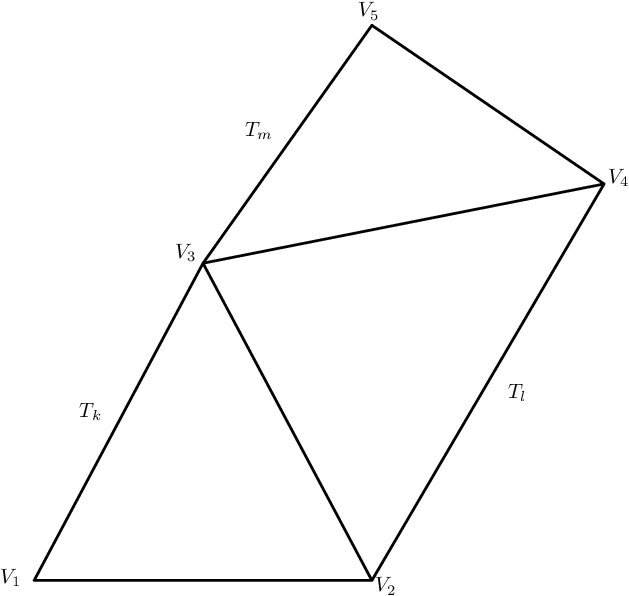}
		\includegraphics[width=0.35\linewidth]{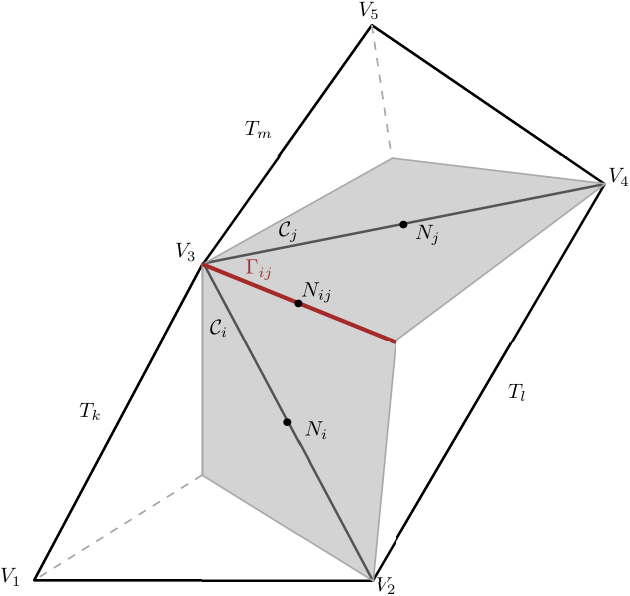}
	\end{center}
	\caption{Sketch of the face-based unstructured grids in 2D. Left: $\Pcell{k}$, $\Pcell{l}$, $\Pcell{m}$ are the triangles of the primal grid and $V_{j},\, j=1,\dots, 5$ are their vertices. Right: Interior (grey elements) and boundary (white elements) elements of the dual mesh. The boundary between the two interior dual cells $\cell{\ci}$ and  $\cell{\cj}$, $\Gamma_{\ci\cj}$, is highlighted in red.}
	\label{fig:mesh}
\end{figure}

To discretise the computational domain we employ the so-called face-based or diamond-shaped staggered grids, \cite{BDDV98,VC99,BFSV14,TD16}. We denote $\boldsymbol{\Upsilon}=\left\lbrace \Pcell{\ck}, \, \ck =1, \dots, M \right\rbrace$  the tessellation corresponding to the primal grid composed of $M$ triangular elements $\Pcell{\ck}$. Then, each triangle $\Pcell{\ck}$ is divided into three subtriangles having as base one of the boundary edges of $\Pcell{\ck}$ and opposite vertex the barycentre of $\Pcell{\ck}$, denoted by $\x_{\ck}$. Merging the two subtriangles related to a boundary edge, we get an interior dual cell $\cell{\ci}$. Meanwhile, for the element edges located at a boundary of the domain, the related dual element is simply taken to be the corresponding subtriangle inside the domain. A sketch of the dual mesh construction in 2D is depicted in Figure~\ref{fig:mesh}. For a detailed description of the unstructured faced-based staggered grids in 3D, we may refer to \cite{Hybrid1}.

The use of these kind of staggered grids is two-folded motivated. On the one hand, having a primal grid made of triangles/tetrahedra eases the tessellation of complex domains if compared with Cartesian grids made of quadrilaterals/hexahedra. On the other hand, the combination of two staggered grids, the primal one to be used within the pressure system discretisation and the dual one employed during the transport stage, avoids stability issues as the well-known checkerboard phenomenon that often appears when collocated grids are employed, \cite{Tyl16}. Further, the use of this grid arrangement becomes useful for the design of a second order finite volume scheme with a very small stencil, as it has already been shown in \cite{BFTVC17,Hybrid1}. 

To complete the description of the spatial domain discretisation in 2D, we still need to introduce some notation to be employed in the description of the algorithm stages. Given a dual cell $\cell{\ci}$, we denote $\left|\cell{\ci}\right|$ its area and $\Gamma_{\ci}\equiv\partial \cell{\ci}$ its boundary. Taking into account the shape of the dual interior elements, $\Gamma_{\ci}$ can be decomposed into four straight edges labelled as $\Gamma_{\ci \cj}$, with $\cell{\ci}$ and $\cell{\cj}$ the two dual cells sharing that edge. Similarly, for a boundary cell we have three edges, two of them of the interior type $\Gamma_{\ci \cj}$ and a boundary edge $\Gamma_{\ci \Gamma}$. Moreover, $\un_{\ci\cj}$ represents the unitary external normal of $\Gamma_{\ci\cj}$, and $\nn_{\ci\cj}$ the length weighted normal, i.e., $\nn_{\ci\cj}=\un_{\ci\cj}\left|\Gamma_{\ci\cj}\right|=\un_{\ci\cj}\left\|\nn_{\ci\cj}\right\|$ with $\left|\Gamma_{\ci\cj}\right|$ the length of edge $\Gamma_{\ci\cj}$. Finally, $V_{\ck\cl}$, $\cl\in\left\lbrace 1,\dots,3\right\rbrace$, refer to the three vertex of a primal element $\Pcell{\ck}$ and $\un_{\ck}$ denotes its outward-pointing unit normal.

\subsection{Transport stage. Finite volume method in the dual grid}
System \eqref{eqn.generaltransportncp} is discretised by employing an explicit finite volume method. Accordingly, we integrate \eqref{eqn.generaltransportncp} on each control volume $\cell{\ci}$ and apply Gauss theorem to transform the integral of the flux term into the integral of the normal flux along the cell boundary yielding to
\begin{equation}
	\tQ_{\ci} = \Q^{n}_{\ci} 
	- \frac{\Dt}{\vcell{\ci}} \left( \int\limits_{\Gamma_{\ci}} \Flux\left(\Q^{n}\right) \cdot \un_{\ci} \dS 
	+  \int\limits_{\cell{\ci}} \NCP{\Q^{n}}\dV 
	- \int\limits_{\cell{\ci}}  \source\left(\Q\right)\dV\right) ,
	\label{eqn.generaltransportncp_disc}
\end{equation}
with $\tQ=\left(\rho^{n+1}, \tbmom, \bA^{n+1}\right)^{T}$ for the incompressible GPR model and $\tQ=\left(\rho^{n+1}, \tbmom, \bA^{n+1}, \bJ^{n+1}, \tpress_{\press}\right)^{T}$ for the weakly compressible GPR model. In what follows, to describe the FV scheme, we focus on the weakly compressible GPR case since the incompressible one can be seen just as a subcase of it.

\subsubsection{Explicit treatment of the convective terms}
The integral of the flux term is decomposed onto the sum of the contributions of the normal flux along each cell boundary as
\begin{equation}
	\int\limits_{\Gamma_{\ci}} \Flux\left(\Q^{n}\right) \cdot \un_{\ci} \dS 
	= \sum\limits_{N_{\cj}\in\mathcal{K}_{\ci}} \left|\Gamma_{\ci \cj}\right| \Flux^{NF}\left(\eQ^{n}_{\ci},\eQ^{n}_{\cj},\un_{\ci\cj}\right)  ,
	\label{eqn.flux}
\end{equation}
with $\mathcal{K}_{\ci}$ the set of neighbours of $\cell{\ci}$ and $\Flux^{NF}$ a numerical flux function. In particular, we employ the Rusanov numerical flux, \cite{Rus62},
\begin{equation}
	 \Flux^{R}\left(\eQ^{n}_{\ci},\eQ^{n}_{\cj},\un_{\ci\cj}\right) = \halb\left(\Flux\left(\eQ_{\ci}^{n}\right)+\Flux\left(\eQ_{\cj}^{n}\right)\right)
	 - \alpha_{\ci\cj}^{n} \left(\eQ_{j}^{n}-\eQ_{i}^{n}\right),
	\label{eqn.fluxRS}
\end{equation}
with the maximum signal speed on the edge
\begin{gather}
\alpha_{\ci\cj}^{n} = 
\max\left\lbrace
\left|\overline{\bvel}_{\ci}^{n}\cdot \un_{\ci\cj} \pm c_s \right|, \left|\frac{3}{2}\overline{\bvel}_{\ci}^{n}\cdot \un_{\ci\cj} \pm c_{\ci} \right|, \left| \overline{\bvel}_{\cj}^{n}\cdot \un_{\ci\cj}\pm c_s\right|,\left| \frac{3}{2}\overline{\bvel}_{\cj}^{n}\cdot \un_{\ci\cj}\pm c_{\cj}\right|
\right\rbrace,
\notag \\
c_{\ci} = \sqrt{\frac{4}{3} c_s^2 +\frac{1}{4}\left| \bvel_{\ci}\right|^2  } ,
\end{gather}
for the incompressible GPR model and
\begin{equation}
	\alpha_{\ci\cj}^{n} = 
	\max\left\lbrace
	\left|\overline{\bvel}_{\ci}^{n}\cdot \un_{\ci\cj} \pm c_{\ci}^{n} \right|, \left| \overline{\bvel}_{\cj}^{n}\cdot \un_{\ci\cj}\pm c_{\cj}^{n}\right|
	\right\rbrace,
	\qquad
	c_{\ci}^{n} = \sqrt{\frac{4}{3} c_s^2 + \frac{2 c_{h}^{2}  \overline{T}_{\ci}^{n}}{\left( \overline{\rho}^{n}_{\ci}\right)^2 c_v}  } ,
\end{equation}
for the weakly compressible GPR model. Besides, if we sought a first order scheme, $\eQ_{\ci}$ and $\eQ_{\cj}$ are simply taken as the values of $\Q$ at the two dual cells related to the dual edge, $\eQ_{\ci}:=\Q_{\ci}$, $\eQ_{\cj}:=\Q_{\cj}$. On the other hand, to attain second order, $\eQ_{\ci}$ and $\eQ_{\cj}$ must correspond to the half in time evolved boundary extrapolated values. More precisely, we consider the local ADER approach proposed in \cite{BFTVC17} and perform the following steps:
\begin{enumerate}
	\item Piece-wise polynomial reconstruction. Given a variable $Q$, we build the left and right reconstruction polynomials related to edge $\Gamma_{\ci\cj}$ as
	\begin{equation}
		P^{L}_{\ci\cj}(\x) = Q_{\ci} + \left(\x-\x_{\ci} \right) \gra Q_{\ci\cj}^{L},\qquad
		P^{R}_{\ci\cj}(\x) = Q_{\cj} + \left(\x-\x_{\cj} \right) \gra Q_{\ci\cj}^{R}.
	\end{equation}
	The slopes $\gra Q_{\ci\cj}^{L}$ and $\gra Q_{\ci\cj}^{R}$ are computed using an ENO interpolation method so that the final scheme is nonlinear and therefore circumvents Godunov's theorem. Accordingly, denoting $\Pcell{\ci\cj}$ the primal element containing the face $\Gamma_{\ci\cj}$ and $\Pcell{\ci\cj}^{L}$ and  $\Pcell{\ci\cj}^{R}$ the two neighbour primal elements containing one halve of the dual cells $\cell{\ci}$ and  $\cell{\cj}$, the slopes are computed as
	\begin{gather}
		\gra Q_{\ci\cj}^{L} = \left\lbrace \begin{array}{lr}
			\gra Q_{\mid\Pcell{\ci\cj}^{L}}  &  \textrm{if}\, \left|\gra Q_{\mid\Pcell{\ci\cj}^{L}}\cdot \left(\x_{\ci\cj}-\x_{\ci} \right) \right| \leq \left|\gra Q_{\mid\Pcell{\ci\cj }}\cdot \left(\x_{\ci\cj}-\x_{\ci} \right) \right|,\\
			\gra Q_{\mid\Pcell{\ci\cj  }}  &  \textrm{otherwise} ;
		\end{array}\right.\\
		\gra Q_{\ci\cj}^{R} = \left\lbrace \begin{array}{lr}
			\gra Q_{\mid \Pcell{\ci\cj}^{R}}  &  \textrm{if}\, \left|\gra Q_{\mid\Pcell{\ci\cj}^{R}}\cdot \left(\x_{\ci\cj}-\x_{\ci} \right) \right| \leq \left|\gra Q_{\mid\Pcell{\ci\cj }}\cdot \left(\x_{\ci\cj}-\x_{\ci} \right) \right|,\\
			\gra Q_{\mid\Pcell{\ci\cj  }}  &  \textrm{otherwise} .
		\end{array}\right.
	\end{gather}
	The gradients $\gra Q_{\mid\Pcell{\ci\cj}}$, $\gra Q_{\mid\Pcell{\ci\cj}^{L}}$ and $\gra Q_{\mid\Pcell{\ci\cj}^{R}}$ are computed in the primal cells using Crouzeix-Raviart finite elements which have as nodes the barycentres of the faces that are identified with the nodes of the dual cells.
	
	\item Computation of boundary extrapolated data. The polynomials are evaluated in the barycentre of the dual edge, $\x_{\ci\cj}$, obtaining the boundary extrapolated values $ Q_{\ci}^{L}$ and $ Q_{\cj}^{R}$.
	\item Half in time evolution. A midpoint rule, combined with the Cauchy-Kovalevskaya procedure to transform the time derivatives into spatial derivatives using the governing equations, provides the approximation of the conservative variables at time $t^{n}+\halb\Dt$. Further details on this methodology and on the original ADER approach can be found, e.g. in \cite{BRMVCD21,Toro,TMN01}. Moreover, for recent advances in ADER methods,  including the ADER-DG approach which avoids the Cauchy-Kovalevskaya procedure by introducing a local space-time predictor, we refer to \cite{DBTM08,BCDGP20}.
\end{enumerate} 
 In the numerical results, Section~\ref{sec:numericalresults}, as an alternative to the ENO-based reconstruction introduced above, we also consider the use of the min-mod limiter of Roe, \cite{Roe85}, and the Barth and Jespersen limiter, \cite{BJ89}.
 
\subsubsection{Pressure gradient in the incompressible GPR model}
Let us note that the incompressible GPR model equation \eqref{eqn.incGPR_semidiscrete_mom} includes a term on the pressure gradient at the previous time step. Contrary to what is done in most Godunov-type methods, where this kind of term is included within the flux, \cite{DPRZ16,Busto2022HTCGPR}, we compute it as if it was a source term, since it does not depend on the pressure at the new time step and thus it does not need to be included in the convective terms nor in the CFL time step restriction of the explicit subsystem. 
Therefore, to approximate its contribution, we interpolate the pressure at the previous time step, which has been computed in the primal vertex, into the dual edges, $\press_{\ci\cj}^{n}$, by simply taking the average between the two vertexes of each edge. In case one of the vertex corresponds to the barycentre of the primal element, its value is first obtained by averaging the pressure at the three vertex of the primal element. Finally, we compute
\begin{equation}
	\int\limits_{\cell{\ci}} \gra \press^{n} \dV = \sum\limits_{N_{\cj}\in\mathcal{K}_{\ci}} \press_{\ci\cj}^{n}\nn_{\ci\cj}.
\end{equation}

\subsubsection{Path conservative discretisation of the non-conservative products}\label{sec.NCP}
The non-conservative products, $\NCP{\Q}$, are discretised employing a path conservative scheme based on the straight line segment path, \cite{Par06,CFFP09,GCD18}. Accordingly, we approximate
\begin{equation}
	\int\limits_{\cell{\ci}} \NCP{\Q^{n}}\dV  = \int\limits_{\Gamma{\ci}} \mathcal{D}\left( \eQ^{n}\right)\cdot\un_{\ci} \dS  + \int\limits_{\cell{\ci}\setminus \Gamma_{\ci}} \NCP{\eQ^{n}}\dV.
	\label{eqn.NCP}
\end{equation}

In Equation \eqref{eqn.NCP}, the first term considers the jumps of the discrete solution across the cell boundaries for which we employ the boundary extrapolated values related to the face: 
\begin{equation}
	\int\limits_{\Gamma{\ci}} \mathcal{D}\left( \eQ^{n}\right)\cdot\un_{\ci} \dS 
	= \halb \sum\limits_{N_{\cj}\in\mathcal{K}_{\ci}} \BNCP{\eQ_{\ci\cj}} \cdot \nn_{\ci\cj} \left(\eQ_{\cj}^{n}-\eQ_{\ci}^{n}\right),\qquad 
	\eQ_{\ci\cj} = \eQ_{\ci}+\eQ_{\cj}.
	\label{eqn.NCP_jumps}
\end{equation}

Meanwhile, the second term in Equation \eqref{eqn.NCP} corresponds to the smooth contribution of the non-conservative product within the cell, which must be taken into account to get high-order accurate schemes. To compute it, we again employ the dual grid structure and approximate the needed gradients using a Galerkin approach in the primal grid. Then, the gradient for the non-conservative product contribution is computed as a weighted average of the contribution from the two primal subtriangles composing it, i.e.,
\begin{equation}
	\int\limits_{\cell{\ci}\setminus \Gamma_{\ci}} \NCP{\eQ^{n}}\dV = \vcell{\ci} \BNCP{\eQ_{\ci}^{n}} \left(\frac{\vcell{\ci_{1}}}{\vcell{\ci}} \gra \eQ_{\ci_{1}}^{n} +\frac{\vcell{\ci_{2}}}{\vcell{\ci}} \gra \eQ_{\ci_{2}}^{n}\right), \label{eqn.NCP_smooth}
\end{equation}
with $\ci_{1}$ and $\ci_{2}$ the two halves of cell $\cell{\ci}$, $\cell{\ci}=\cell{\ci_{1}}\cup \cell{\ci_{2}}$, $\cell{\ci_{1}}\subset \Pcell{\ci_{1}}$, $\cell{\ci_{2}}\subset \Pcell{\ci_{2}}$, $\Pcell{\ci_{1}},\Pcell{\ci_{2}}\in\boldsymbol{\Upsilon}$, and $\gra \eQ_{\ci_{1}}$, $\gra \eQ_{\ci_{2}}$ the gradients of $\eQ$ computed in $\Pcell{\ci_{1}}$ and $\Pcell{\ci_{2}}$, respectively.

\subsubsection{Source term in the momentum equations}
The source term of the momentum equations is integrated on each dual cell employing the density at the previous time step as
\begin{equation}
	\int\limits_{\cell{\ci}} \overline{\rho}^{n} \g \dV = \vcell{\ci} \overline{\rho}_{\ci}^{n}\g.
\end{equation}

\subsubsection{Implicit discretisation of the algebraic source terms for the distortion and heat conduction fields}
The algebraic source terms related to the relaxation times may become very stiff in the fluid limit of the equations. Consequently, the needed time step to treat them explicitly may become very restrictive. To avoid this issue, an implicit discretisation of those algebraic source terms can be performed. We assume that convective and non-conservative terms in \eqref{eqn.WGPR_semidiscrete} have already been computed. Then, we get the following system for the algebraic source terms:
\begin{subequations}\label{eqn.sourcesystem_sd}
\begin{eqnarray}
	\frac{1}{\Delta t}\left( \A_{ik}^{n+1}\!-\!\A_{ik}^{\ast}\right) 
	&=& -\frac{1}{\theta_{1}^{n+1} \left(\tau_{1}\right)} \E_{\A_{ik}}^{n+1}, \\
	\frac{1}{\Delta t}\!\left( \J_{k}^{n+1}\!-\!\J_k^{\ast}\right) 
	&=&  -\frac{1}{\theta_{2}^{\star} \left(\tau_{2}\right)} \E_{\J_{k}}^{n+1},
\end{eqnarray}
\end{subequations}
with
\begin{eqnarray}
	\A_{ik}^{\ast}=\A_{ik}^{n} 
	-	 \pxk \left(\vel_{m}^{n}\A_{im}^{n}\right) 
	- \vel_{j}^{n} \pxj \A_{ik}^{n} 
	+\vel_{j}^{n}\pxk \A_{ij}^{n}, \\
	 \J_{k}^{\ast} = \J_k^{n}
	- \pxk \left(\J_{m}^{n}\vel_{m}^{n} \right) -\pxk T^{n} 
	- \vel_{j}^{n}\left(\pxj J_{k}^{n} + \pxk J_{j}^{n}\right),
\end{eqnarray}
that have already been computed explicitly using the explicit finite volume scheme. 
We note that system \eqref{eqn.sourcesystem_sd} can also be seen as a second splitting of the original model \eqref{eqn.WGPR} with corresponding continuous source term subsystem
\begin{subequations}\label{eqn.sourcesystem}
\begin{eqnarray}
	\partial_{t} \A_{ik} = -\frac{1}{\theta_{1} \left(\tau_{1}\right)} \E_{\A_{ik}}, \\
	 \partial_{t}  \J_{k}=  -\frac{1}{\theta_{2} \left(\tau_{2}\right)} \E_{\J_{k}}.
\end{eqnarray}
\end{subequations}

Denoting $\bY = (\bA,\bJ)^{T}=(\A_{1},\dots, \A_{9},J_{1},\dots,J_{3})^{T}$ the vector of unknowns, the system \eqref{eqn.sourcesystem} can be recast into an ordinary differential equation system of the form
\begin{equation}
	\bY ' (t) = \mathbf{g}(t,\bY(t)),
\end{equation}
that can be solved using classical Runge-Kutta methods as the implicit Euler scheme or the DIRK scheme of Pareschi and Russo, \cite{PR05}. This methodology requires the computation of the root of $\mathbf{G}(\bY(t)) = \bY ' (t) - \mathbf{g}(t,\bY(t))$ which is performed using an inexact Newton algorithm, \cite{DES82}.

Focusing on the weakly compressible GPR system, once the distortion and heat flux fields are obtained at the new time step, they are used to approximate the source term contribution of the pressure equation \eqref{eq:disc_pres_tilde} at each cell $\cell{\ci}$ as
\begin{gather}
	\int\limits_{\cell{\ci}} \left[  \frac{\rho^{n+1}\left(\! \gamma\!-\! 1\!\right)}{\theta_{1}^{n+1} \left(\tau_{1}\right)} \E_{\A_{ik}}^{n+1}\!\E_{\A_{ik}}^{n+1}
	\!+ \!\frac{\rho^{n+1}\left(\! \gamma\! -\! 1\!\right)}{\theta_{2}^{\star} \left(\tau_{2}\right)} \E_{J_k}^{n+1}\! \E_{\J_{k}}^{n+1}\right]  \dV = \notag\\
	\left|\cell{\ci} \right|\left[ \frac{ 3\rho_0^2 c_{s}^2  (\gamma-1) }{\tau_{1} \rho^{n+1}_{\ci} \left| \bA_{\ci}^{n+1} \right|^{\frac{11}{3}} }
	\left( \bA_{\ci}^{n+1} \devG_{\ci}^{n+1}\right) \cdot \left( \bA_{\ci}^{n+1} \devG_{\ci}^{n+1}\right) 
	+ \frac{ c_{h}^{2}  (\gamma-1)  \rho_{\ci}^{n+1}  T_{\ci}^{n}}{\tau_{2}\rho_{0}  T_{0}} \bJ_{\ci}\cdot\bJ_{\ci} \right].
\end{gather}

\subsection{Interpolation stage}
In the weakly compressible GPR model the pressure is computed in two steps. First, an intermediate value gathering the contributions of the convective terms and non-conservative products is obtained from solving \eqref{eq:disc_pres_tilde} and \eqref{eq:disc_pres_tilderho}. The obtained value is then employed within system \eqref{eqn.pressuresystem_WGPR} to calculate the pressure at the new time step. Hence, before the projection stage, the intermediate pressure  $\widetilde{\press}_{\press}$ is interpolated from the dual cells, $\cell{\ci}$, to the primal elements, $\Pcell{\ck}$, as
\begin{equation}
	\widetilde{\press}_{\press\, \ck} = \sum\limits_{\ci \in \mathcal{K}_{\ck}} \frac{\left| \Pcell{\ck\ci}\right|}{\left| \Pcell{\ck}\right|} \widetilde{\press}_{\press\, \ci}, \qquad \Pcell{\ck\ci}=\Pcell{\ck} \cap \cell{\ci} ,\label{eqn.interpolation_dualcellprimalcell}
\end{equation}
with $\mathcal{K}_{\ck}$ the set of dual cell index identifying the dual elements generated from the primal faces of element $\Pcell{\ck}$.
Next, the intermediate pressure $\widetilde{\press}$ is computed as
\begin{equation}
	\widetilde{\press}_{\ck} = \widetilde{\press}_{\press\, \ck} + \widetilde{\press}_{\rho\, \ck},
\end{equation}
where the contribution of the non-conservative product on the density, $\widetilde{\press}_{\rho}$, to the pressure equation \eqref{eq:disc_pres_tilderho}, is computed using a finite volume approach in the primal grid. More precisely, we have 
\begin{gather}
	\widetilde{\press}_{\rho\, \ck} = -\frac{\Dt}{\left|\Pcell{\ck}\right|} \int\limits_{\Pcell{\ck}} \left( c^{n} \right)^{2} \bvel^{n}\cdot \gra \rho^{n} \dV
	=\frac{\Dt}{\left|\Pcell{\ck}\right|}\left( c^{n}_{\ck} \right)^{2} \bvel^{n}_{\ck} \cdot \sum\limits_{\ci\in\mathcal{K}_{\ck}}  \int\limits_{\Gamma_{\ck\ci}} \rho^{n}_{\ci} \un_{\ck\ci} \dV \notag\\
	=\frac{\Dt}{\left|\Pcell{\ck}\right|}\left( c^{n}_{\ck} \right)^{2} \bvel^{n}_{\ck} \cdot \sum\limits_{\ci\in\mathcal{K}_{\ck}}  \rho^{n}_{\ci} \nn_{\ck\ci},
\end{gather}
with 
\begin{equation}
	 c^{n}_{\ck} = \sum\limits_{\ci\in\mathcal{K}_{\ck}}  \frac{\gamma \left| \Pcell{\ck\ci}\right| \press_{\ck\ci}^{n}}{\left| \Pcell{\ck}\right| \rho_{\ci}^{n}}, \qquad  \press_{\ck\ci}^{n} = \frac{1}{2}\sum_{\cm=1}^{2}\press_{\ck\ci\cm}^{n},
\end{equation}
$\press_{\ck\ci\cm}^{n}$ the pressure at vertex $\cm$ of edge $\Gamma_{\ck\ci}$, $\Gamma_{\ck\ci}$ the primal edge of element $\Pcell{\ck}$ used to generate $\cell{\ci}$,
$\bvel^{n}_{\ck}$ the velocity interpolated from the dual grid to the primal element $\Pcell{\ck}$, following Equation \eqref{eqn.interpolation_dualcellprimalcell}, 
$\un_{\ck\ci}$ the unitary outward pointing normal of  $\Gamma_{\ck\ci}$, and $\nn_{\ck\ci}=\left| \Gamma_{\ck\ci} \right| \un_{\ck\ci}$.

\subsection{Projection stage. Finite element method in the primal grid}
The pressure subsystem associated with the incompressible or the weakly compressible GPR models is solved by employing continuous finite element methods in the primal grid. Focusing on the weakly compressible GPR model, we substitute \eqref{eqn.incGPR_disc_press2} into \eqref{eqn.incGPR_disc_press1}, obtaining
\begin{equation*}
	\frac{1}{\Delta t}\left( \press^{n+1}- \widetilde{\press}\,\right) + c^{2} \dive \left(  \tbmom\right)  - c^{2} \Dt \lapla \press^{n+1}= 0.
\end{equation*}
Next, multiplying the former equation by a test function $z\in V_{0}$, $$V_{0}= \left\lbrace z\in H^{1}(\Omega) \mid \int\limits_{\Omega} z \dV = 0\right\rbrace,$$ integrating in the computational domain $\Omega$, and using Green theorem, we get the weak problem
\begin{weakproblem}
	Find $\press \in V_{0}$ such that
\begin{equation*}
	\frac{1}{c^2}\int\limits_{\Omega} \press^{n+1} z \dV 
	+ \Dt^2 \int\limits_{\Omega} \gra \press^{n+1} \cdot \gra z \dV
	= \frac{1}{c^2} \int\limits_{\Omega} \widetilde{\press} z \dV 
	+ \Dt \int\limits_{\Omega}  \tbmom \cdot \gra z \dV
	- \Dt \int\limits_{\partial \Omega}  \bmom^{n+1} \cdot \un z \dS
\end{equation*}
for all $z\in V_{0}$.
\end{weakproblem}
Similarly, for the incompressible system \eqref{eqn.incGPR_disc_press}, we have
\begin{weakproblem}
	Find $\delta\press^{n+1} \in V_{0}$ such that
	\begin{equation*}
		\int\limits_{\Omega} \gra \delta \press^{n+1} \cdot \gra z \dV
		=  \frac{1}{\Dt} \int\limits_{\Omega}  \tbmom \cdot \gra z \dV
		- \frac{1}{\Dt} \int\limits_{\partial \Omega}  \bmom^{n+1} \cdot \un z \dS
	\end{equation*}
	for all $z\in V_{0}$.
\end{weakproblem}
Finally, we employ the second order $\mathbb{P}_{1}$ continuous finite element method to discretise the weak problems, and the resulting algebraic systems are solved using a matrix-free conjugate gradient method. 
Let us note that the obtained weak problems correspond to the ones arising for the incompressible and weakly compressible Navier-Stokes equations when the splitting procedure \cite{TV12} is considered; further details on the applied methodology can be found in \cite{BFTVC17,Hybrid1}.

\subsection{Correction stage}
Once the new pressure has been obtained as the solution of the projection stage, the intermediate momentum is updated at each dual cell $\cell{\ci}$ using \eqref{eqn.incGPR_disc_press2} and \eqref{eqn.WGPR_semidiscrete_mompress2},
for the incompressible GPR model and the weakly compressible GPR model, respectively. The involved pressure gradients are computed as 
\begin{equation}
	\left( \gra  \press \right)_{\ci} = \frac{1}{\left|\Pcell{\ck}\right|}\sum\limits_{\ck\in\mathcal{T}_{\ci}} 
	\left|\Pcell{\ck\ci}\right| \left( \gra \press\right)_{\ck},
\end{equation}
with $\mathcal{T}_{\ci}$ the set of primal elements related to $\cell{\ci}$ and $\left( \gra \press\right) _{\ck}$ calculated in the primal cells using the $\mathbb{P}_1$ finite element basis functions.

\subsection{Boundary conditions}\label{sec.bc}
Before assessing the proposed methodology, we briefly introduce the main types of boundary conditions employed in Section~\ref{sec:numericalresults}.

If a periodic solution is expected, we may simply employ periodic boundary conditions, which are implemented assuming each pair of boundaries to have a periodic mesh. Consequently, the vertex of the primal elements where the pressure is computed can be merged. On the other hand, we defined each couple of dual elements related through the boundary as a unique dual cell for the dual grid. Then, the solution at the boundary cells is computed as if they were interior elements.

Regarding Dirichlet boundary conditions, two different subcases are considered: strong and weak boundary conditions. If strong boundary conditions are selected, the values of the conservative variables are directly imposed as the solution in the boundary cells. Alternatively, weakly Dirichlet boundary conditions assume the exact solution to be located at the boundary. Hence the fluxes and gradients needed to compute the explicit stage are approximated by setting the given values in the neighbouring ghost cells. For both kinds of boundary conditions, we can further set the pressure as a Dirichlet boundary condition by defining its values at the boundary vertex. Nevertheless, in most cases, we simply employ Neumann boundary conditions for the pressure field. 

Focusing on the fluid limit of the model and considering the presence of walls, the velocity and pressure fields at the boundary cells are computed as for a Navier-Stokes solver, see e.g. \cite{BRMVCD21}. Nevertheless, special care must be paid to the approximation of the distortion field, which is left ``free'' at the boundary. So that, a specific approach to compute the distortion field at the neighbouring of the wall is required.
We first rewrite the distortion field equations as
\begin{equation}
	\frac{1}{\Delta t}\left( \A_{ik}^{n+1}\!-\!\A_{ik}^{n}\right) 
	+ \A_{im}^{n} \pxk \vel_{m}^{n}
	+ \vel_{m}^{n} \pxm \A_{ik}^{n}
	= -\frac{1}{\theta_{1}^{n+1} \left(\tau_{1}\right)} \E_{\A_{ik}}^{n+1},\label{eqn.distortionncp}
\end{equation}
which corresponds to an asymptotic preserving scheme in the fluid relaxation limit of the equations, i.e., as $\tau_{1}\rightarrow 0$ we recover the Navier-Stokes equations, \cite{Boscheri2021SIGPR}. Further, this allows the flux and the non-conservative product contributions to be reordered into two terms. The first one is linear with respect to the distortion field, $\bA \cdot \gra \bvel$, and is simply discretised using the gradients obtained at the dual cells applying the Galerkin approach, i.e., as introduced for the smooth part of the non-conservative products \eqref{eqn.NCP_smooth}. Meanwhile, for the second term, $\bvel \gra \bA$, the path conservative methodology proposed in Section~\ref{sec.NCP} is employed. Then, the term $\bA \cdot \gra \bvel$ can be treated implicitly by including it within the source term system. Consequently, the explicit stage of the algorithm neglects the presence of this term. Then, instead of solving \eqref{eqn.sourcesystem_sd} for the distortion field, at the wall boundary cells we have the implicit system
\begin{subequations}\label{eqn.sourcesystemboundary_sd}
	\begin{eqnarray}
		\A_{ik}^{n+1}   +\frac{ \Delta t}{\theta_{1}^{n+1} \left(\tau_{1}\right)} \E_{\A_{ik}}^{n+1} =  \A_{ik}^{\star}, \\
		\A_{ik}^{\star}=\A_{ik}^{n}- \Delta t \,\A_{im}^{n} \pxk \tilde{\vel}_{m}-\Delta t\,\vel^{n}_{m}\pxm\A_{ik}^{n},
	\end{eqnarray}
\end{subequations}
with the velocity gradients, $\gra \tilde{\bvel}$, approximated using the Galerkin approach by setting the known velocity at the Crouzeix-Raviart node located at the boundary and the previous time step velocity values at the internal nodes. Let us remark that the use of \eqref{eqn.distortionncp} for the distortion field computation instead of the original version in terms of a convective flux term, \eqref{eqn.ingpr_disc_transp.A}, leads to an asymptotic preserving scheme in the fluid limit of the model, further details can be found in \cite{Boscheri2021SIGPR}.

\section{Numerical results} \label{sec:numericalresults} 
In this section, we assess the proposed methodology through a set of classical test problems for both incompressible and low Mach number flows. The test cases are described employing the international system of units.
Besides, unless stated the contrary, all test cases are run with a variable time-stepping attending to a CFL $= 0.5$ so that stability of the explicit part of the scheme is guaranteed. Taking into account that the pressure subsystem and the source terms for the distortion, heat flux and pressure equations are treated implicitly, the time step restriction reads
\begin{equation*}
	\Dt = \min_{\cell{\ci}} \left\lbrace \Dt_{\ci}\right\rbrace, \qquad  \Dt_{\ci} = \frac{ \mathrm{CFL} \, r_{\ci} }{ \left|\lambda_{\ci} \right|_{\max} },
\end{equation*}
with $r_{\ci}$ the incircle diameter of $\cell{\ci}$ and $\left|\lambda_{\ci} \right|_{\max}$ the maximum absolute eigenvalue associated to the explicit subsystem.

\subsection{Convergence study: Taylor-Green Vortex}
As a first test case, we consider the 2D Taylor-Green vortex benchmark whose known exact solution for the Euler equations in $\Omega=[0,2\pi]^2$ is given by
\begin{gather*}
	\vel \left(\mathbf{x},t\right) = \left( \begin{array}{r} 
		\sin(x)\cos(y) \\ 
		-\cos(x)\sin(y) \end{array} \right), \quad 
	\press \left(\mathbf{x},t\right) = \frac{\press_{0}}{\gamma-1} + \frac{1}{4} \left(\cos(2x)+\cos(2y) \right).\vspace{0.5cm}
\end{gather*}

We run a set of simulations for the successively refined triangular grids described in Table~\ref{TGV_mesh} with both the first and second order approaches for the convective terms using the incompressible and weakly compressible GPR models.
In particular, in the incompressible regime, we simply set $p_0=0$, and the model parameters are $c_s = 0$ and $\mu = 0$. Meanwhile, 
to test the weakly compressible code, we take $\press_{0}=10^{5}$, $\gamma=1.4$, $c_{v}=2.5$, $c_s = c_h = 0$, and $\mu = \kappa = 0$, yielding to a characteristic Mach number of $M \approx 1.7 \cdot10^{-3}$. 
In all cases, the expected convergence orders are attained, as shown in Tables~\ref{TGV_INCGPR_errors}-\ref{TGV_WGPR_errors}.
\begin{table}[H]
	\renewcommand{\arraystretch}{1.2}
	\begin{center}
		\begin{tabular}{cccc}
			\hline 
			Mesh & Elements & Vertices & Dual elements \\\hline
			$M_1$ & $128 $ & $81 $ & $208 $ \\ 
			$M_2$ & $512 $ & $289 $ & $800 $ \\ 
			$M_3$ & $2048 $ & $1089 $ & $3136 $ \\ 
			$M_4$ & $8192 $ & $4225 $ & $12416 $ \\ 
			$M_5$ & $32768 $ & $16641 $ & $49408 $ \\ 
			$M_6$ & $131072 $ & $66049 $ & $197120 $ \\ 
			\hline \\[0.01cm]
		\end{tabular}
		\caption{2D Taylor-Green vortex. Main features of the primal triangular grids used to run the convergence table.} \label{TGV_mesh}
	\end{center}
\end{table}

\begin{table}[H]
	\renewcommand{\arraystretch}{1.2}
\begin{center}
	\begin{tabular}{ccccc}
		\hline               
		\multirow{2}{*}{Mesh} & \multicolumn{4}{c}{First order scheme}\\\hhline{~----}
		& $L^{2}_{\Omega}\left(\vel\right)$ & $\mathcal{O}\left(\vel\right)$ & $L^{2}_{\Omega}\left(\press\right)$ & $\mathcal{O}\left(\press\right)$ \\ \hline
		M1  & $3.45\cdot 10^{-1} $ & $     $ & $7.42\cdot 10^{-1} $ & $     $  \\
		M2  & $1.99\cdot 10^{-1} $ & $0.79 $ & $2.02\cdot 10^{-1} $ & $1.88 $  \\
		M3  & $1.09\cdot 10^{-1} $ & $0.87 $ & $9.44\cdot 10^{-2} $ & $1.10 $  \\
		M4  & $5.71\cdot 10^{-2} $ & $0.93 $ & $4.64\cdot 10^{-2} $ & $1.03 $  \\
		M5  & $2.93\cdot 10^{-2} $ & $0.97 $ & $2.31\cdot 10^{-2} $ & $1.01 $  \\
		M6  & $1.48\cdot 10^{-2} $ & $0.98 $ & $1.15\cdot 10^{-2} $ & $1.00 $  \\
		M7  & $7.46\cdot 10^{-3} $ & $0.99 $ & $5.75\cdot 10^{-3} $ & $1.00 $  \\
		\hline \\[0.01cm]
	\end{tabular}
	\begin{tabular}{cccc}
	\hline 
	\multicolumn{4}{c}{Local ADER scheme}\\ \hline
	$L^{2}_{\Omega}\left(\vel\right)$ & $\mathcal{O}\left(\vel\right)$ & $L^{2}_{\Omega}\left(\press\right)$ & $\mathcal{O}\left(\press\right)$ \\ \hline
		 $1.23\cdot 10^{-1} $ & $     $ & $4.21\cdot 10^{-1} $ & $     $ \\
		 $3.06\cdot 10^{-2} $ & $2.01 $ & $2.11\cdot 10^{-1} $ & $1.00 $ \\
		 $7.62\cdot 10^{-3} $ & $2.01 $ & $6.35\cdot 10^{-2} $ & $1.73 $ \\
		 $1.90\cdot 10^{-3} $ & $2.00 $ & $1.66\cdot 10^{-2} $ & $1.93 $ \\
		 $4.75\cdot 10^{-4} $ & $2.00 $ & $4.21\cdot 10^{-3} $ & $1.98 $ \\
		 $1.19\cdot 10^{-4} $ & $2.00 $ & $1.06\cdot 10^{-3} $ & $1.99 $ \\
		 $2.97\cdot 10^{-5} $ & $2.00 $ & $2.64\cdot 10^{-4} $ & $2.00 $ \\
		\hline \\[0.01cm]
	\end{tabular}
	\caption{2D Taylor-Green vortex. Incompressible GPR model. Spatial $L_{2}$ error norms and convergence rates at time $t=0.1$.} \label{TGV_INCGPR_errors}
\end{center}
\end{table}

\begin{table}[H]
	\renewcommand{\arraystretch}{1.2}
	\begin{center}
		\begin{tabular}{ccccccc}
			\hline 
			\multirow{2}{*}{Mesh} & $L^{2}_{\Omega}\left(\rho\right)$ & $\mathcal{O}\left(\rho\right)$                  
			&$L^{2}_{\Omega}\left(\vel\right)$ & $\mathcal{O}\left(\vel\right)$ & $L^{2}_{\Omega}\left(\press\right)$ & $\mathcal{O}\left(\press\right)$ \\ \hhline{~------}
			& \multicolumn{6}{c}{Local ADER scheme}\\ \hline
			M1 & $2.30\cdot 10^{-2} $ & $     $ &$1.04\cdot 10^{-1} $ & $     $ & $9.40\cdot 10^{1} $ & $     $ \\
			M2 & $3.23\cdot 10^{-3} $ & $2.83 $ &$2.94\cdot 10^{-2} $ & $1.82 $ & $1.55\cdot 10^{1} $ & $2.60 $ \\
			M3 & $4.54\cdot 10^{-4} $ & $2.83 $ &$7.62\cdot 10^{-3} $ & $1.95 $ & $1.78\cdot 10^{0} $ & $3.12 $ \\
			M4 & $7.71\cdot 10^{-5} $ & $2.56 $ &$1.90\cdot 10^{-3} $ & $2.00 $ & $1.74\cdot 10^{-1} $ & $3.35 $ \\
			M5 & $1.63\cdot 10^{-5} $ & $2.24 $ &$4.75\cdot 10^{-4} $ & $2.00 $ & $1.67\cdot 10^{-2} $ & $3.38 $ \\
			M6 & $4.00\cdot 10^{-6} $ & $2.02 $ &$1.19\cdot 10^{-4} $ & $2.00 $ & $1.71\cdot 10^{-3} $ & $3.29 $ \\
			\hline \\[0.01cm]
		\end{tabular}
		\caption{2D Taylor-Green vortex. Weakly compressible GPR model. Spatial $L_{2}$ error norms and convergence rates at time $t=0.1$.} \label{TGV_WGPR_errors}
	\end{center}
\end{table}

\subsection{Lid-driven cavity}
\begin{figure}[H]
	\begin{center}
		\includegraphics[trim=15 -70 15 100,clip,width=0.45\textwidth]{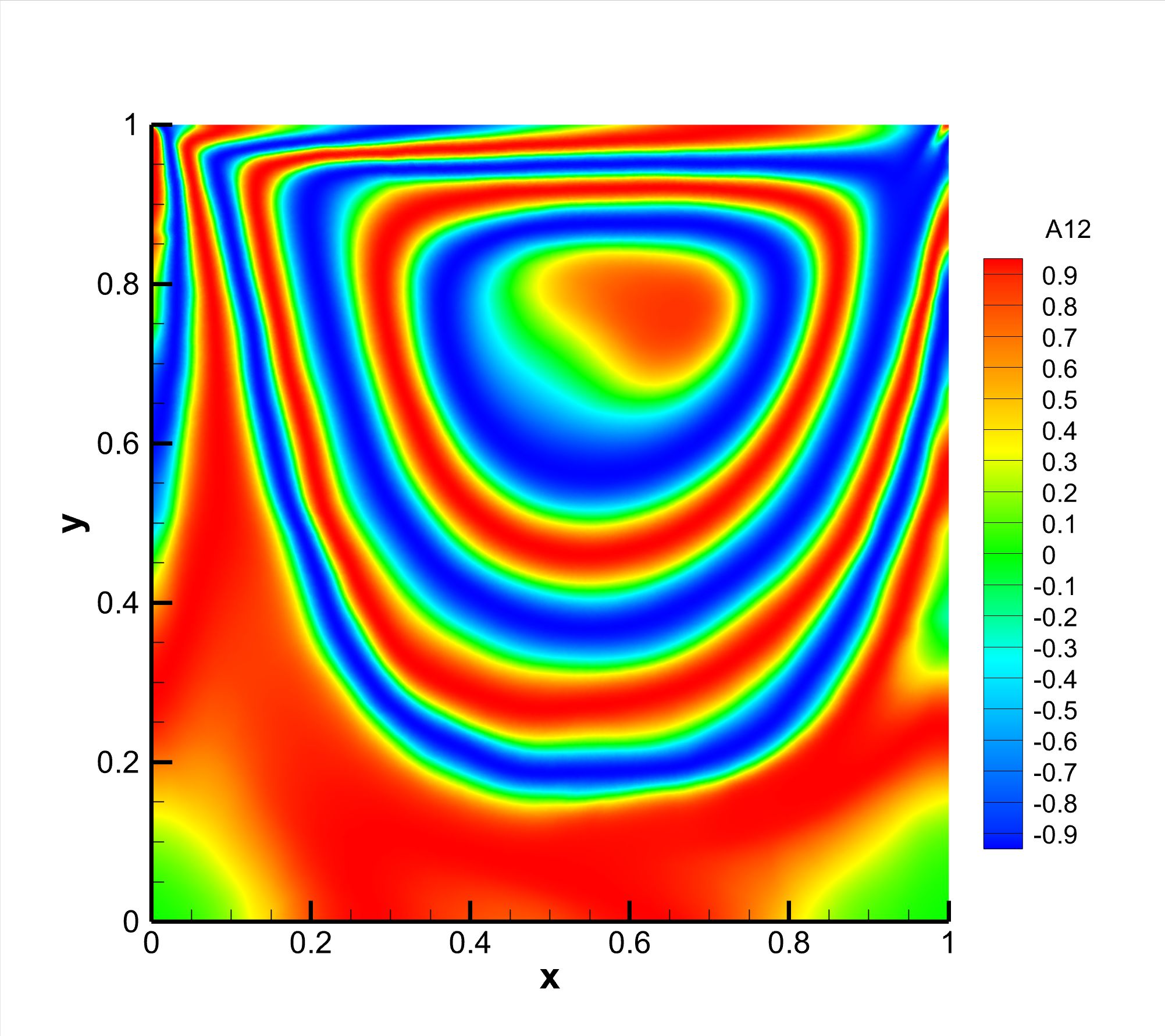}  
		\includegraphics[width=0.45\textwidth]{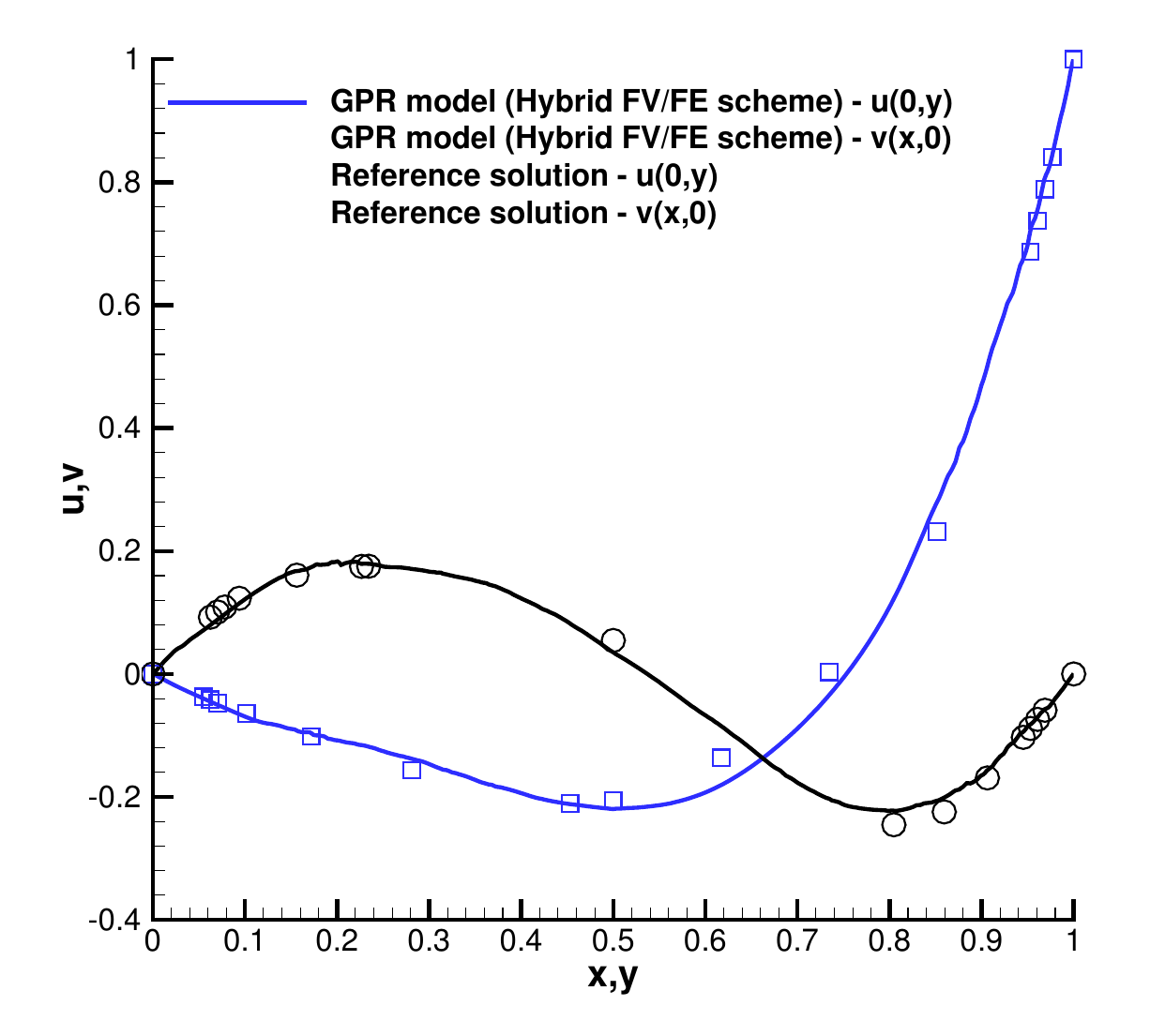}    
		\caption{Lid-driven cavity. Left: contour plot of the distortion component $A_{12}$. Right: 1D cut in $x-$ and $y-$directions of the velocity components $\vel_{2}$ and $\vel_{1}$ computed using the new hybrid FV/FE method for the incompressible GPR model (blue solid line - $\vel_{1}$ and dark grey solid line -  $\vel_{2}$) and reference solutions reported in \cite{GGS82} (blue squares - $\vel_{1}$  and black circles - $\vel_{2}$).}  
		\label{fig.LDC}
	\end{center}
\end{figure}

To analyse the behaviour of the proposed methodology in the incompressible fluid limit, we study the lid-driven cavity benchmark, \cite{GGS82}. We consider an initial fluid at rest with $\bvel=\boldsymbol{0}$, $\press=\rho=1$, $\mathbf{A}=\mathbf{I}$, $c_{s}=8$, and $\mu=10^{-2}$. 
Homogeneous wall boundary conditions are set in the bottom and laterals of the computational domain $\Omega=[-0.5,0.5]^2$, while the upper bound is assumed to be moving horizontally with a lid velocity $\vel_{\mathrm{lid}}=1$. The new hybrid FV/FE methodology is employed to solve the incompressible GPR model up to time $t=10$. The obtained results are depicted in Figure~\ref{fig.LDC}. We can observe a good qualitative agreement with former results available in the bibliography, see, e.g. \cite{Busto2022HTCGPR}. Moreover, Figure~\ref{fig.LDC} also reports the 1D cuts of the velocity field along the horizontal and vertical centerlines of the domain, which compare well with the reference data in \cite{GGS82}.

\subsection{Shear motion}
We now study four shear motion tests in the computational domain $ \Omega=[-0.5,0.5]\times[-0.05,0.05]$ with the initial condition
\begin{gather*}
	\rho\left(\x,0\right) = 1,\quad
	\press \left(\x,0\right) = \frac{1}{\gamma}, \quad
	{u}_{1} \left(\x,0\right) = 0, \quad
	{u}_{2} \left(\x,0\right) = \left\lbrace \begin{array}{lc}
		-0.1 & \mathrm{ if } \; y \le 0,\\
		0.1 & \mathrm{ if } \; y > 0,
	\end{array}\right. \\ \mathbf{A}\left(\x,0\right)=\mathbf{I}, \quad \mathbf{J}=\boldsymbol{0}.
\end{gather*}

\begin{figure}[H]
	\begin{center}
		\includegraphics[width=0.46\textwidth]{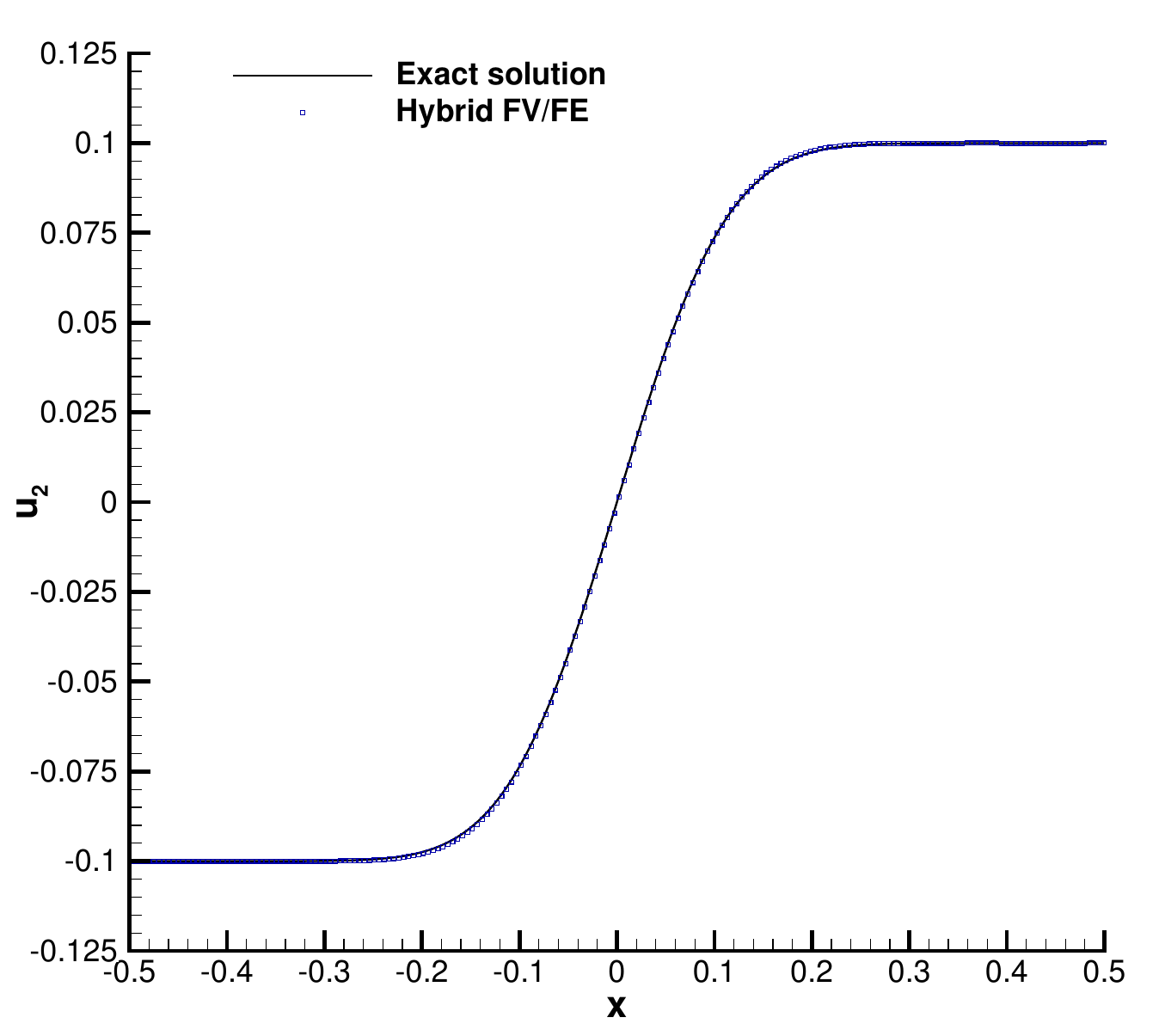}  
		\includegraphics[width=0.46\textwidth]{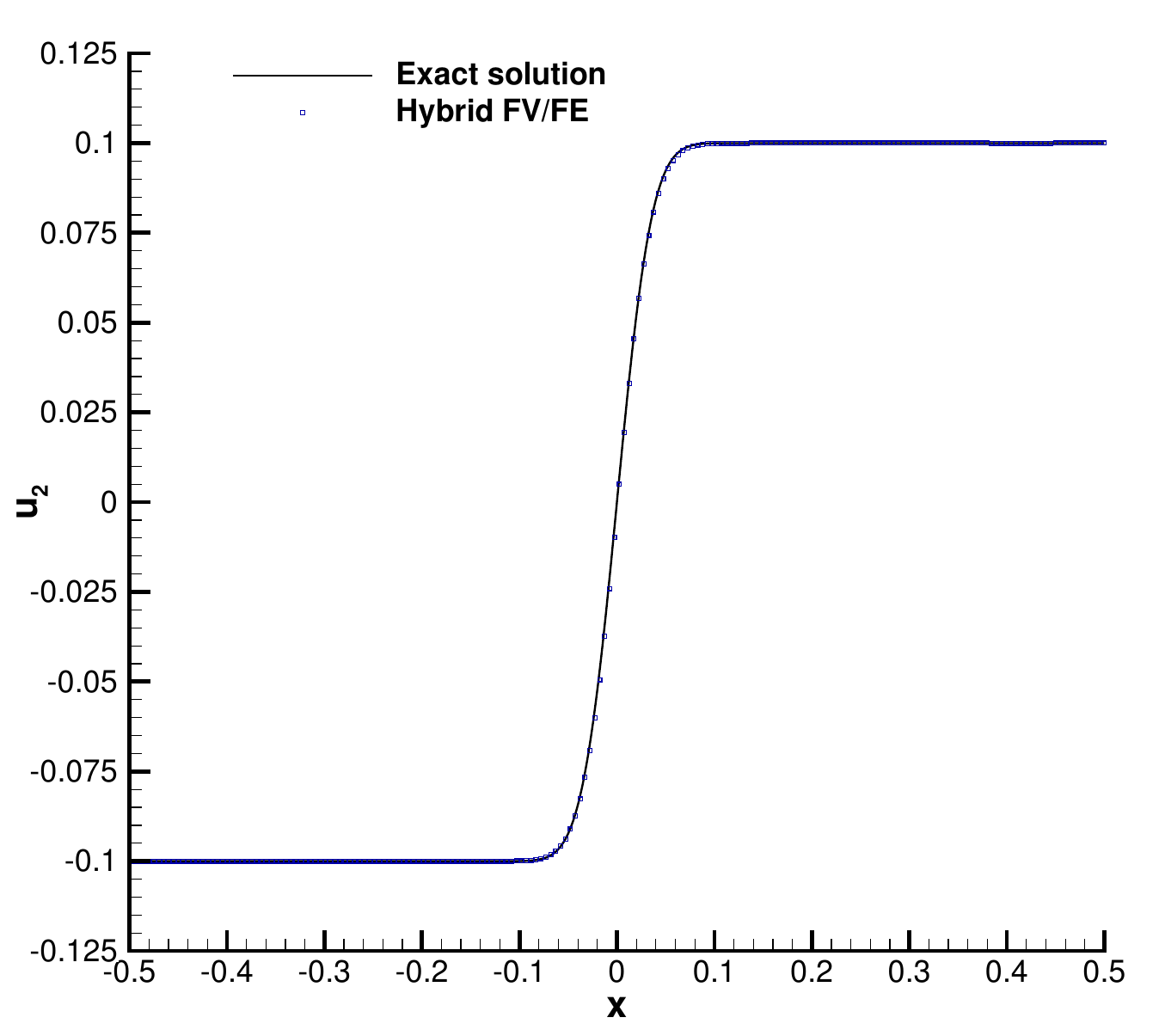}    \\ 
		\includegraphics[width=0.46\textwidth]{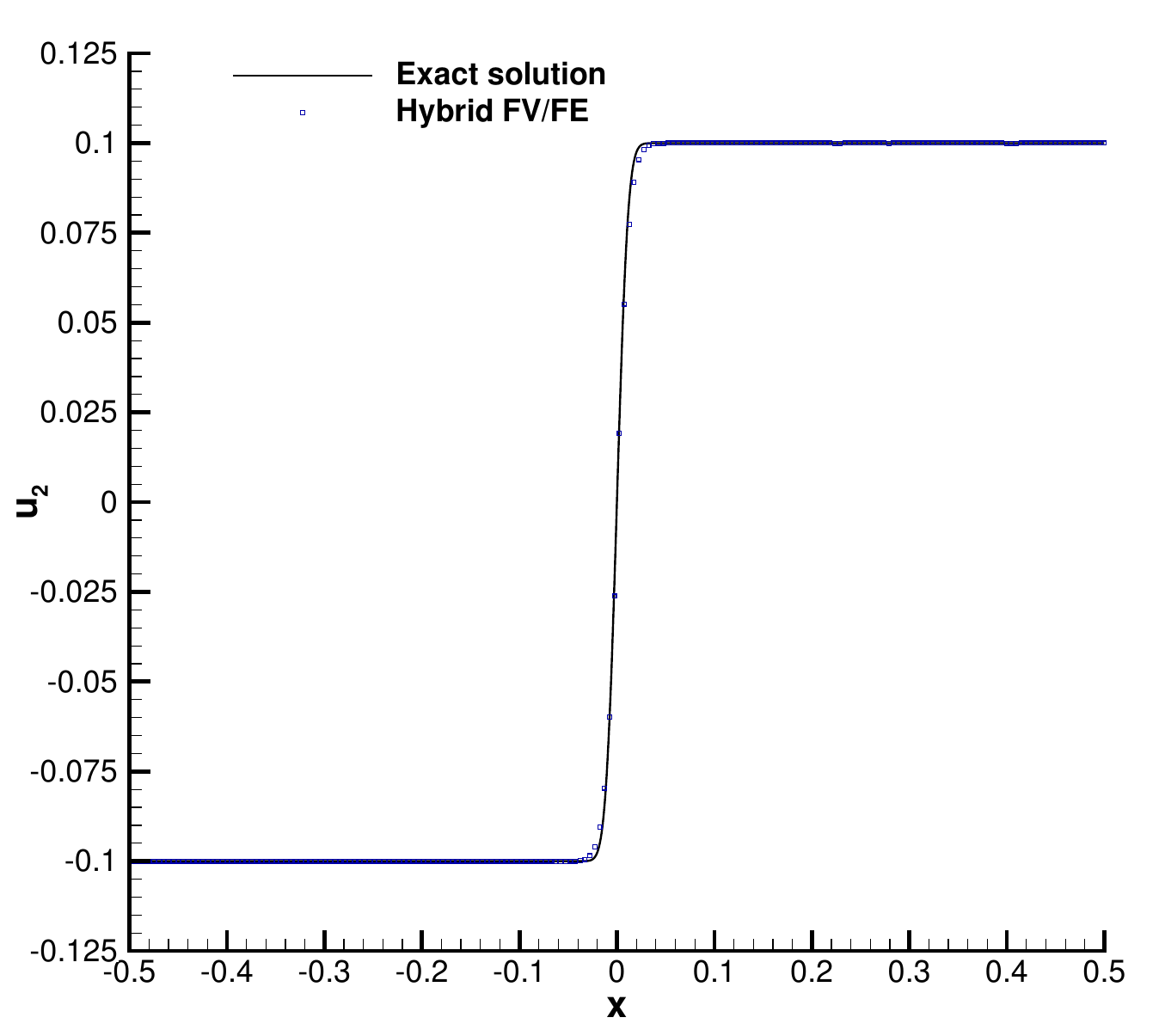}   
		\includegraphics[width=0.46\textwidth]{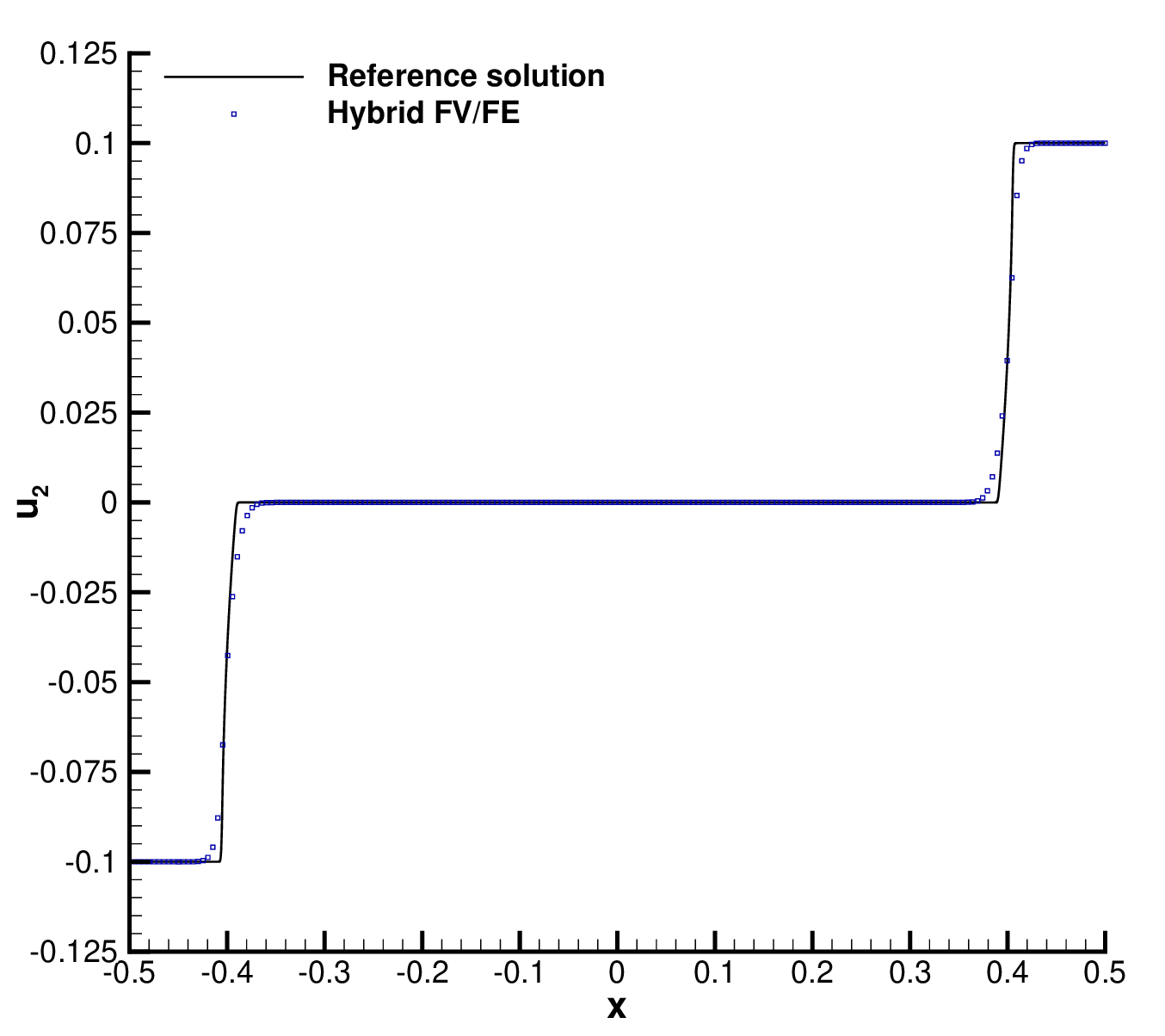}    \\ 
		\caption{Shear motion. 1D cut in $x-$direction of the velocity component $\vel_{2}$ of the numerical solution obtained using the hybrid FV/FE method for the weakly compressible GPR model with the local ADER-ENO approach (blue squares). Reference solution computed with a TVD-FV scheme on a mesh of $1000$ cells (black solid line). From left top to right bottom: first Stokes with $\mu=10^{-2}$, first Stokes with $\mu=10^{-3}$, first Stokes with $\mu=10^{-4}$, shear solid.}  
		\label{fig.FS_SF}
	\end{center}
\end{figure}

Taking the incompressible fluid limit of the GPR model by setting $c_s=c_h=1$, $c_v=2.5$, $\mu\in\left\lbrace 10^{-4}, 10^{-3}, 10^{-2}\right\rbrace$, and $\kappa=\mu$, we recover the well-known first problem of Stokes with known exact analytical solution for the incompressible Navier-Stokes equations given by
\begin{equation*}
	\vel_{2} \left(\x,t\right) = \frac{1}{10} \mathrm{erf}\left( \frac{x}{2\sqrt{\mu t}}\right).\vspace{0.3cm}
\end{equation*}
Since these test cases are run in a 2D domain, we set periodic boundary conditions in $y$-direction while strong Dirichlet boundary conditions are imposed in the left and right boundaries. A primal triangular grid of $N_{x}=200$ divisions is employed for $\mu=10^{-2}$ and $\mu=10^{-3}$, while $N_{x}=400$ is defined for $\mu=10^{-4}$.
The results obtained at time $t=0.4$ using the hybrid FV/FE scheme with the local ADER-ENO approach for the convective terms are reported in Figure~\ref{fig.FS_SF}. Excellent agreement is observed regarding the exact solutions for the three viscosities considered.

As the fourth shear motion test, we set $\tau_1=\tau_2=10^{20}$ to obtain a 1D shear solid benchmark. Then, the GPR model is solved on a triangular grid of $N_{x}=400$ divisions along the $x$-direction. The obtained 1D cut along the centerline of the domain for the velocity component $\vel_{2}$ is depicted in Figure~\ref{fig.FS_SF}. We observe a good agreement with the reference solution computed employing a second order MUSCL-Hancock TVD-FV scheme on a one-dimensional grid of $1000$ control volumes.

\subsection{Double shear layer}
\begin{figure}[H]
	\begin{center}
		\includegraphics[trim=15 50 15 100,clip,width=0.38\textwidth]{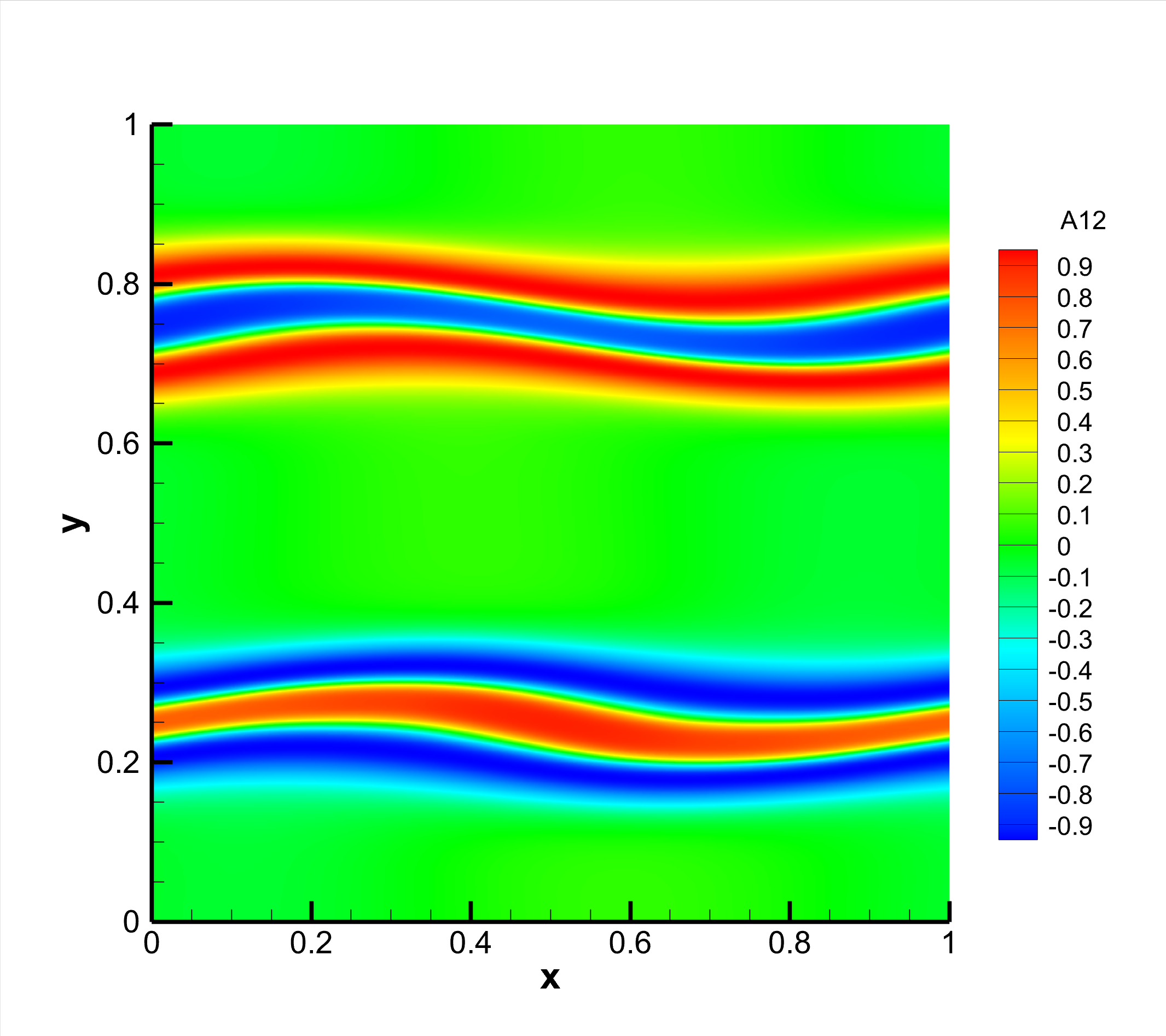} 
		\includegraphics[trim=15 50 15 100,clip,width=0.38\textwidth]{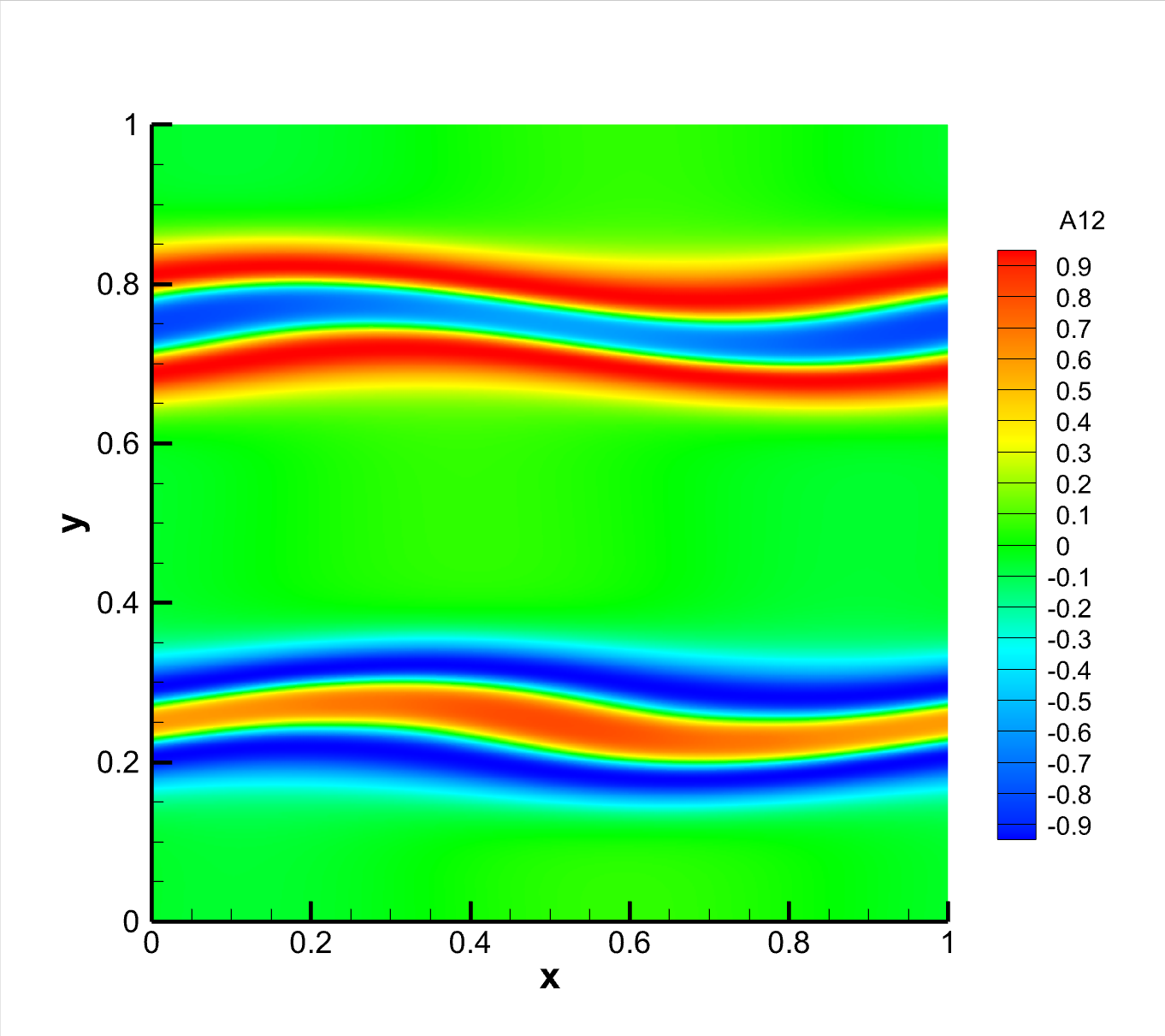}
		\\
		\includegraphics[trim=15 50 15 100,clip,width=0.38\textwidth]{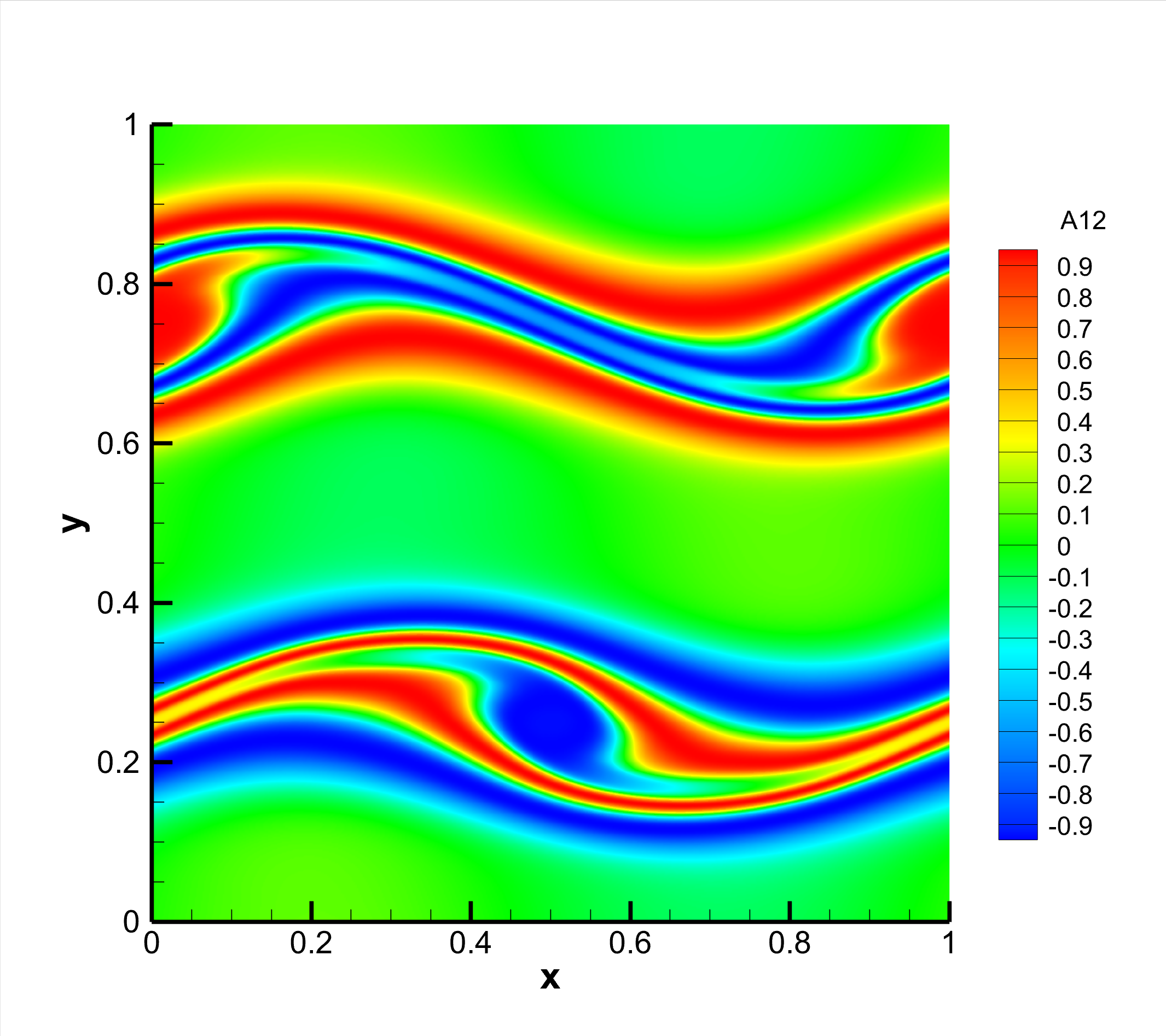}
		\includegraphics[trim=15 50 15 100,clip,width=0.38\textwidth]{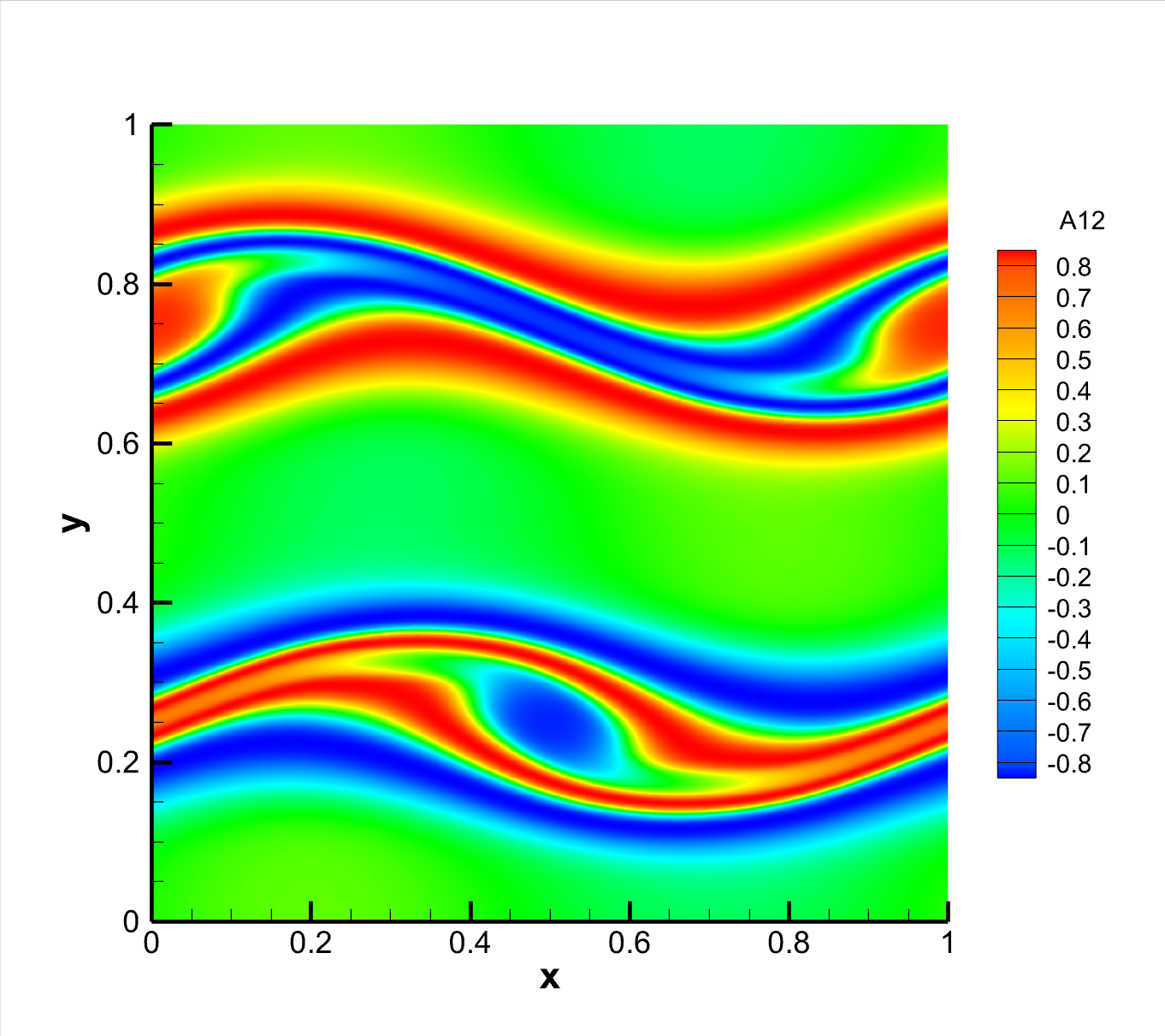}    \\ 
		\includegraphics[trim=15 50 15 100,clip,width=0.38\textwidth]{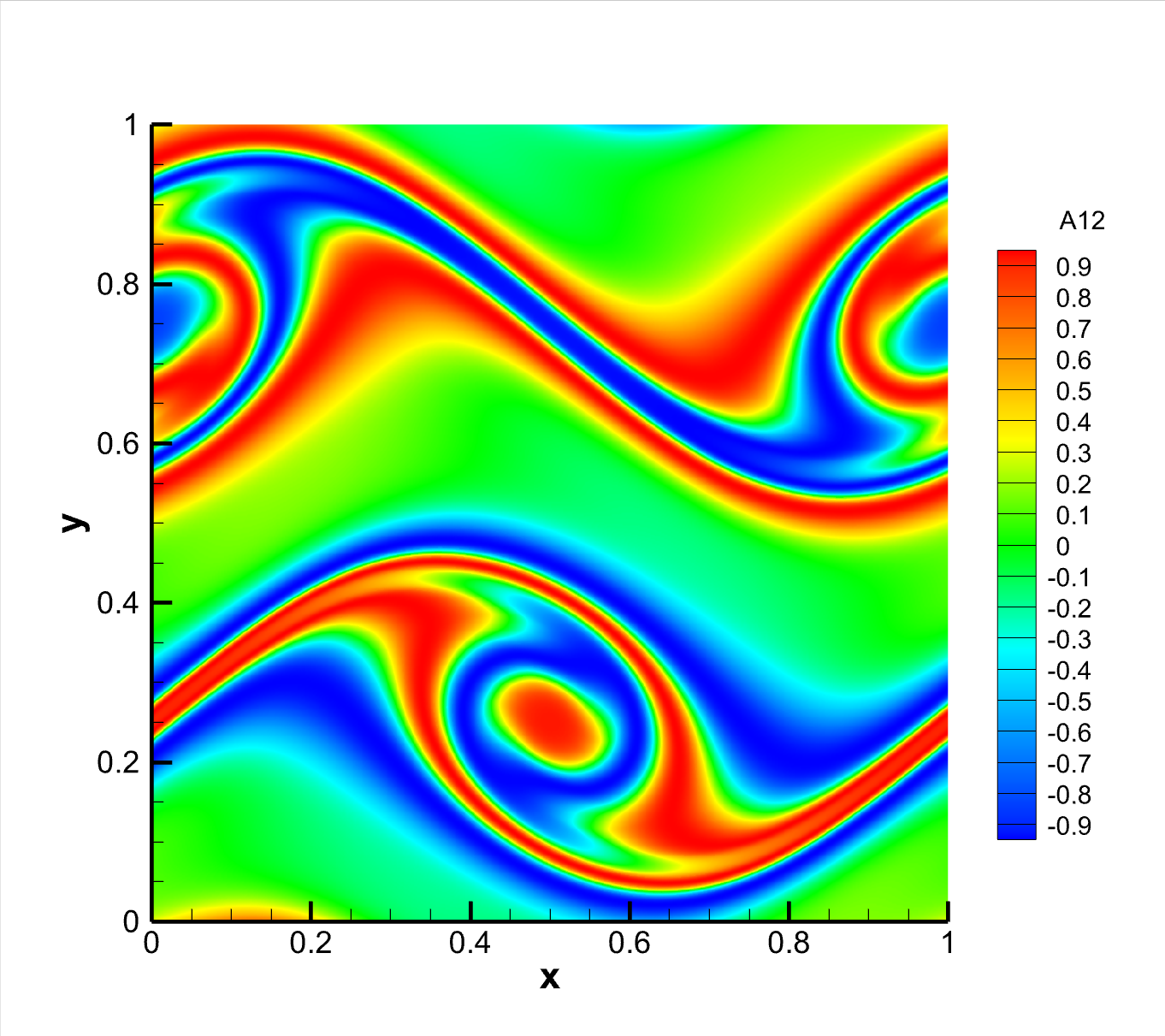} 
		\includegraphics[trim=15 50 15 100,clip,width=0.38\textwidth]{DSL-HybridGPR-12-weakly}  \\ 
		\includegraphics[trim=15 50 15 100,clip,width=0.38\textwidth]{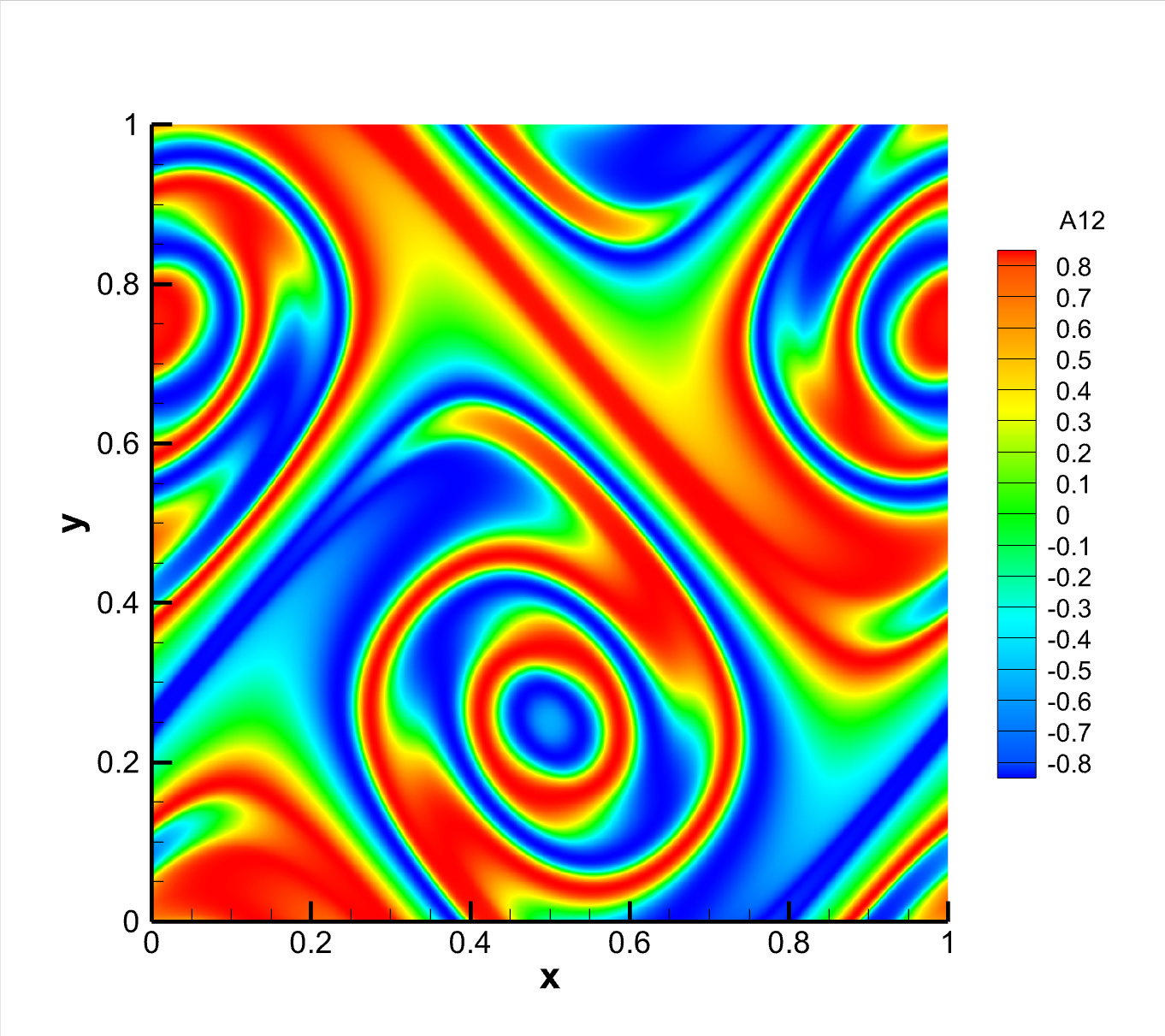}   
		\includegraphics[trim=15 50 15 100,clip,width=0.38\textwidth]{DSL-HybridGPR-18-weakly}    
		\caption{Double shear layer. Contour plots of the distortion field component $A_{12}$ obtained with the new semi-implicit hybrid FV/FE solver at times $t\in\left\lbrace 0.4,0.8,1.2,1.8\right\rbrace$ (from top to bottom). Numerical results were obtained with the incompressible GPR (left) and the weakly compressible GPR (right).}  
		\label{fig.DSLdistortion}
	\end{center}
\end{figure}

As the fourth test case, we consider the double shear layer benchmark whose initial condition, defined in $\Omega=\left[0,1\right]^2$, reads
\begin{gather*}
	\rho\left(\mathbf{x},0\right)  = 1,   \qquad
	\vel_1\left(\mathbf{x},0\right) =\left\lbrace \begin{array}{ll}
		\tanh\left( \tilde{\rho} (y-0.25) \right) & \textrm{if } y \leq 0.5,\\
		\tanh\left( \tilde{\rho} (0.75-y) \right)  & \textrm{if } y > 0.5,
	\end{array}
	\right. \qquad
	\vel_2\left(\mathbf{x},0\right)= \delta \sin(2\pi x) \notag\\
	p\left(\mathbf{x},0\right) = 0,  \qquad
	\mathbf{A}\left(\mathbf{x},0\right)=\mathbf{I},\qquad \mathbf{J}\left(\mathbf{x},0\right)=\mathbf{0},\qquad 
	\delta=0.05,\qquad \tilde{\rho}=30.
\end{gather*}
The numerical solution is computed using both the incompressible and the weakly compressible GPR methodologies with parameters $\mu = 2 \cdot 10^{-3}$, $\kappa = 4 \cdot 10^{-2}$,  $c_{v}=2.5$, $c_{h}=2$ and $c_s=8$. Further, the CFL is set to $0.1$. Periodic boundary conditions are defined everywhere, and a computational grid formed by $2097152$ primal triangular elements is considered. Figure~\ref{fig.DSLdistortion} shows the contour plots of the distortion field component $A_{12}$ at times $t\in\left\lbrace 0.4,0.8,1.2,1.8\right\rbrace$. The method captures well the very thin structures reported in the literature, e.g. \cite{DPRZ16,HTCTotalAbgrall} where high order explicit finite volume methods and thermodynamically compatible schemes have been used, respectively.


\subsection{Solid rotor}
\begin{figure}[H]
	\begin{center}
		\includegraphics[width=0.49\textwidth]{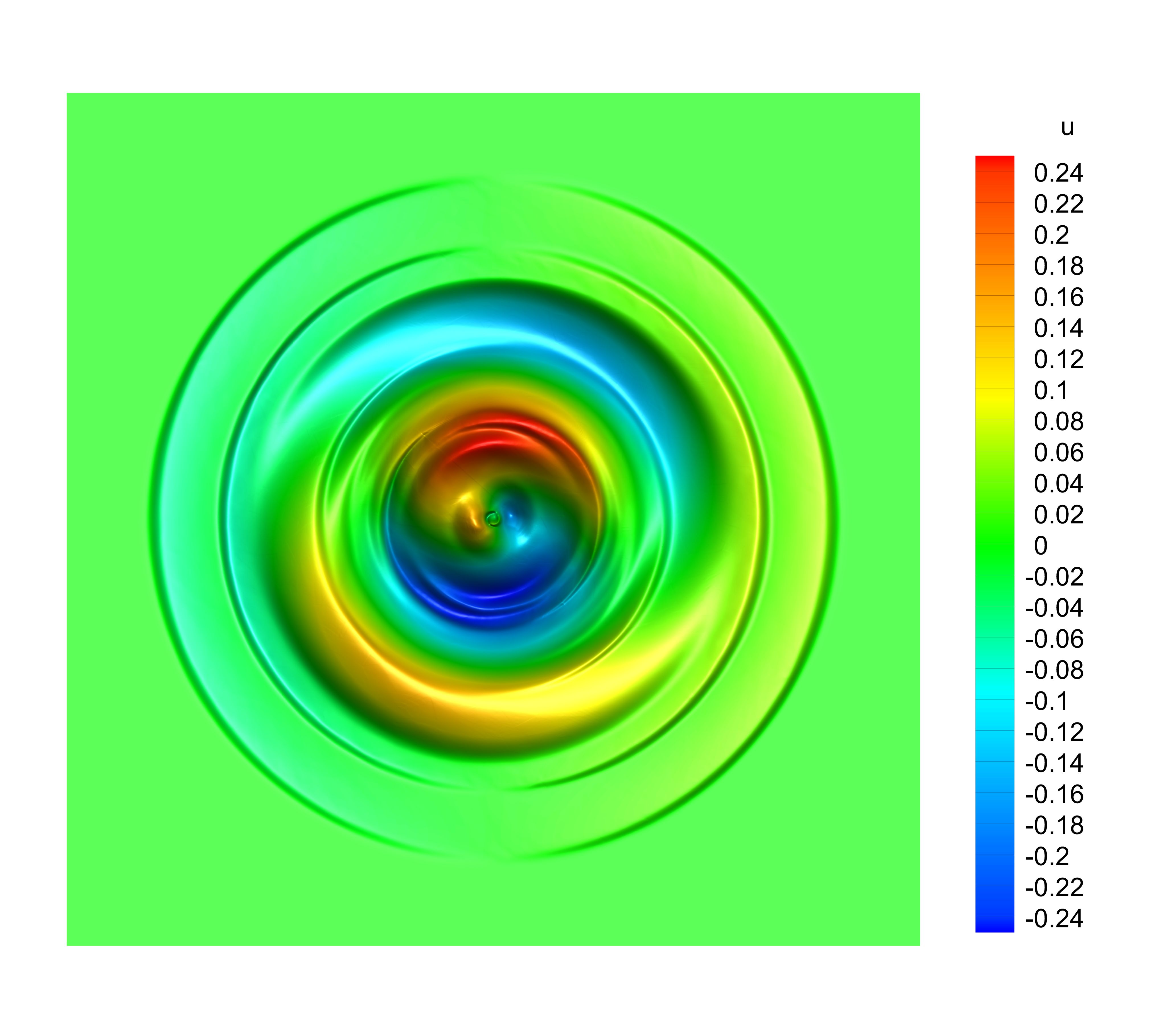}  		\includegraphics[width=0.49\textwidth]{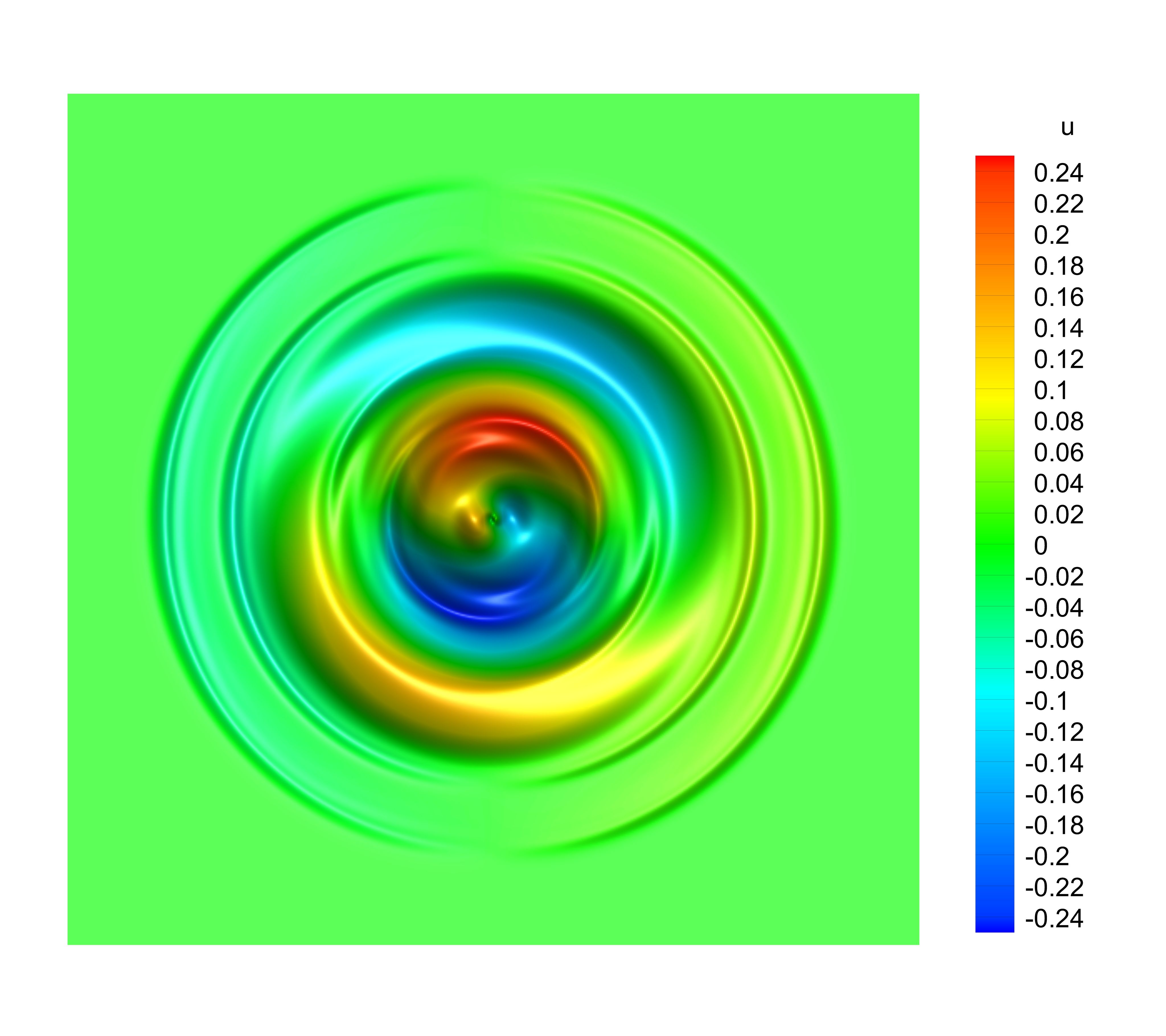}  
		\caption{Solid rotor. Contour plots of the velocity field component $\vel_{1}$ obtained with the hybrid FV/FE approach (left) and the HTC-FV scheme in \cite{HTCTotalAbgrall} (right).}  
		\label{fig.SolidRotor}
	\end{center}
\end{figure}

\begin{figure}[H]
	\begin{center}
		\includegraphics[width=0.45\textwidth]{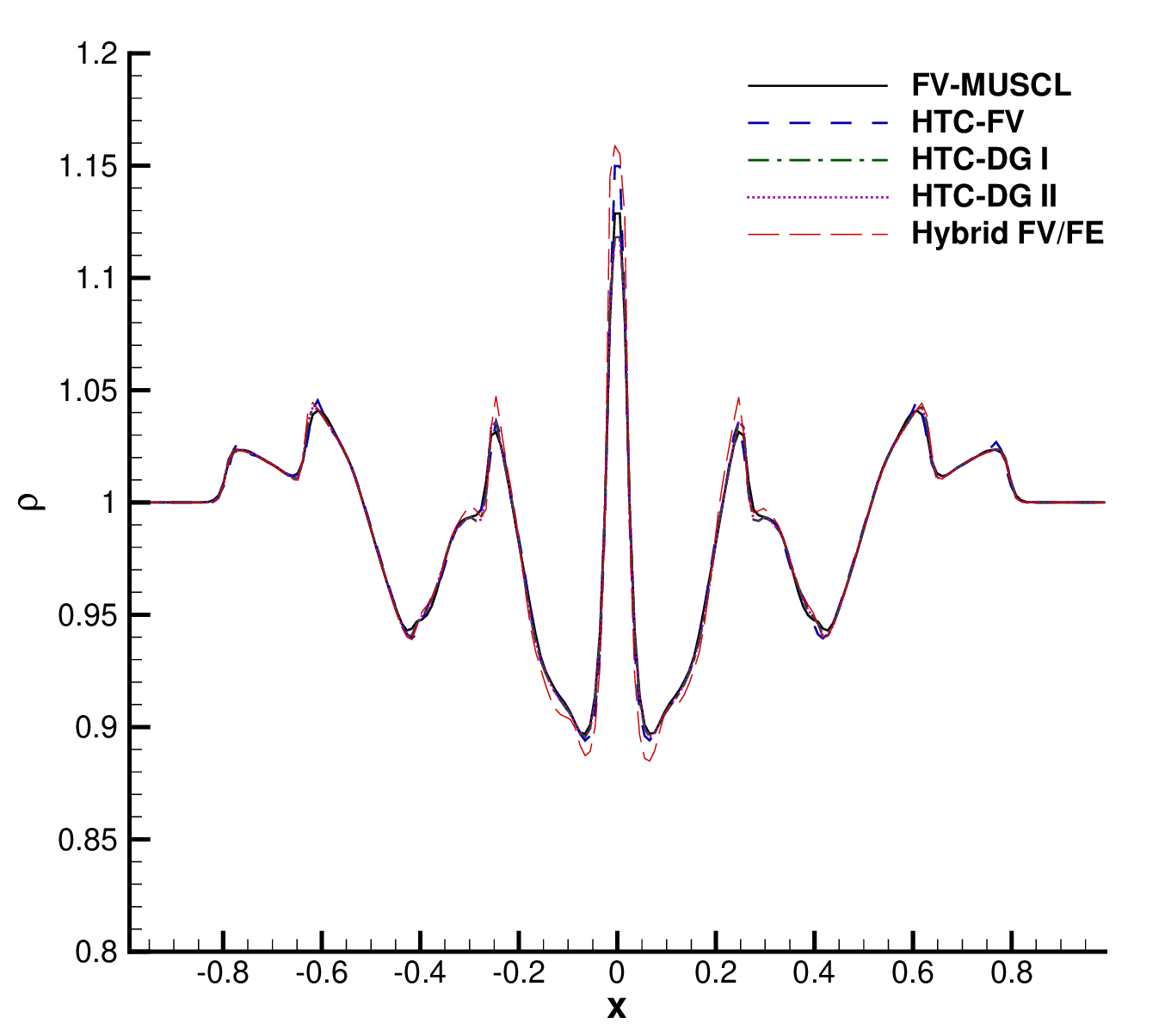}  
		\includegraphics[width=0.45\textwidth]{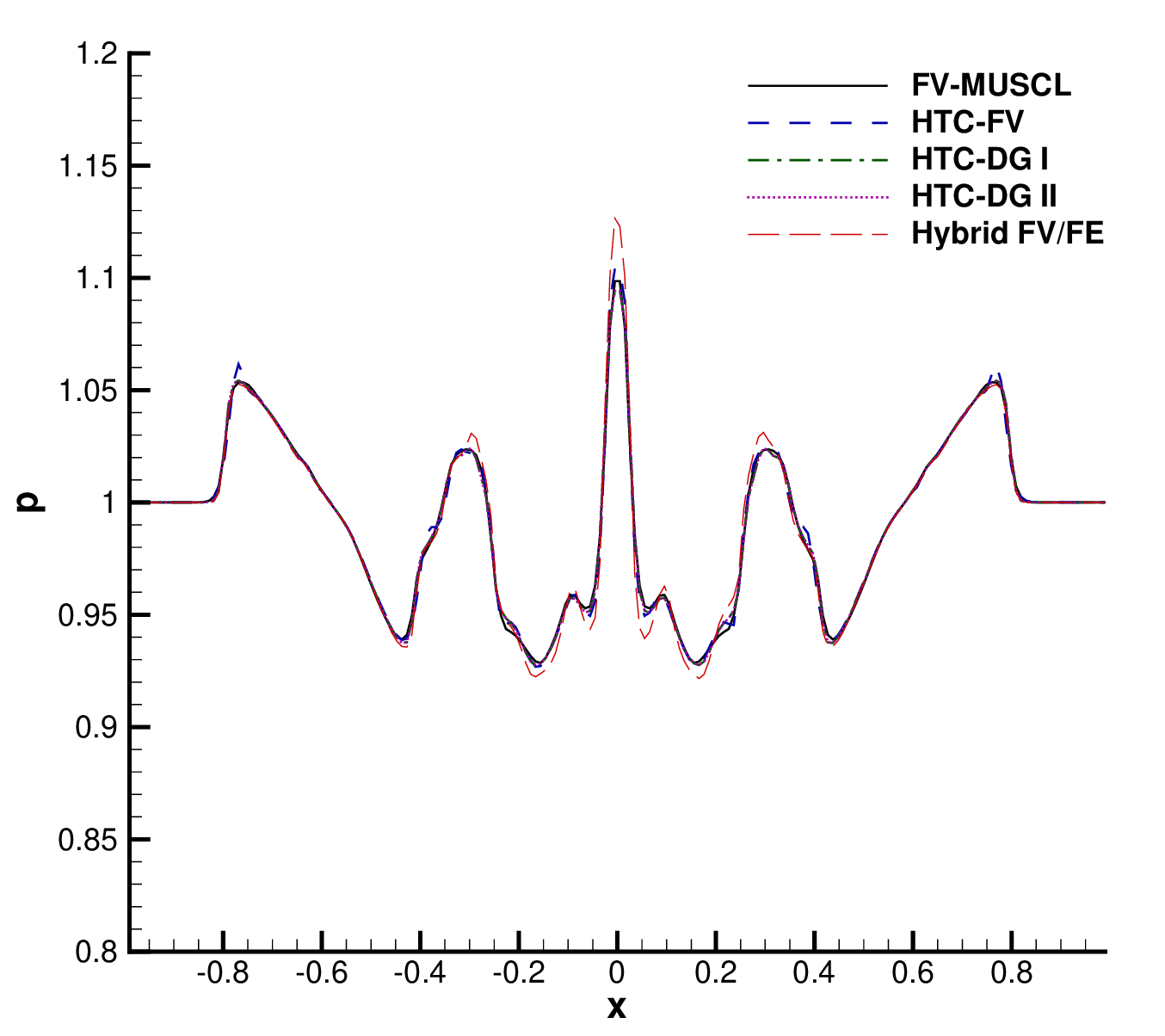}    \\  
		\includegraphics[width=0.45\textwidth]{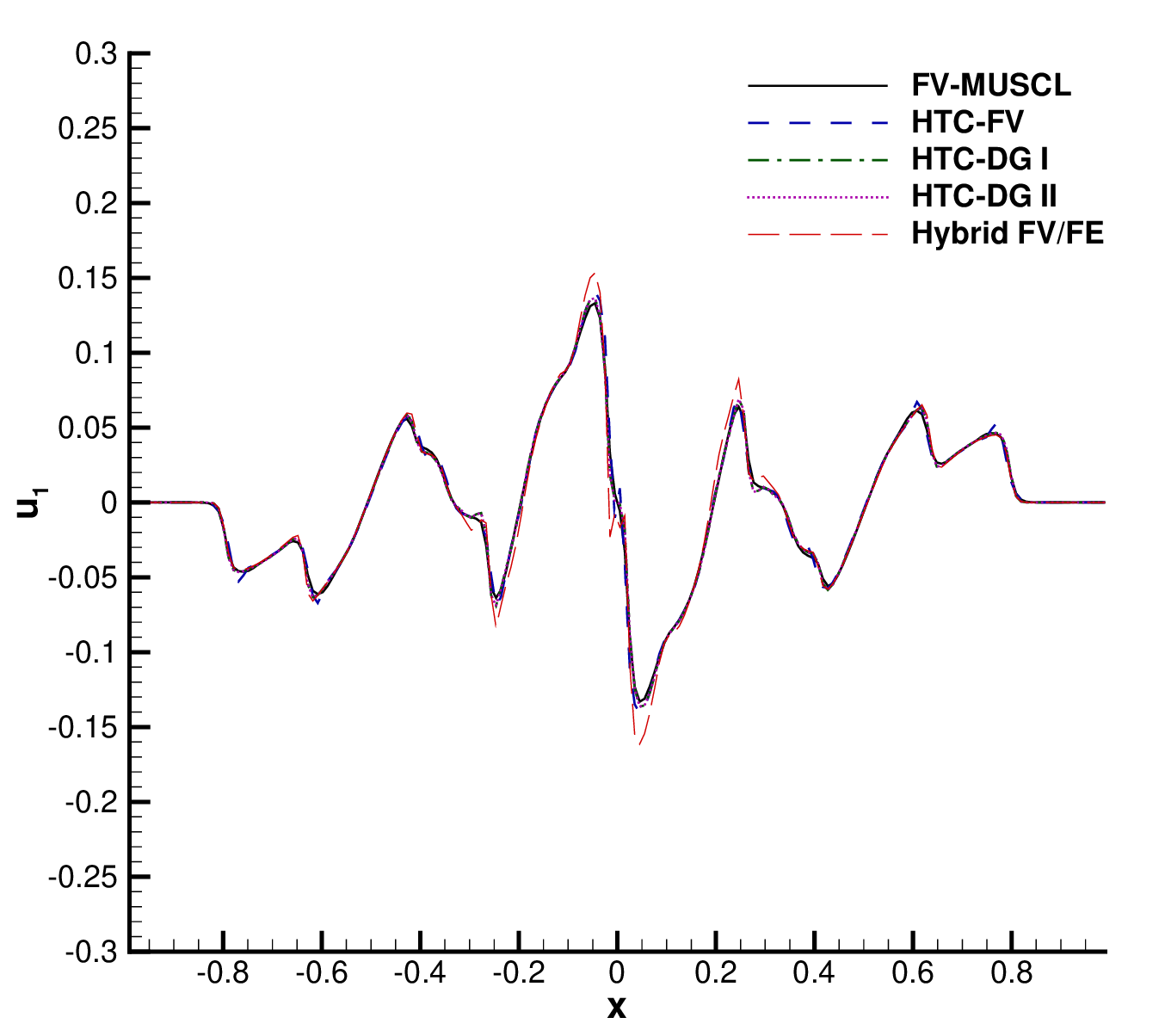}     
		\includegraphics[width=0.45\textwidth]{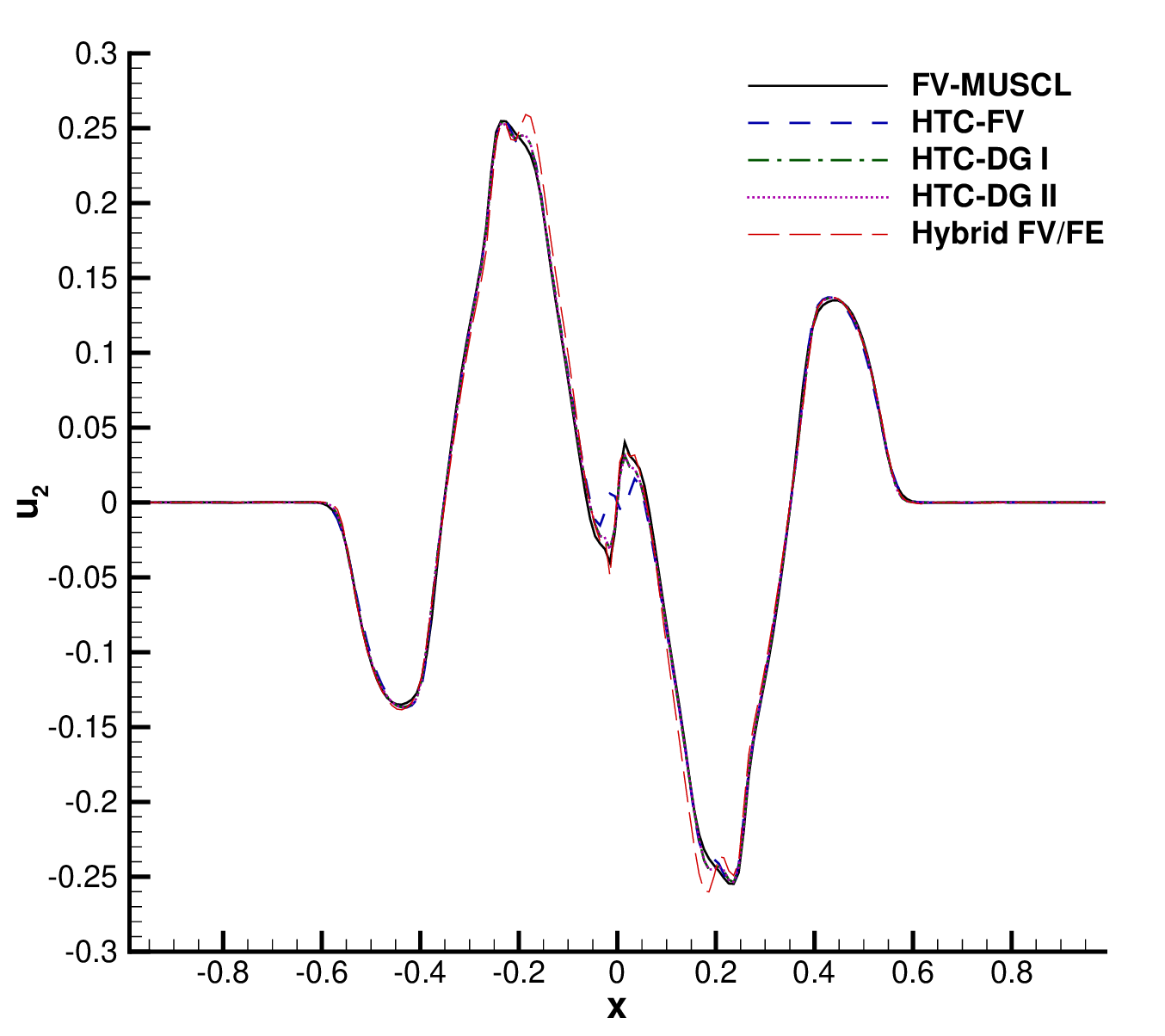}  \\
		\includegraphics[width=0.45\textwidth]{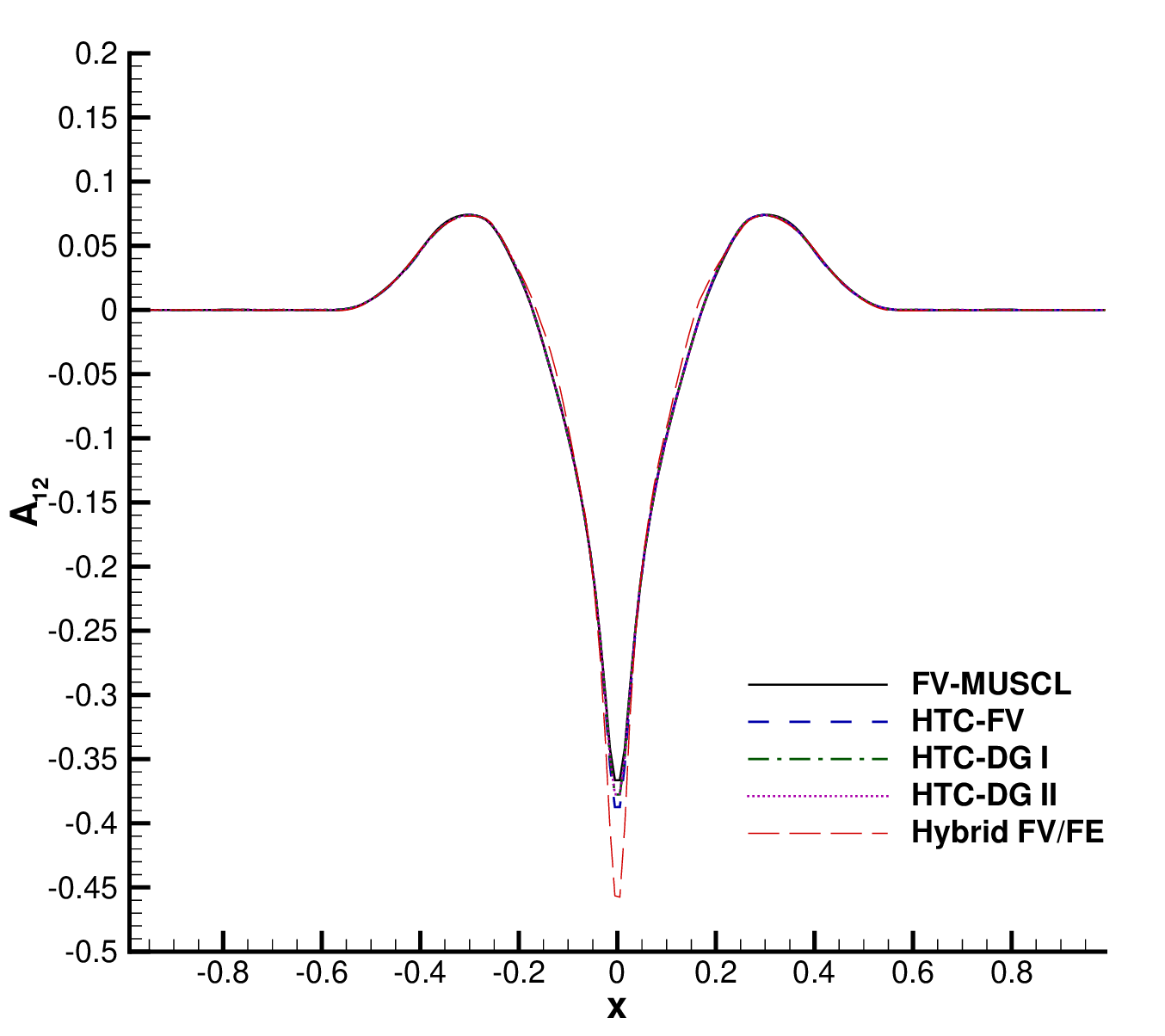}     
		\includegraphics[width=0.45\textwidth]{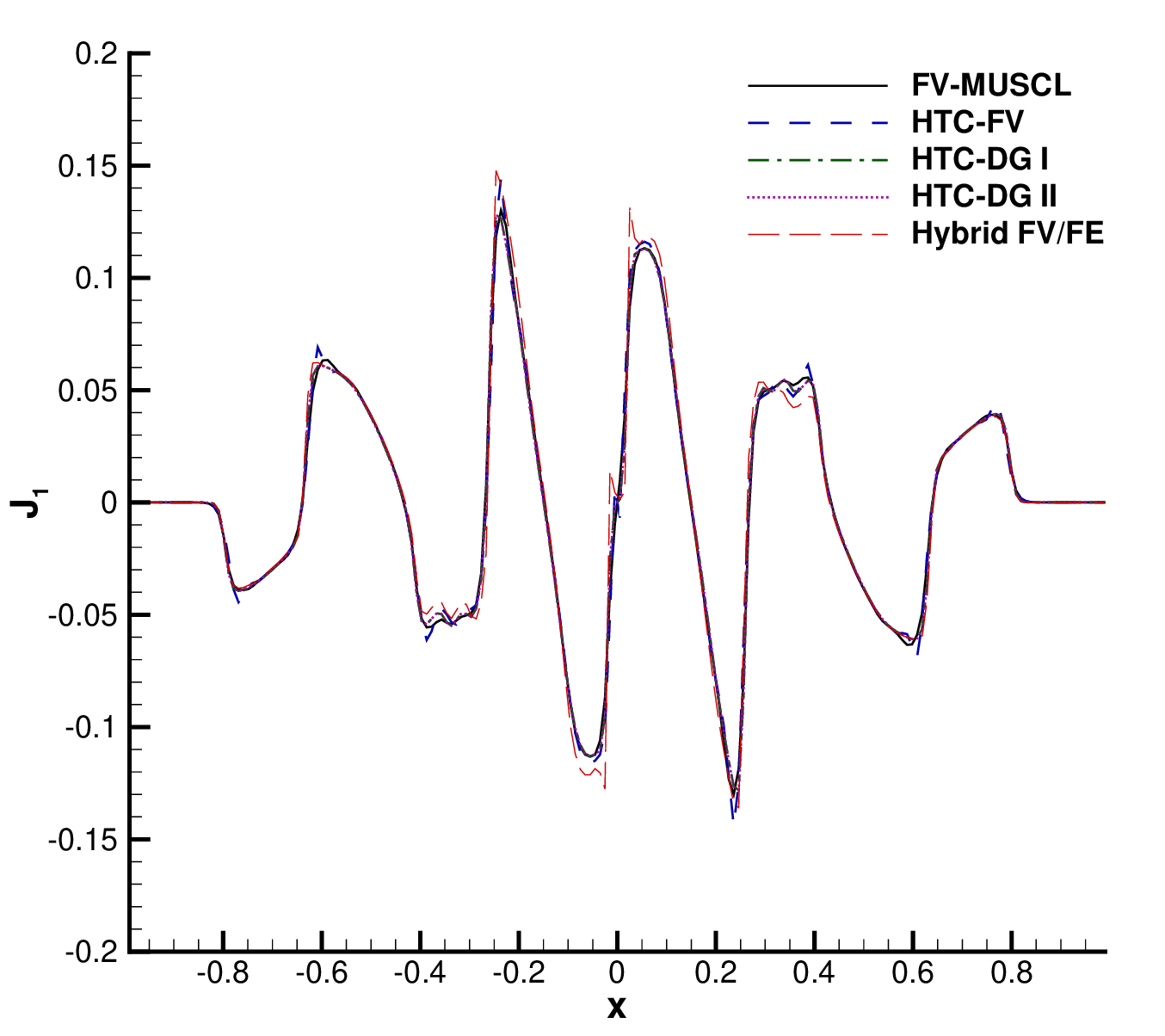}  
		\caption{Solid rotor. From top left to right bottom: 1D cuts of the density, pressure, velocity, $A_{12}$ and $J_1$ fields along $y=0$ at $t=0.3$. The five schemes compared correspond to a second order MUSCL-Hanchock FV scheme (solid black line), the HTC-FV scheme in \cite{HTCTotalAbgrall} (dashed blue line), the HTC-DG I (dash-dotted green line) and HTC-DG-II (dotted purple line) methods with $N=5$ in \cite{Busto2022HTCDG}  and the new hybrid FV/FE approach (long dashed red line). }  
		\label{fig.SolidRotor1D}
	\end{center}
\end{figure}

We now consider the solid rotor test case to analyse further the behaviour of the proposed methodology in the GPR model's solid limit, \cite{Boscheri2021SIGPR}. The initial condition 
\begin{gather*}
	\rho\left(\mathbf{x},0\right) = 1,\quad \bvel\left(\mathbf{x},0\right)=\left\lbrace\begin{array}{ll}
		\left(\frac{-y}{0.2}, \frac{x}{0.2}, 0\right)^{T}  & \mathrm{if} \; \left\| \mathbf{x} \right\| \leq 0.2,\\
		\mathbf{0} & \mathrm{if} \; \left\| \mathbf{x} \right\| > 0.2,
	\end{array} \right.
	\\ \press\left(\mathbf{x},0\right) = 1,\quad  \mathbf{A}\left(\mathbf{x},0\right) = \mathbf{I},\quad \mathbf{J}\left(\mathbf{x},0\right)=\boldsymbol{0}.
\end{gather*}
is defined in the computational domain $\Omega=[-1,1]^2$. Moreover, the model parameters are $\tau_{1}=\tau_{2}=10^{20}$, $\mu =\kappa = 0$, and $c_{s}=c_{h}=c_{v}=1.0$. In Figure~\ref{fig.SolidRotor}, we show the solution obtained with the hybrid FV/FE weakly compressible GPR scheme at time $t=0.3$ employing the local ADER min-mod approach on an unstructured grid with $2975744$ primal triangular elements. The provided reference solution has been obtained using the thermodynamically compatible finite volume scheme presented in \cite{HTCTotalAbgrall}, which solves the entropy-based formulation of the model. 
For comparison, the 1D profiles of the density, velocity, pressure, $A_{12}$ and $J_1$ fields along $y=0$ are reported in Figure~\ref{fig.SolidRotor1D}. A good agreement is observed for all the schemes presented, namely the hybrid FV/FE approach, the HTC-FV scheme in \cite{HTCTotalAbgrall}, the HTC-DG I and HTC-DG-II methods with $N=5$ in \cite{Busto2022HTCGPR} and a second order MUSCL-Hanchock FV approach run on a very fine grid, \cite{Toro}.

\subsection{Riemann problems}
The behaviour of the proposed methodology in presence of strong waves including shocks is analysed through a set of Riemann problems. We consider the computational domain $\Omega=[-0.5,0.5]\times[-0.05,0.05]$ and the initial conditions given by 
\begin{equation*}
	\mathbf{V}\left(\mathbf{x},0\right) = \left\lbrace \begin{array}{lc}
		\mathbf{V}^{L} & \mathrm{ if } \; x \le  x_{c},\\
		\mathbf{V}^{R} & \mathrm{ if } \; x>  x_{c},
	\end{array}\right.
\end{equation*}
with the left and right states for the velocity, density and pressure fields defined in Table~\ref{tab.RP_IC} and $\mathbf{A}_{L}=\mathbf{A}_{R}=\mathbf{I}$, $\mathbf{J}_{L}=\mathbf{J}_{R}=\boldsymbol{0}$. The parameters of the weakly compressible GPR model are set to $c_{s}=c_{h}=0$ and $\mu =\kappa = 0$ for the three classical Riemann problems of the Navier-Stokes equations: RP1, RP2 and RP3. In RP4, a weak viscous fluid with shear is considered by defining $\mu =\kappa = 10^{-5}$ and $c_{s}=c_{h}=1$. On the other hand, Riemann problems RP5 and RP6 correspond to the solid limit of the equations so $\tau_{1}=\tau_{2}=10^{20}$, and we consider $c_{s}=1.0$, $c_{v}=2.5$. Moreover, heat conduction effects are neglected in RP5 by taking $c_{h}=0$, while  $c_{h}=1$ is set for RP6.

\begin{figure}[H]
	\begin{center}
		\includegraphics[trim=5 10 10 10,clip,width=0.32\textwidth]{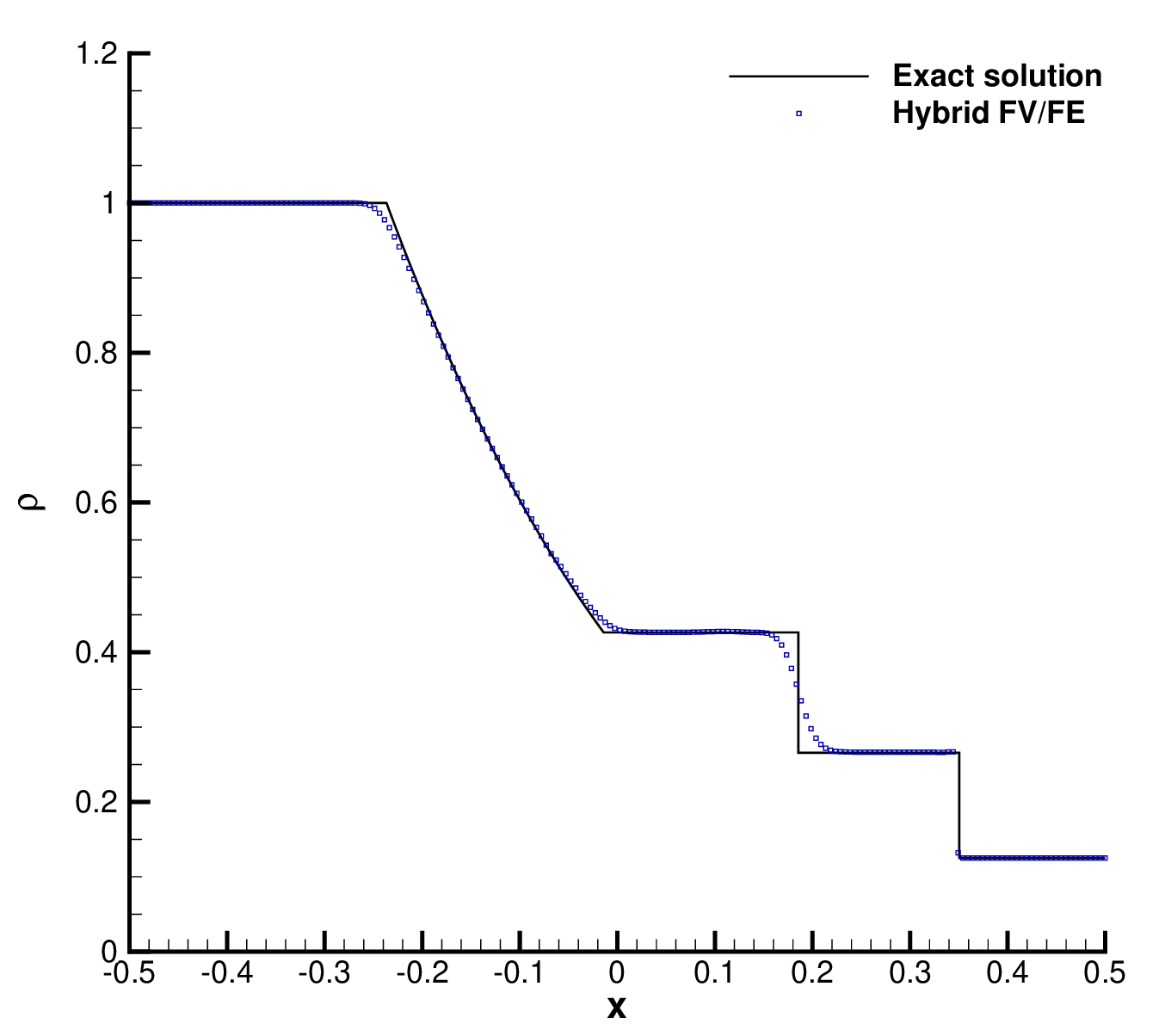}  
		\includegraphics[trim=5 10 10 10,clip,width=0.32\textwidth]{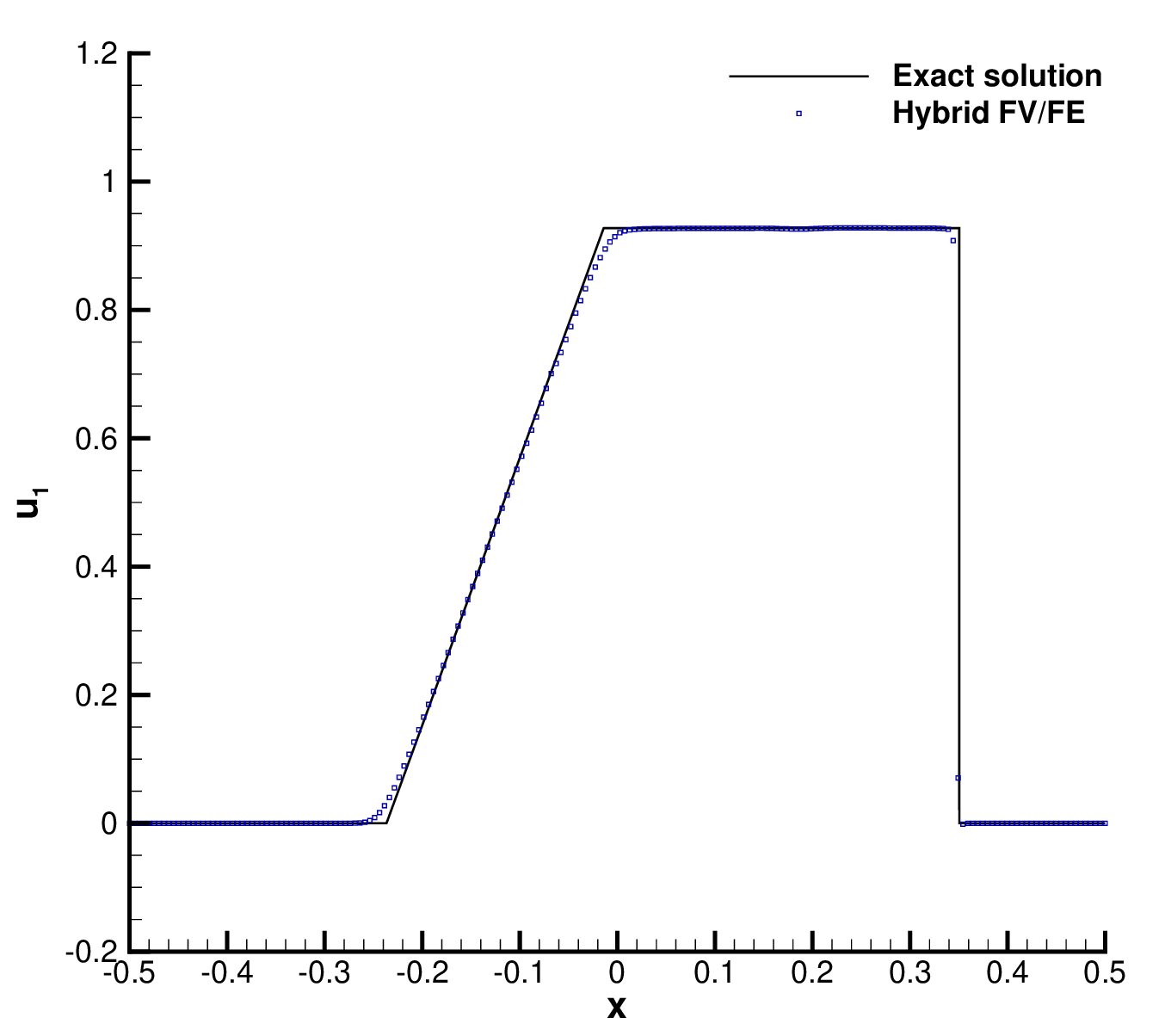}    
		\includegraphics[trim=5 10 10 10,clip,width=0.32\textwidth]{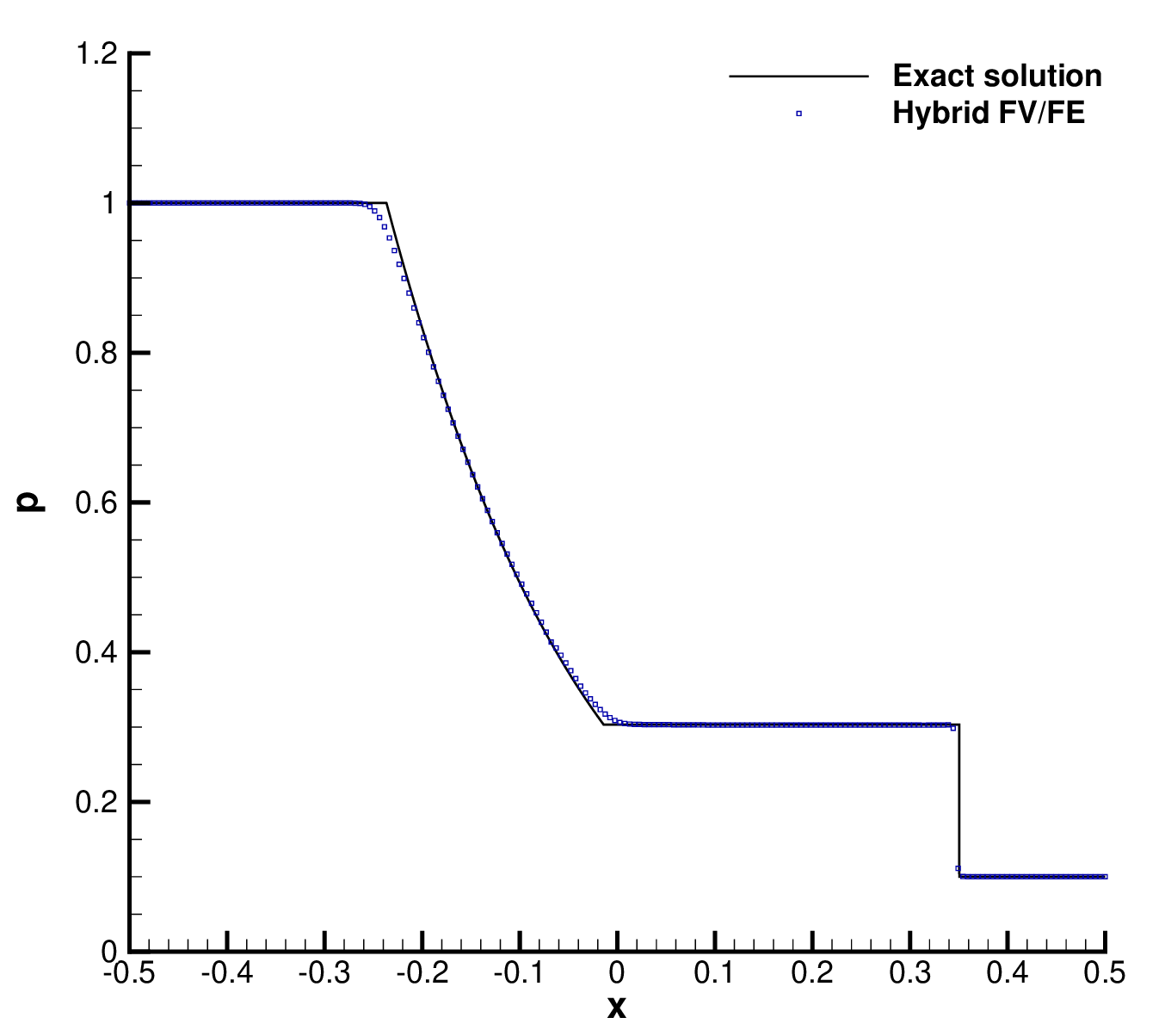}    
		\caption{RP1 Sod. 1D cut in $x-$direction of the numerical solution obtained using the hybrid FV/FE method for the weakly compressible GPR model with the local ADER-BJ approach and auxiliary artificial viscosity $c_{\alpha}=0.2$ (blue squares). Exact solution for the compressible Euler equations (black solid line). From left to right: density, velocity component $u_{1}$, and pressure fields.}  
		\label{fig.RP1sod}
	\end{center}
\end{figure}

\begin{figure}
	\begin{center}
		\includegraphics[trim=5 0 10 10,clip,width=0.32\textwidth]{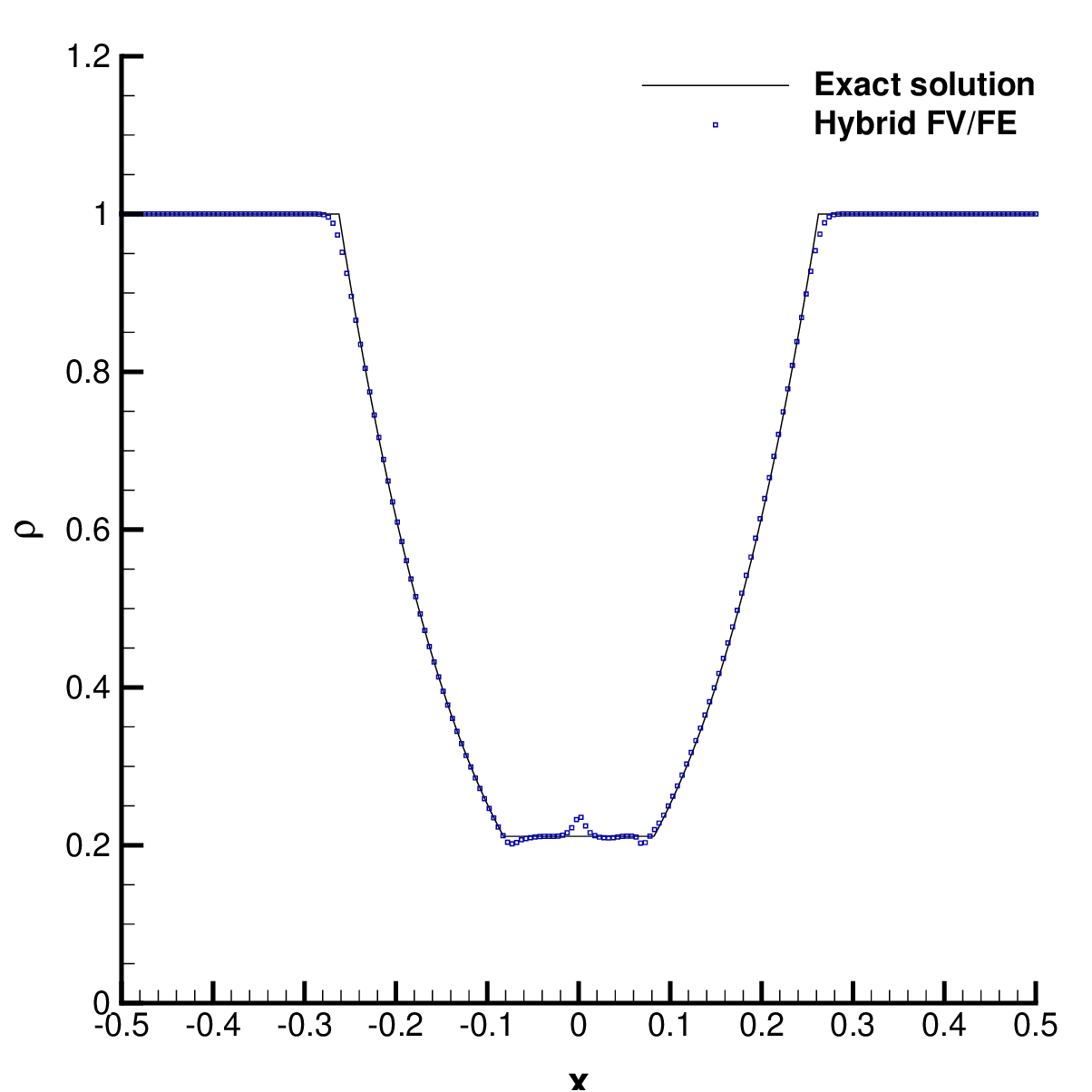}  
		\includegraphics[trim=5 0 10 10,clip,width=0.32\textwidth]{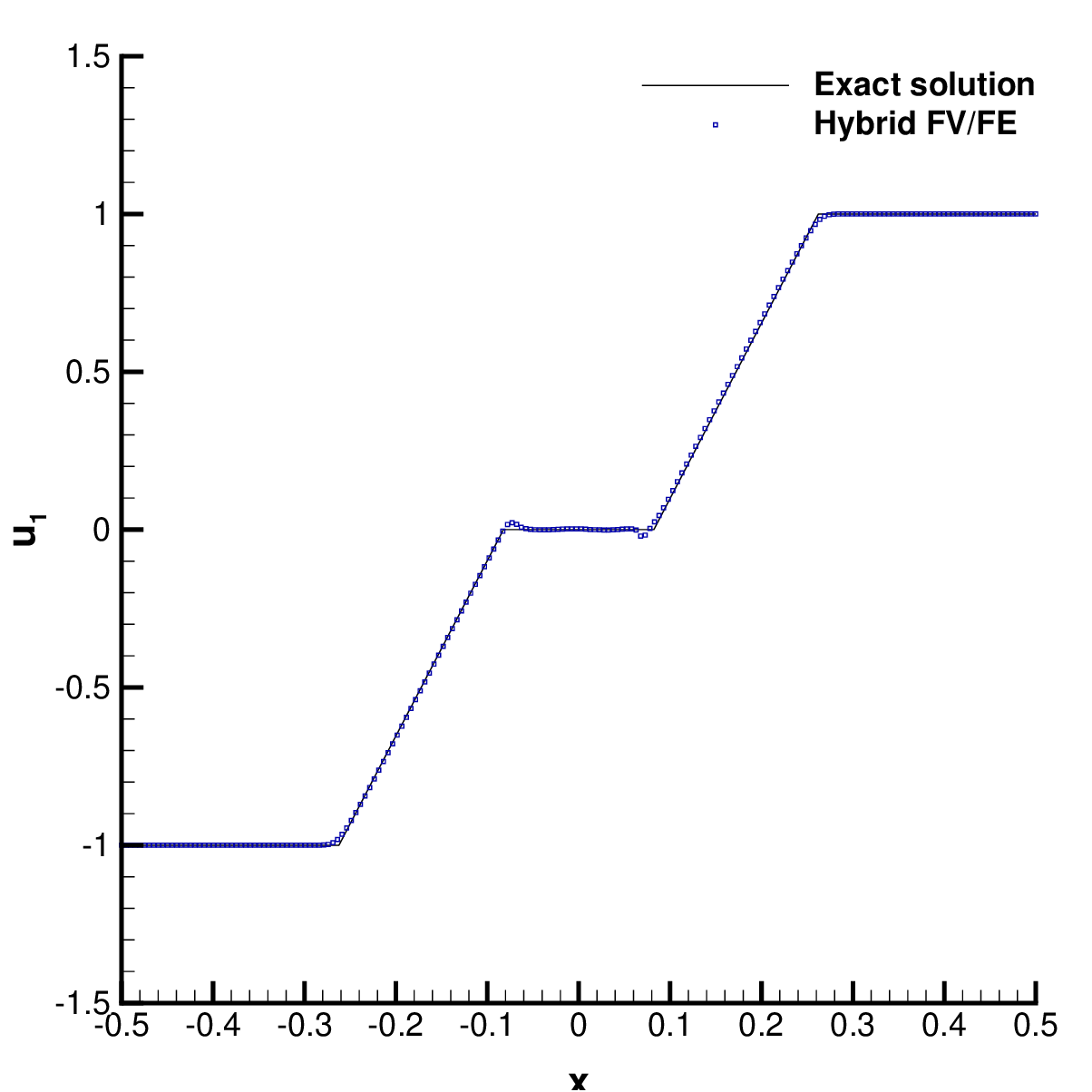}    
		\includegraphics[trim=5 0 10 10,clip,width=0.32\textwidth]{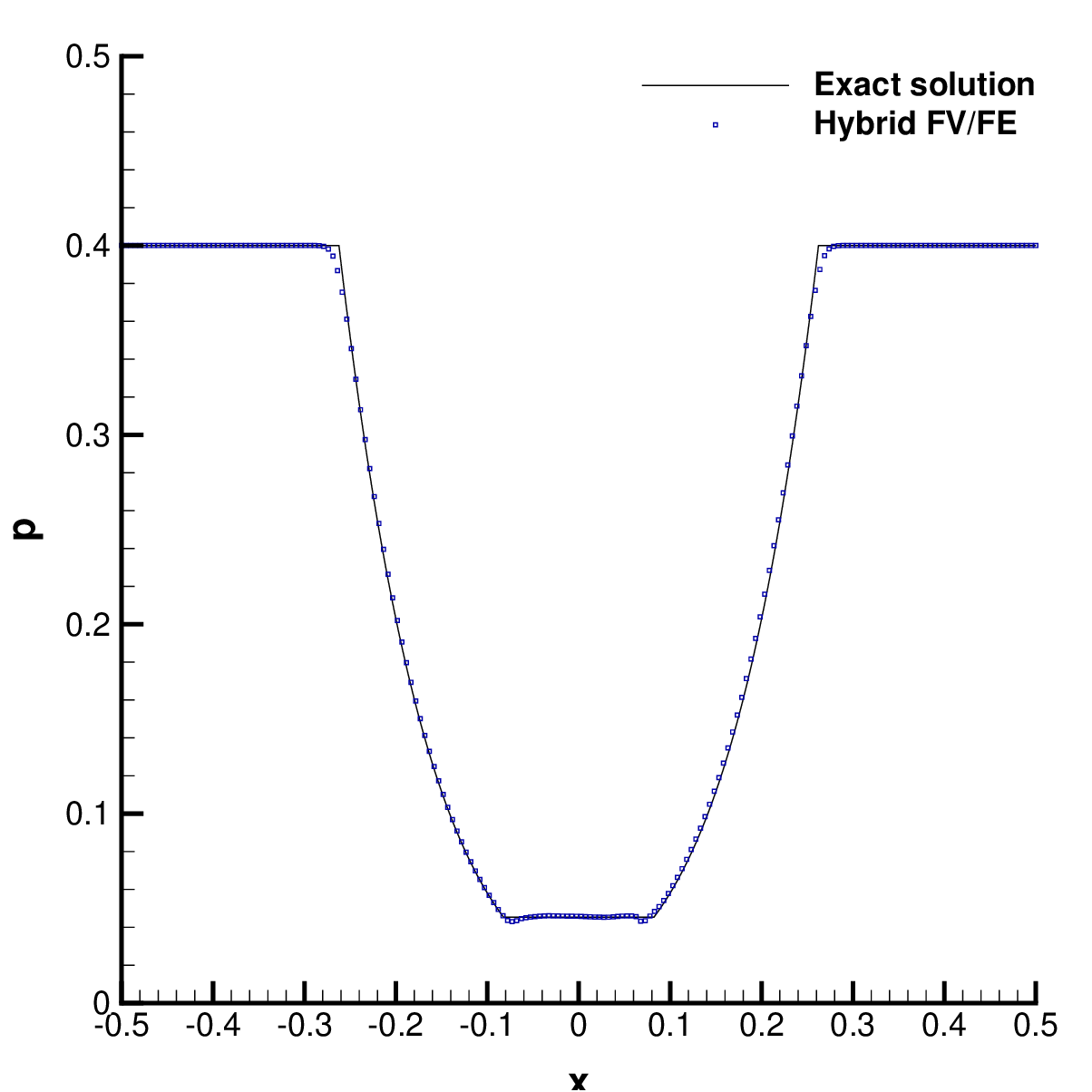}    
		\caption{RP2 double rarefaction. 1D cut in $x-$direction of the numerical solution obtained using the hybrid FV/FE method for the weakly compressible GPR model with the local ADER-BJ approach and auxiliary artificial viscosity $c_{\alpha}=0.1$ (blue squares). Exact solution for the compressible Euler equations (black solid line). From left to right: density, velocity component $u_{1}$, and pressure fields.}  
		\label{fig.RP2}
	\end{center}
\end{figure}

\begin{figure}
	\begin{center}
		\includegraphics[trim=5 10 10 10,clip,width=0.32\textwidth]{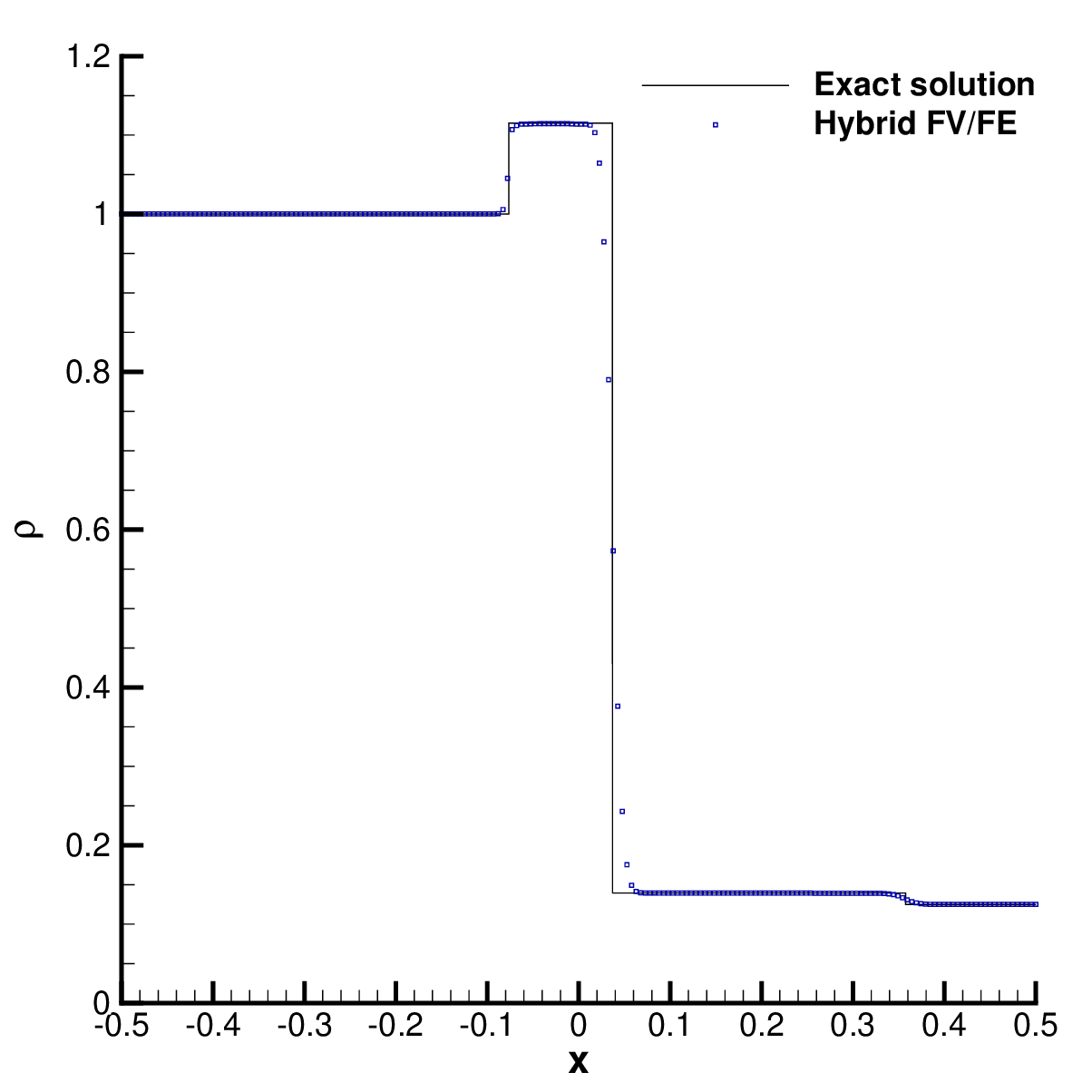}  
		\includegraphics[trim=5 10 10 10,clip,width=0.32\textwidth]{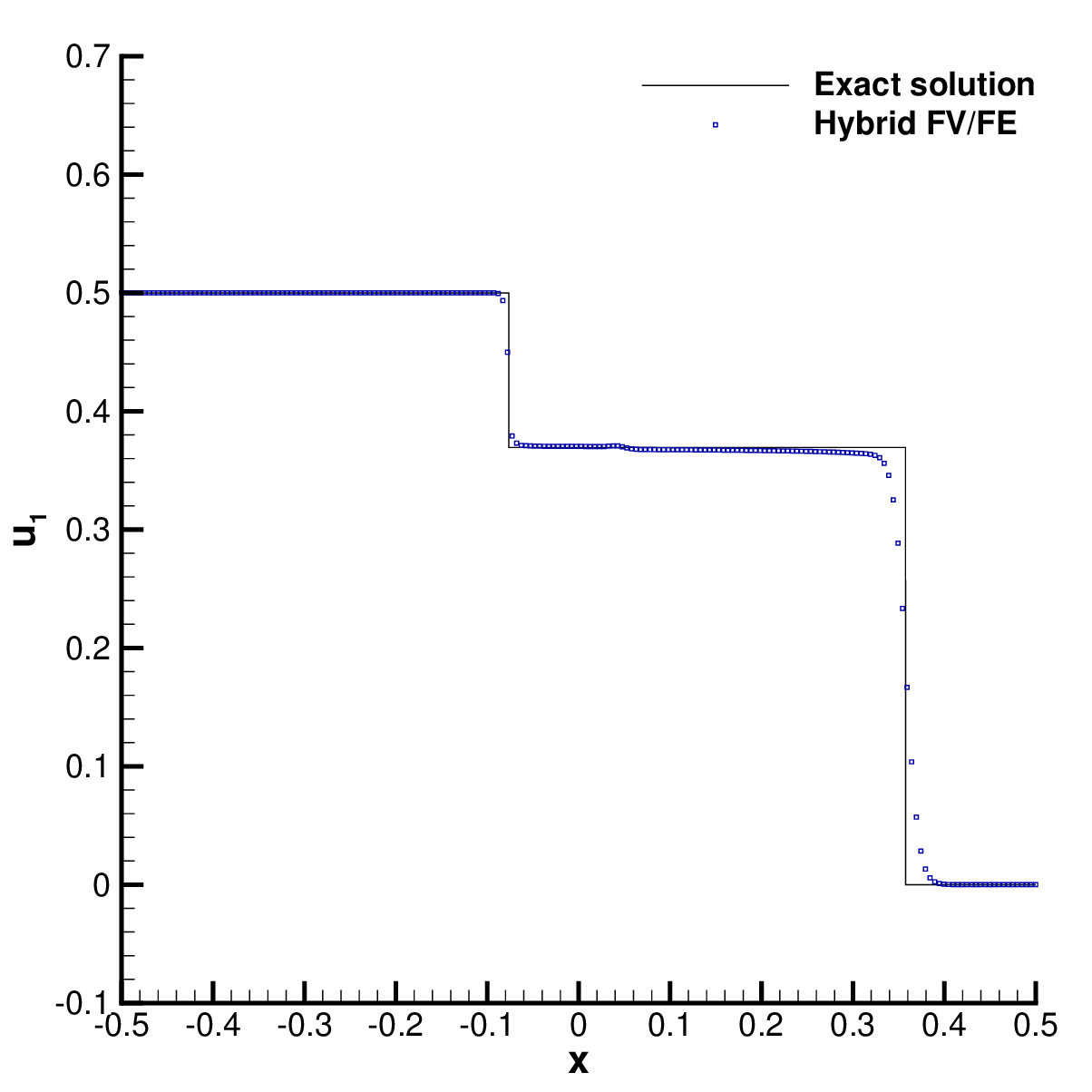}    
		\includegraphics[trim=5 10 10 10,clip,width=0.32\textwidth]{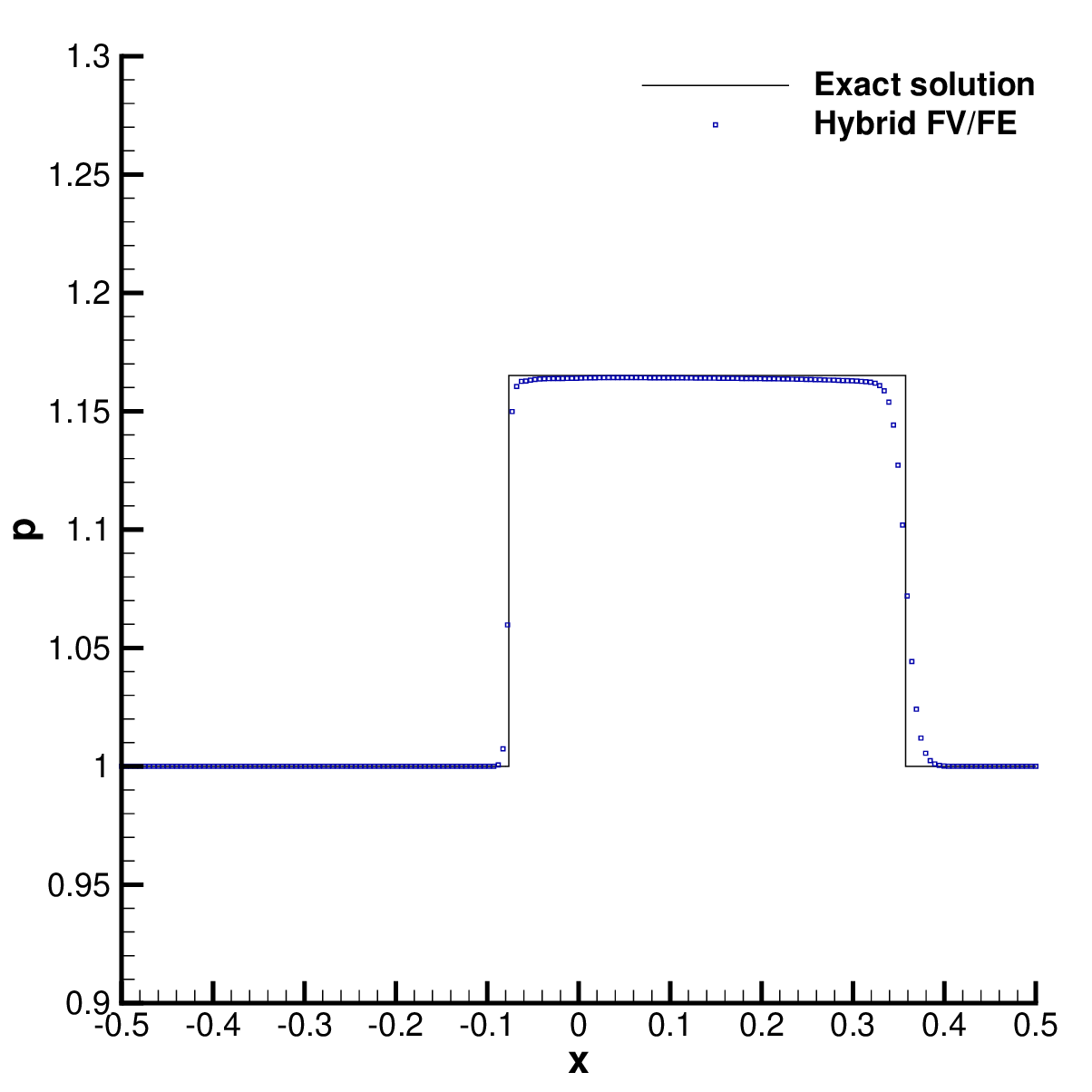}    
		\caption{RP3. 1D cut in $x-$direction of the numerical solution obtained using the hybrid FV/FE method for the weakly compressible GPR model with the local ADER-BJ approach and auxiliary artificial viscosity $c_{\alpha}=0.2$ (blue squares). Exact solution for the compressible Euler equations (black solid line). From left to right: density, velocity component $u_{1}$, and pressure fields.}  
		\label{fig.RP3}
	\end{center}
\end{figure}

Figures~\ref{fig.RP1sod}, \ref{fig.RP2} and \ref{fig.RP3} report the 1D cuts of the density, first component of the velocity vector and pressure along $y=0$ for RP1, RP2 and RP3. A good agreement is observed with the known exact solution of the 1D compressible Euler equations, \cite{Toro}. 
\begin{table}[H]
	\renewcommand{\arraystretch}{1.2}
	\begin{center}
		\begin{tabular}{cccccccccccc}
			Test &  $\rho^{L}$ &  $\rho^{R}$  &  $u_{1}^{L}$ &  $u_{1}^{R}$ &  $u_{2}^{L}$ &  $u_{2}^{R}$ &  $\press^{L}$ &  $\press^{R}$ & $x_{c}$ & $t_{\mathrm{end}}$ & $N_{x}$ \\ \hline
			RP1 & $ 1 $ & $ 0.125 $ & $ 0 $ & $ 0 $ & $ 0 $ & $ 0 $ & $ 1 $ & $ 0.1 $& $ 0 $ &  $ 0.2 $  & $400$\\  
			RP2 & $ 1 $ & $ 1 $ & $ -1 $ & $ 1 $& $ 0 $ & $ 0 $  & $ 0.4 $ & $ 0.4 $& $ 0 $ &  $ 0.15$  & $400$\\
			RP3 & $ 1 $ & $ 0.125 $ & $ 0.5 $ & $ 0 $ & $ 0 $ & $ 0 $ & $ 1 $ & $ 1 $& $ 0 $ &  $ 0.1 $  & $400$\\
			RP4 & $ 1 $ & $ 0.5 $ & $ 0 $ & $ 0 $ & $ -0.2 $ & $ 0.2 $ & $ 1 $ & $ 0.5 $& $ 0 $ &  $ 0.2 $  & $400$\\
			RP5 &  $ 1 $ & $ 0.5 $ & $ 0 $ & $ 0 $ & $ -0.2 $ & $ 0.2 $ & $ 1 $ & $ 0.5 $& $ 0 $ &  $ 0.2 $  & $400$\\
			RP6 &  $ 1 $ & $ 0.5 $ & $ 0 $ & $ 0 $ & $ -0.2 $ & $ 0.2 $ & $ 1 $ & $ 0.5 $& $ 0 $ &  $ 0.2 $  & $400$\\
		\end{tabular} 
	\end{center}
	\caption{Riemann problems. Initial condition, location of the initial discontinuity, $x_{c}$, final time, $t_{\mathrm{end}}$, and number of mesh divisions on $x$-direction, $N_{x}$.}
	\label{tab.RP_IC}
\end{table}

\begin{figure}
	\begin{center}
		\includegraphics[trim=5 10 10 10,clip,width=0.32\textwidth]{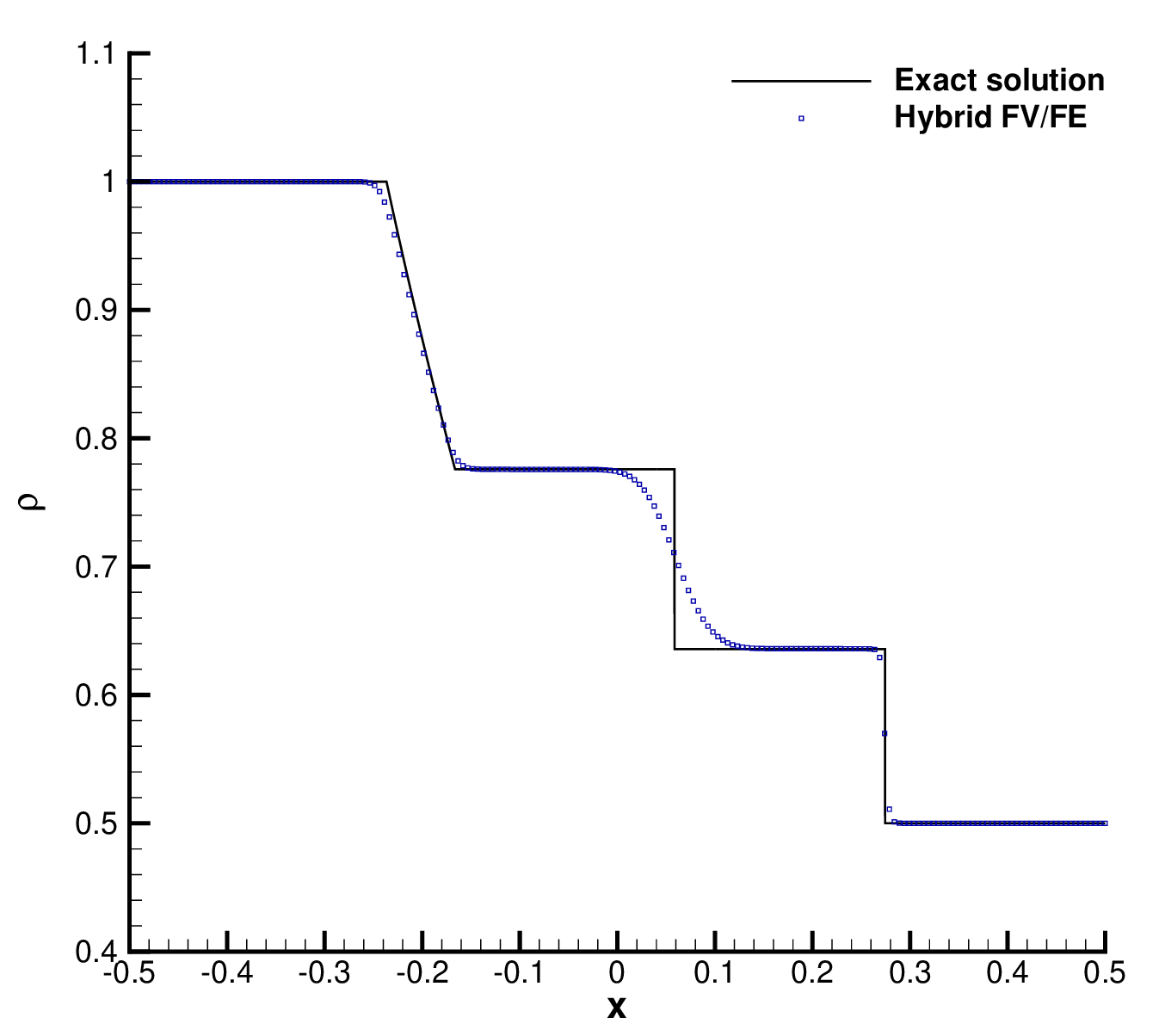}  
		\includegraphics[trim=5 10 10 10,clip,width=0.32\textwidth]{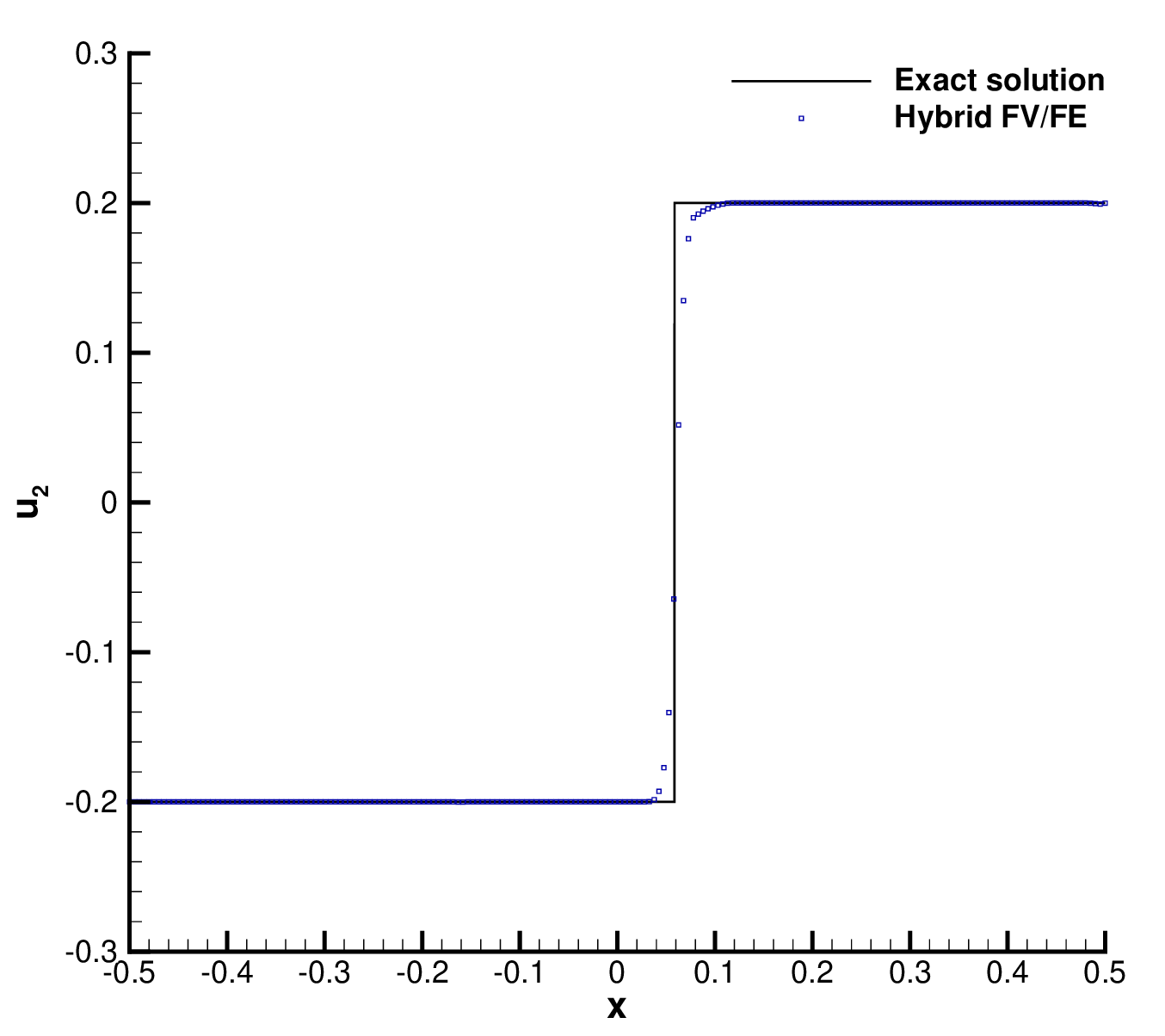}    
		\includegraphics[trim=5 10 10 10,clip,width=0.32\textwidth]{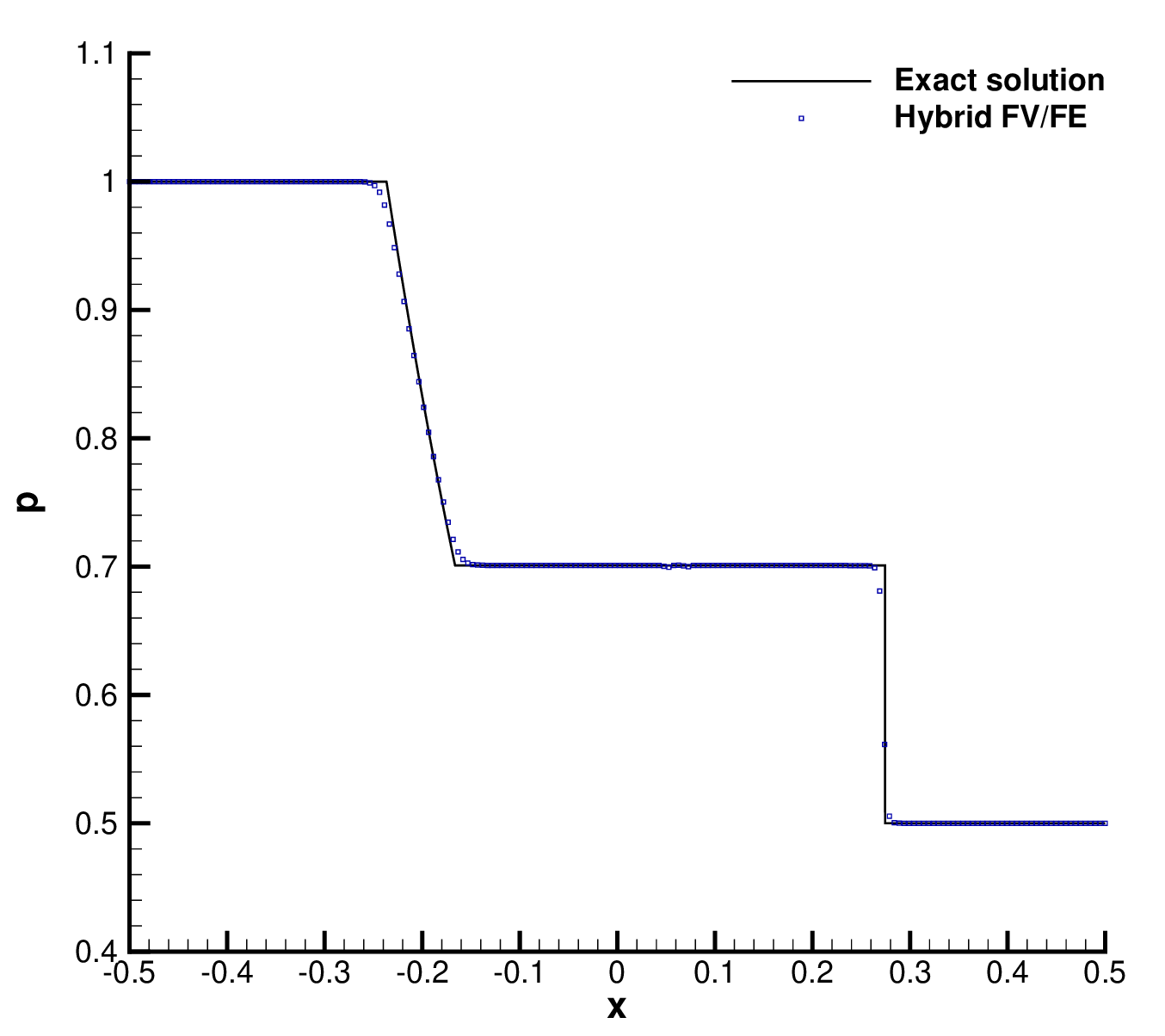}    
		\caption{RP4. 1D cut in $x-$direction of the numerical solution obtained using the hybrid FV/FE method for the weakly compressible GPR model with the local ADER-BJ approach and auxiliary artificial viscosity $c_{\alpha}=1$ (blue squares). Reference solution computed with a TVD-FV scheme on a mesh of $128000$ cells (black solid line). From left to right: density, velocity component $u_{2}$, and pressure fields.}  
		\label{fig.RP4}
	\end{center}
\end{figure}

\begin{figure}
\begin{center}
	\includegraphics[clip,width=0.45\textwidth]{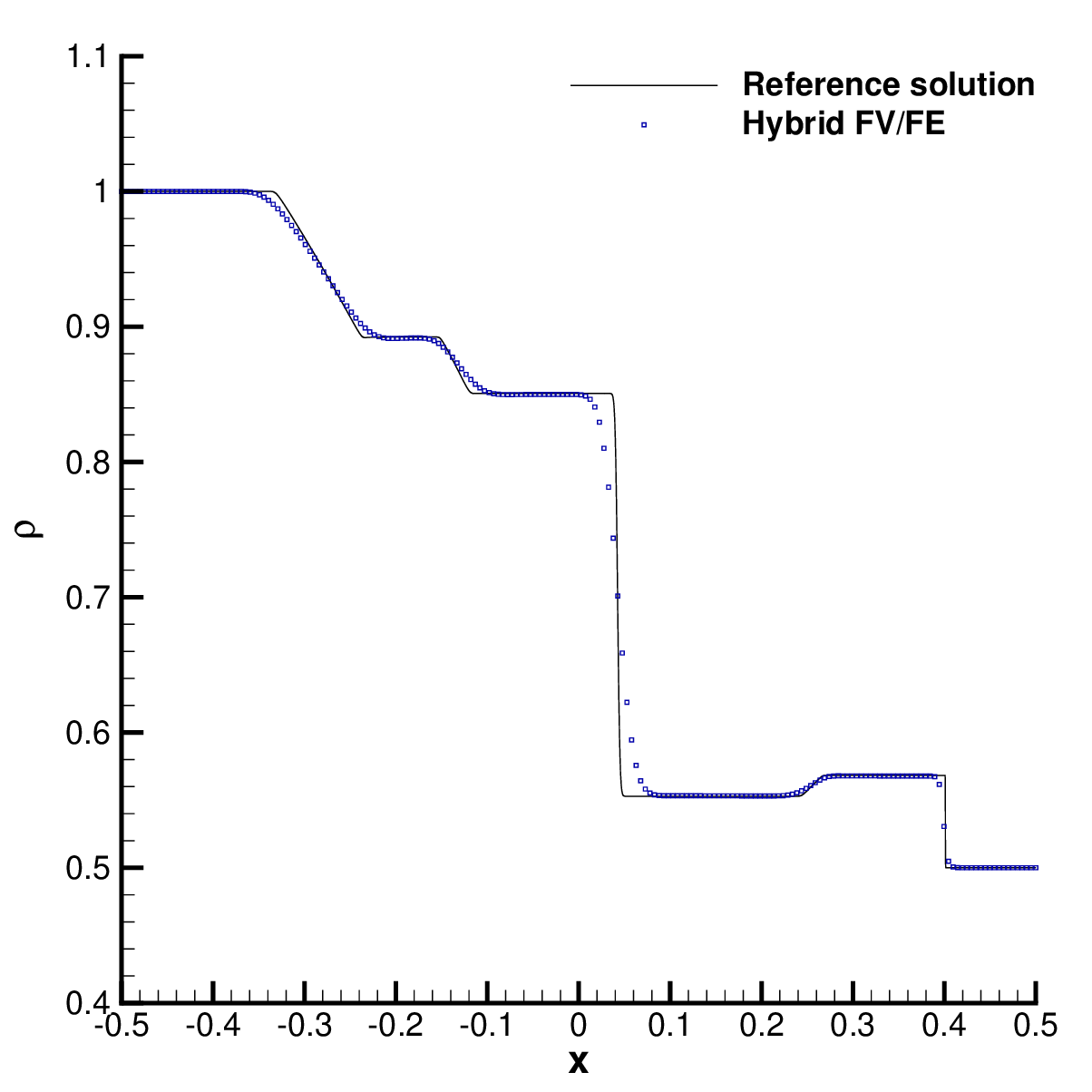}  
	\includegraphics[width=0.45\textwidth]{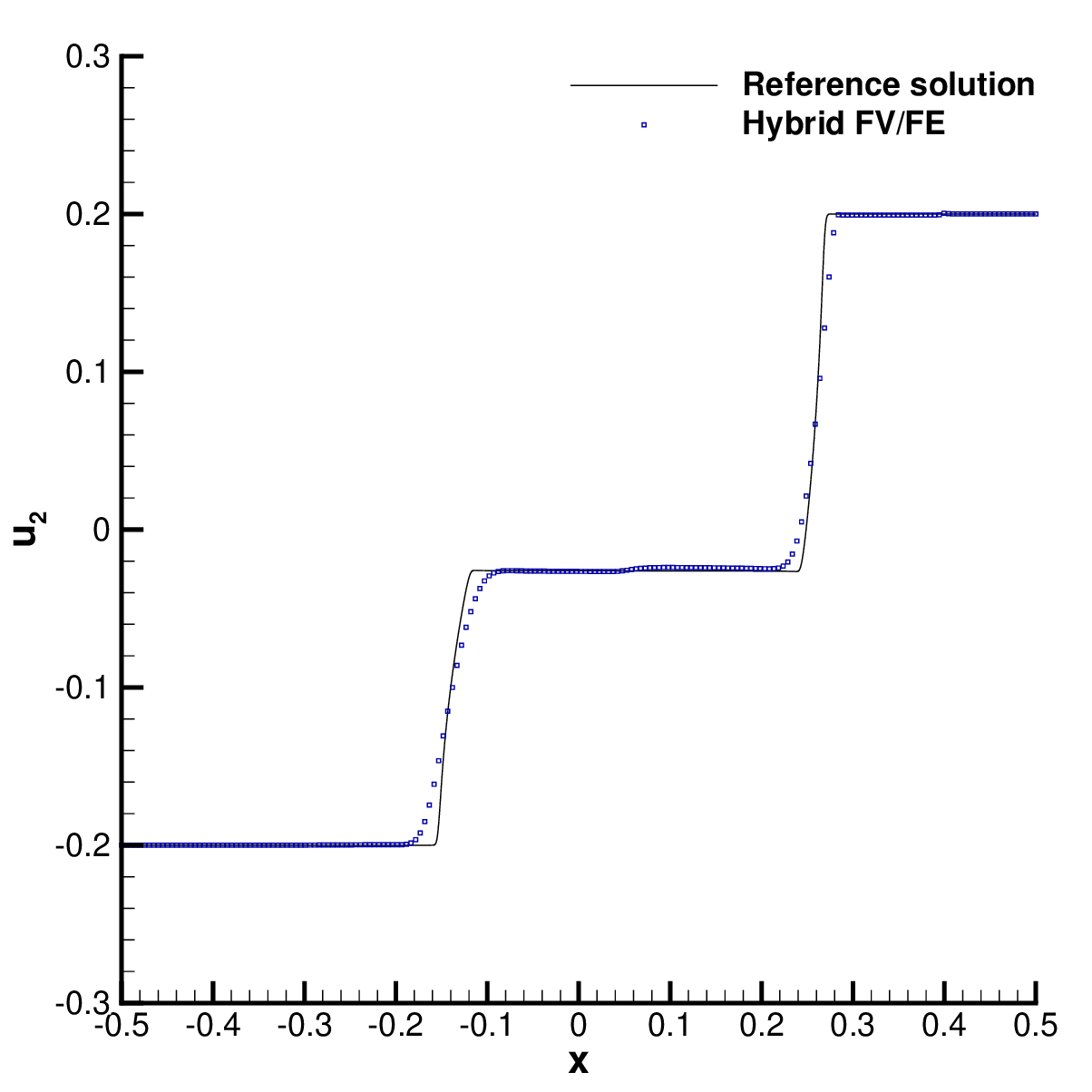}    \\ 
	\includegraphics[width=0.45\textwidth]{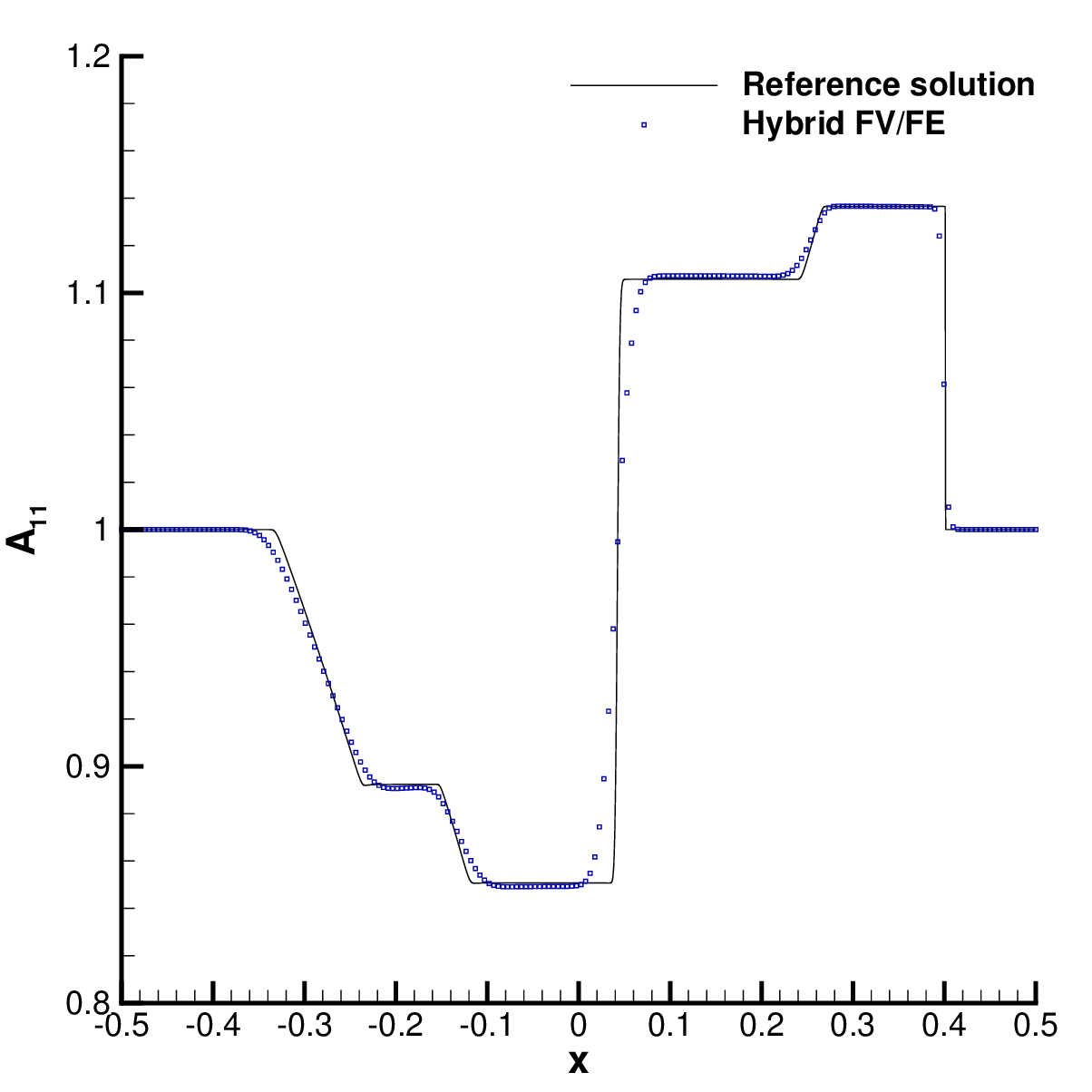}   
	\includegraphics[width=0.45\textwidth]{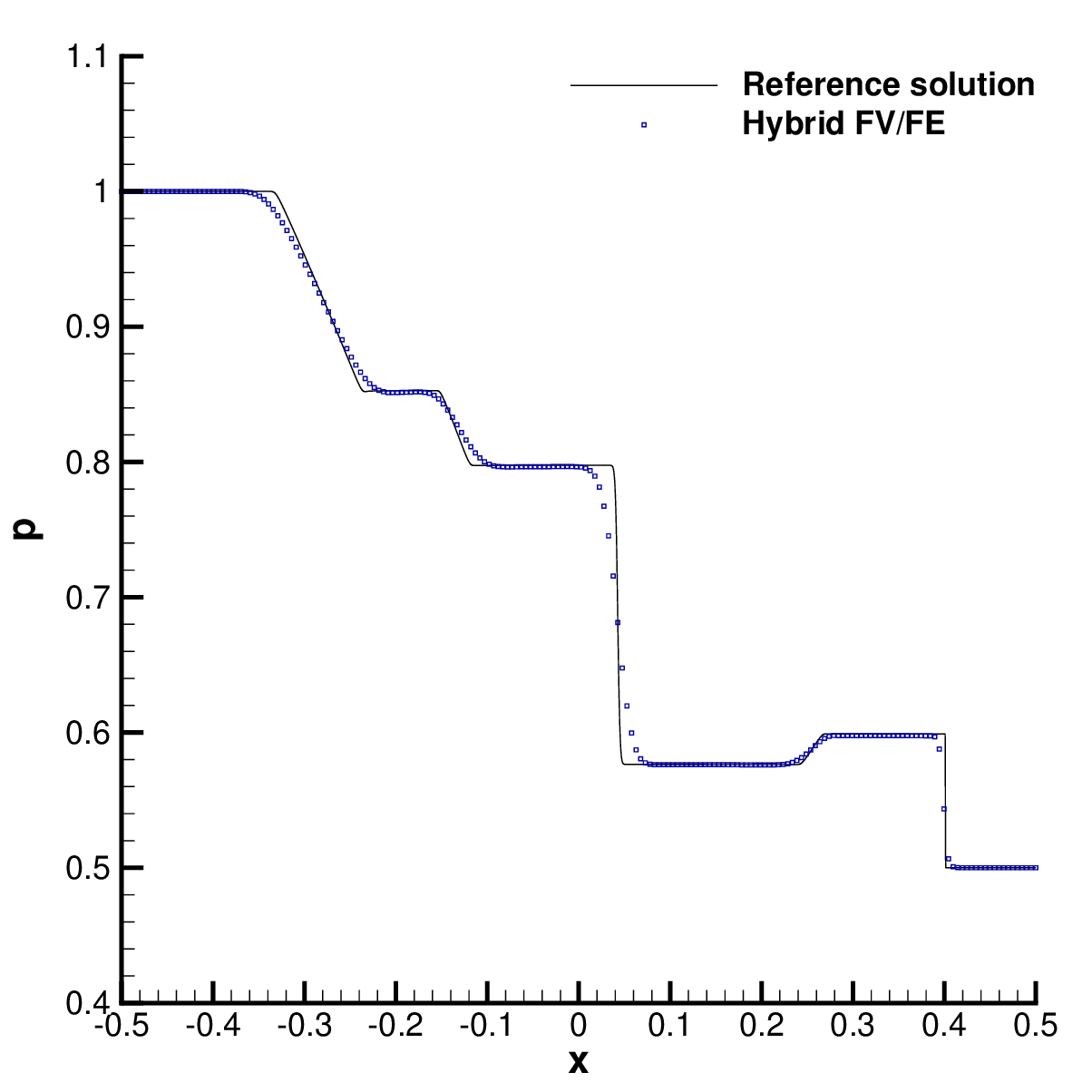}    \\ 
	\caption{RP5. 1D cut in $x-$direction of the numerical solution obtained using the hybrid FV/FE method for the weakly compressible GPR model with the local ADER-MM approach and auxiliary artificial viscosity $c_{\alpha}=1$ (blue squares). Reference solution computed with a HTC-FV scheme on a mesh of $25000$ cells (black solid line). From left top to right bottom: density, velocity component $u_{2}$, distortion field component $A_{11}$, and pressure fields.}  
	\label{fig.RP5ch0}
\end{center}
\end{figure}

\begin{figure}[H]
	\begin{center}
		\includegraphics[width=0.45\textwidth]{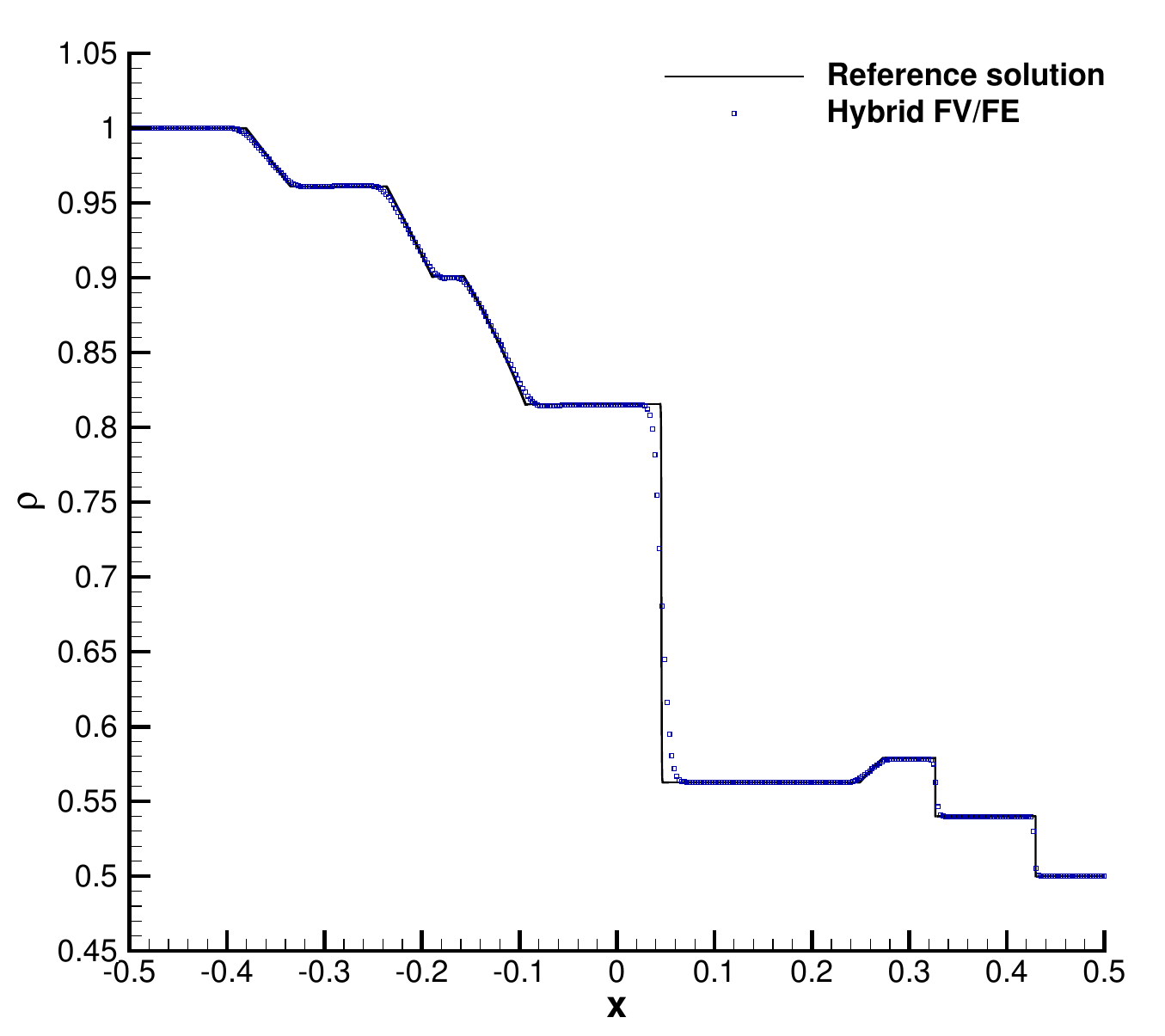}  
		\includegraphics[width=0.45\textwidth]{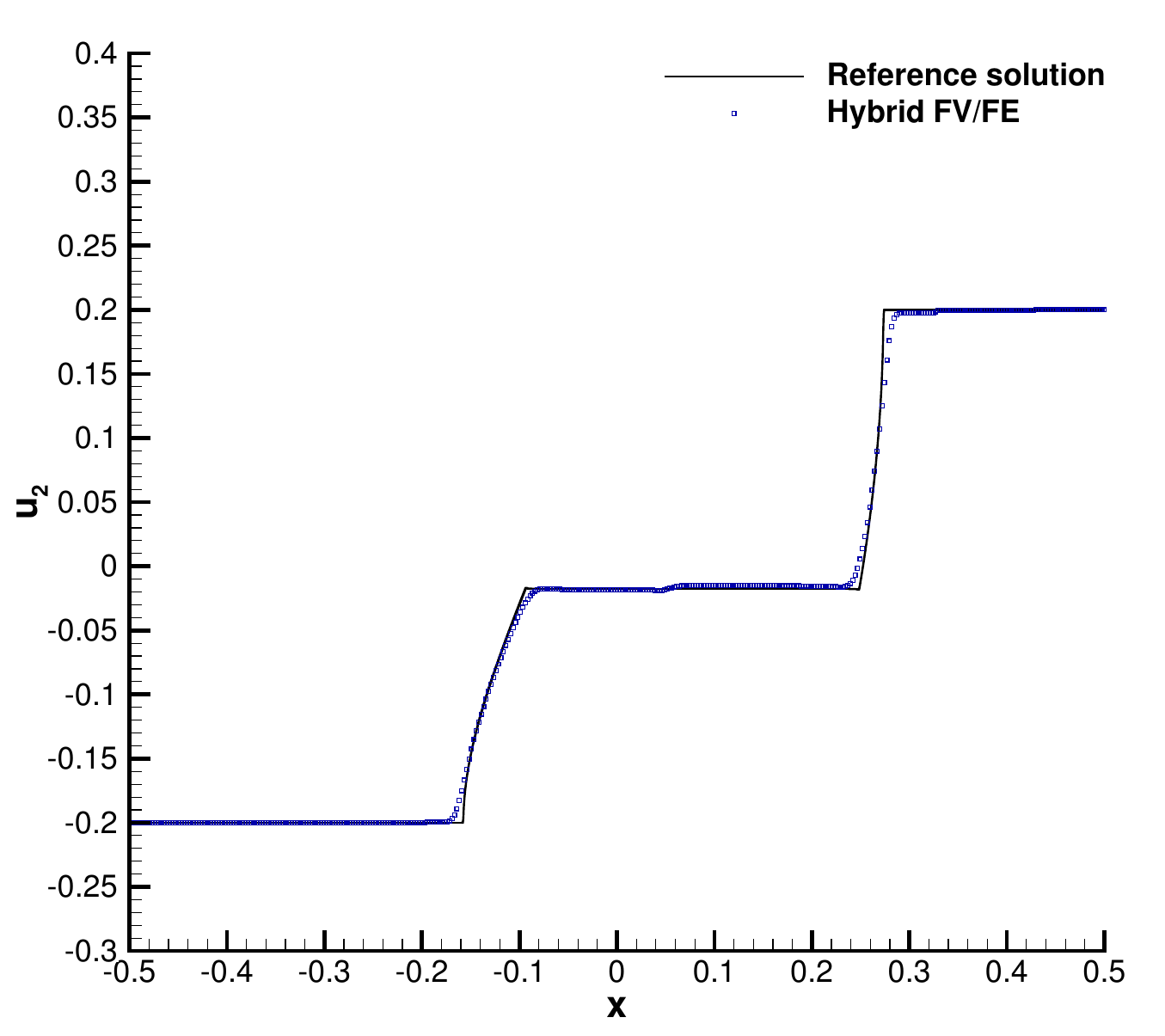}    \\ 
		\includegraphics[width=0.45\textwidth]{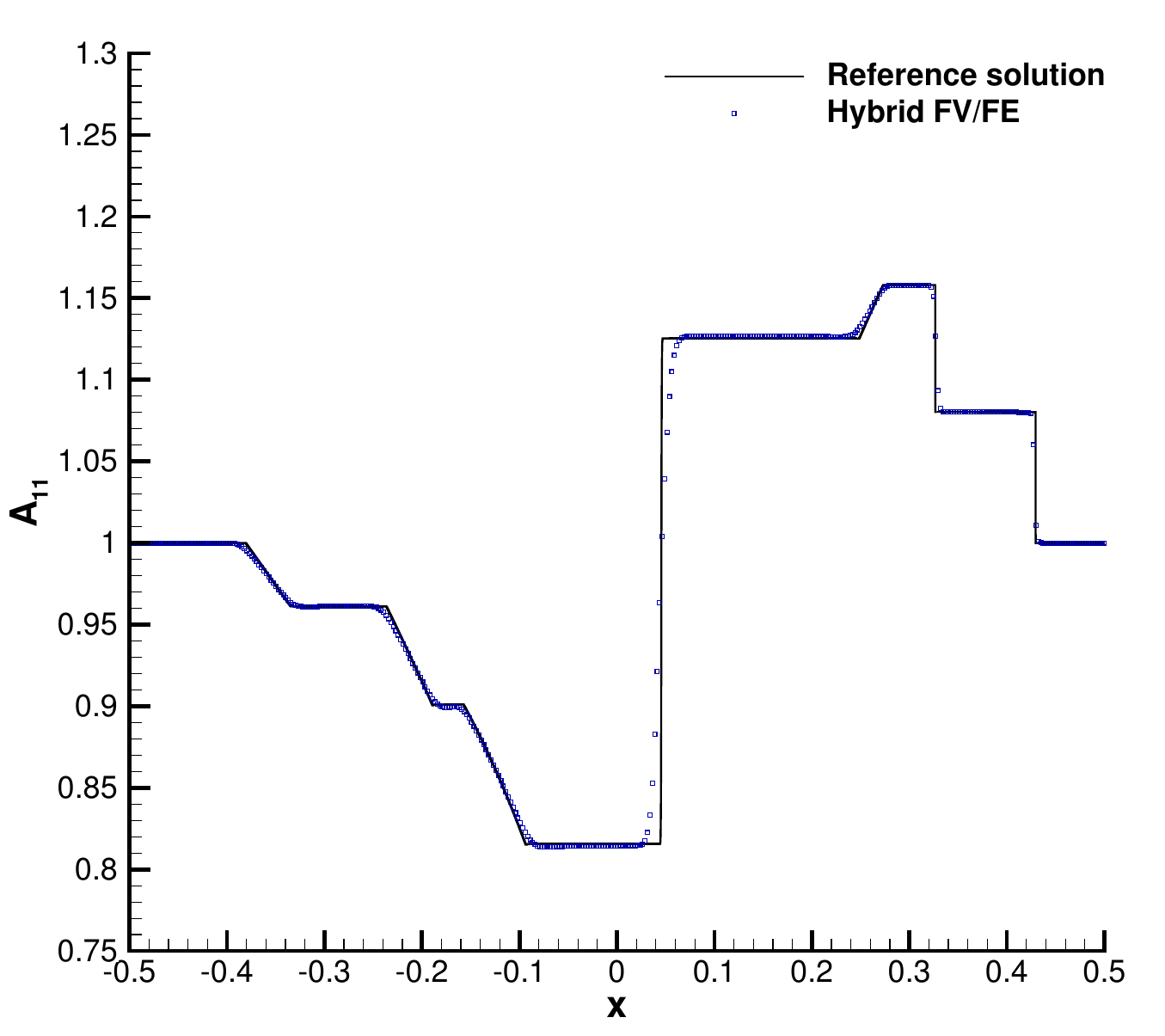}   
		\includegraphics[width=0.45\textwidth]{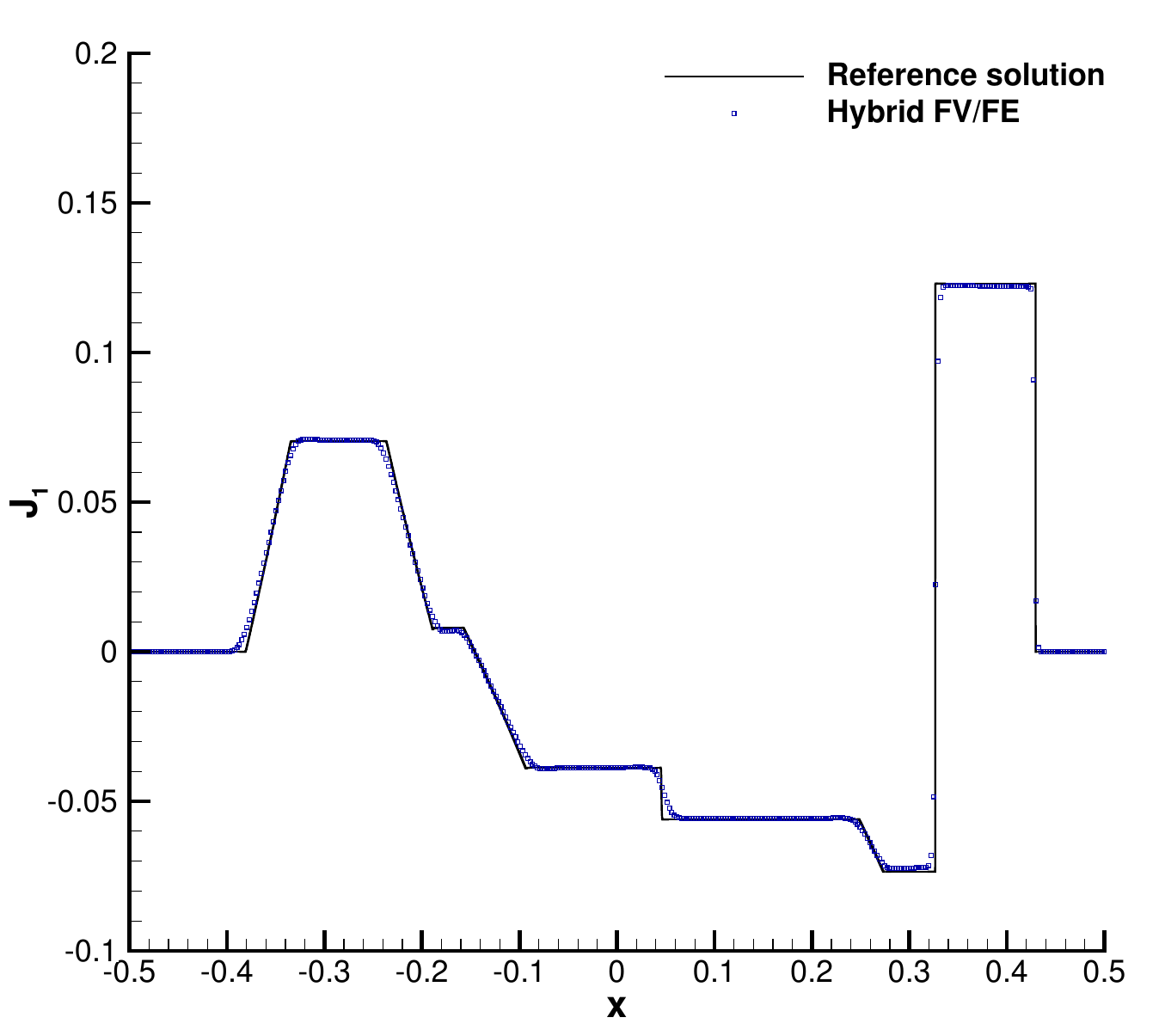}    \\ 
		\caption{RP6. 1D cut in $x-$direction of the numerical solution obtained using the hybrid FV/FE method for the weakly compressible GPR model with the local ADER-MM approach and auxiliary artificial viscosity $c_{\alpha}=1$ (blue squares). Reference solution computed with a HTC-FV scheme on a mesh of $25000$ cells (black solid line). From left top to right bottom: density, velocity component $u_{2}$, distortion field component $A_{11}$, and heat flux component $J_{1}$ fields.}  
		\label{fig.RP5}
	\end{center}
\end{figure}

RP4, RP5 and RP6 are tests specifically designed to assess the complete GPR model, and are characterised by having the same initial conditions but consider different types of materials, \cite{Boscheri2021SIGPR}. As observed in Figure~\ref{fig.RP4}, the ideal fluid studied in RP4 leads to one contact discontinuity, one shear wave and two acoustic waves. The obtained results agree well with the solution computed employing a second order TVD finite volume scheme on a 1D mesh formed by 128000 cells. 
On the other hand, for an ideal elastic solid without heat conduction, we obtain two acoustic waves (a left rarefaction and a right shock), two shear waves (one left and one right going) and one contact discontinuity. 
The obtained results are compared in Figure~\ref{fig.RP5ch0} against a reference numerical solution computed in a 1D grid of $25000$ cells with a second order finite volume thermodynamically compatible scheme for the entropy-based formulation of the GPR model, i.e., the scheme solves the entropy equation instead of the total energy, see \cite{HTCTotalAbgrall} for more details. 
Finally, Figure~\ref{fig.RP5} shows the results obtained for RP5, where the effect of the heat flux is taken into account and, as a consequence, a couple of new left and right thermo-acoustic waves arise. Also in this case the TVD-FV scheme for the compressible GPR model is employed to provide a reference solution.

\subsection{2D circular explosions}
We now address two circular explosion problems, one in the fluid framework and the other one in the solid limit of the weakly compressible GPR model.
\begin{figure}[H]
	\begin{center}
		\includegraphics[width=0.45\textwidth]{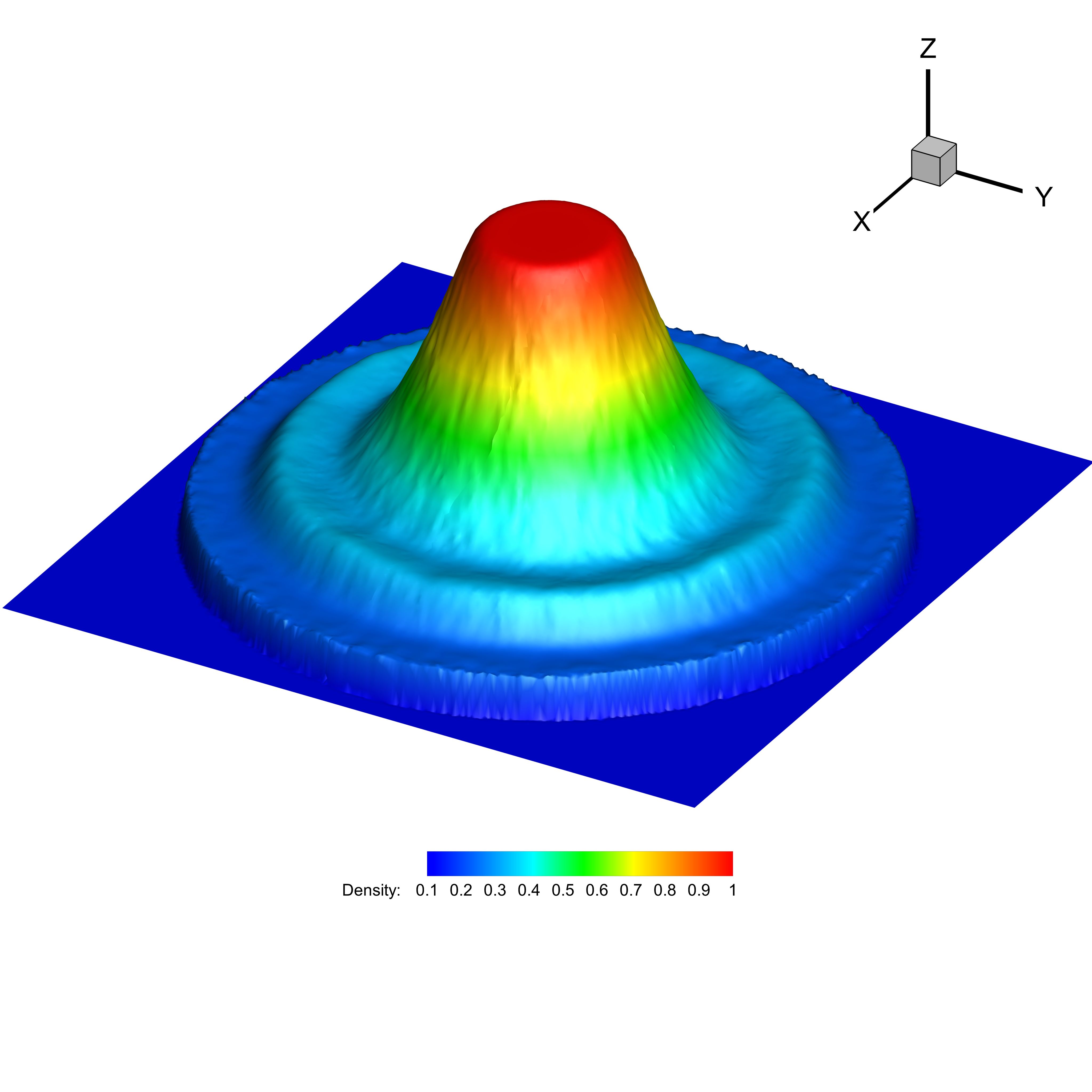}  
		\includegraphics[width=0.45\textwidth]{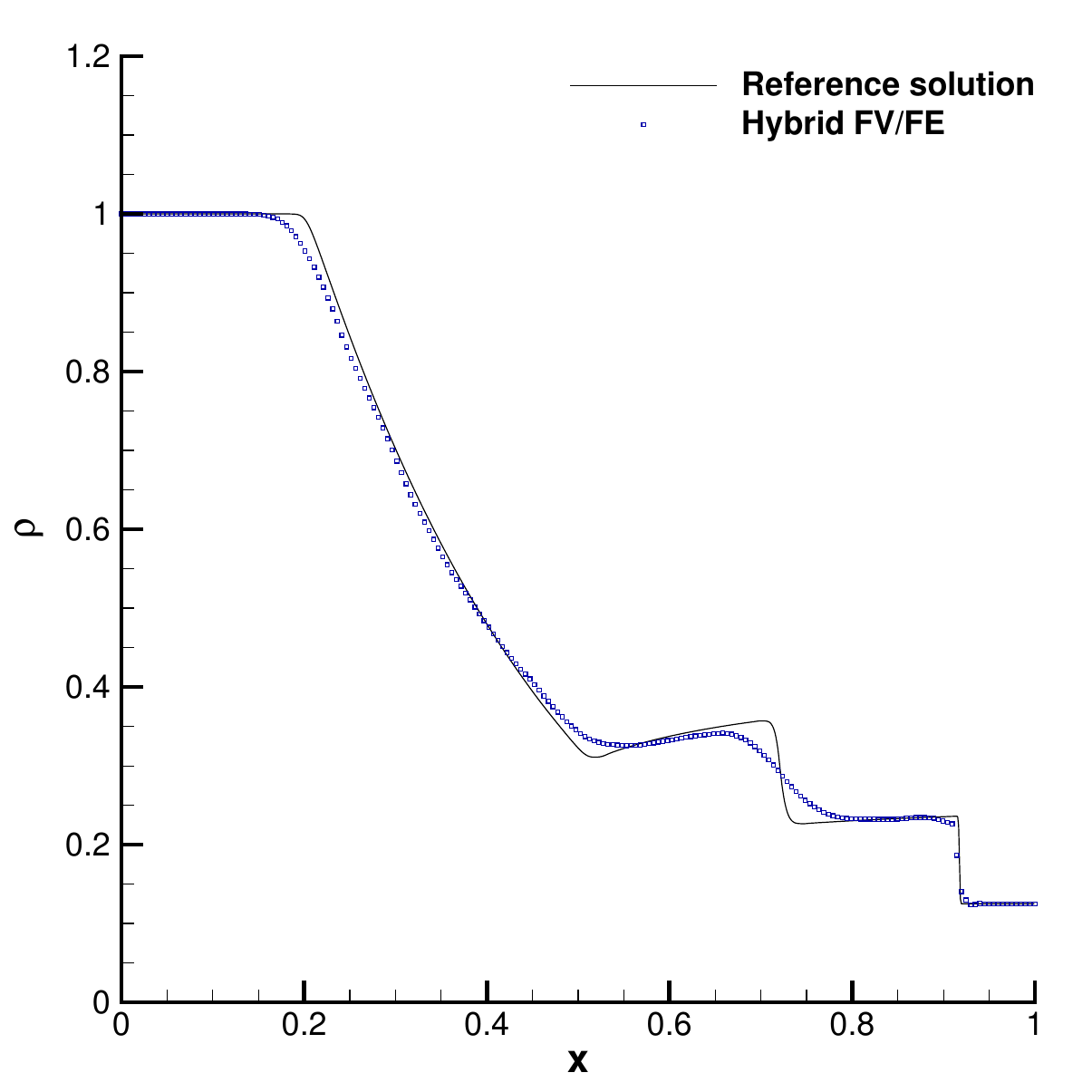}    \\ 
		\includegraphics[width=0.45\textwidth]{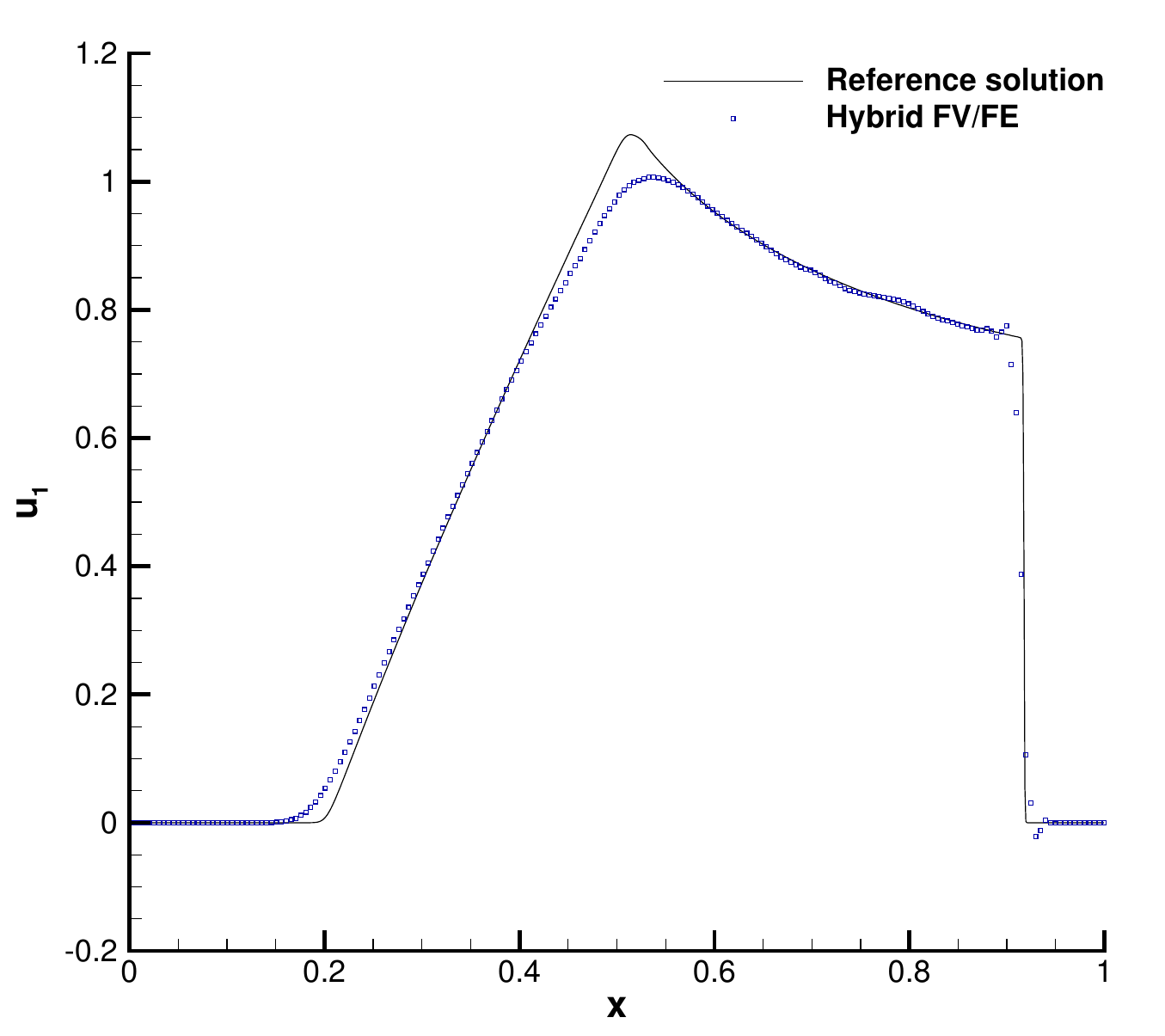}   
		\includegraphics[width=0.45\textwidth]{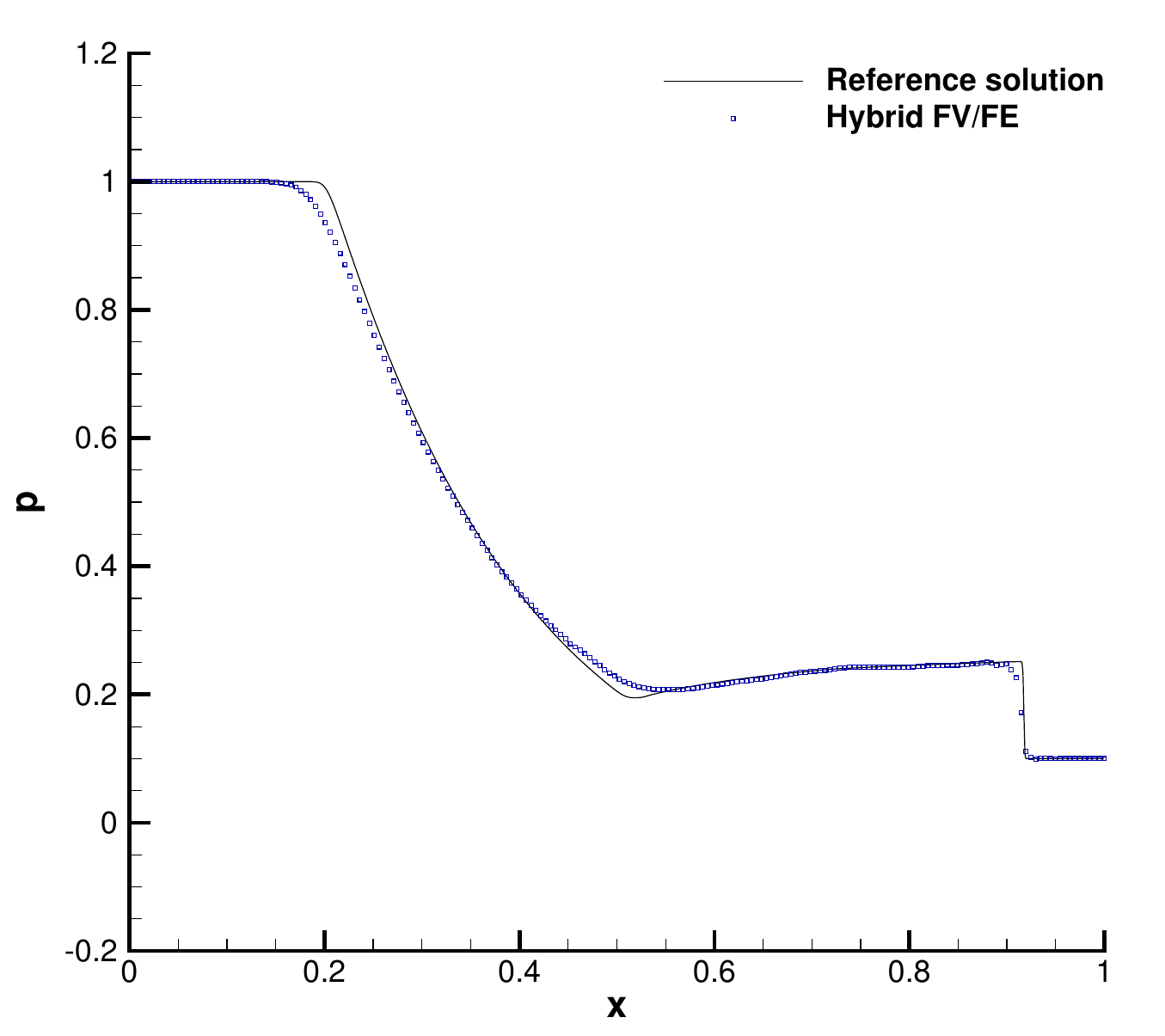}    \\ 
		\caption{Fluid circular explosion. Left top: elevated contour plot of the density field. From right top to right bottom: 1D cut in $x-$direction of the density, velocity component $\vel_{1}$ and pressure obtained using the hybrid FV/FE method for the weakly compressible GPR model with the local ADER-ENO approach and auxiliary artificial viscosity $c_{\alpha}=0.5$ (blue squares). Reference solution (solid black line). }  
		\label{fig.CEfluid}
	\end{center}
\end{figure}
\subsubsection{Fluid circular explosion}
First, we consider the classical circular explosion problem based on the extension to radial flows of the 1D Sod shock tube benchmark, \cite{Lagrange2D}, and whose initial condition is given by
\begin{gather*}
	\rho\left(\x,0\right) =  \left\lbrace \begin{array}{lr}
		1 & \mathrm{ if } \; r \le 0.5,\\
		0.125 & \mathrm{ if } \; r > 0.5,
	\end{array}\right. \qquad
	\mathbf{u}\left(\x,0\right) =  \boldsymbol{0}, \qquad
	\press\left(\x,0\right) = \left\lbrace \begin{array}{lr}
		1 & \mathrm{ if } \; r \le 0.5,\\
		0.1 & \mathrm{ if } \; r > 0.5,
	\end{array}\right. \notag\\
	\mathbf{A}\left(\x,0\right)=\mathbf{I}, \qquad \mathbf{J}\left(\x,0\right) = \boldsymbol{0},
	\qquad r = \sqrt{x^2+y^2}.
\end{gather*}
The parameters of the GPR model are set to $c_{s}=c_{h}=0$ and $\mu =\kappa = 0$. 
The computational domain $\Omega=[-1,1]^2$ is discretised with a primal mesh of $85344$ triangles, and periodic boundary conditions are imposed everywhere. The numerical results obtained with the new hybrid FV/FE scheme using the local ADER-ENO approach are reported in Figure~\ref{fig.CEfluid}. The obtained solution shows a good agreement with the reference solution computed employing the 1D partial differential equation in radial direction with geometrical source terms equivalent to the compressible Euler system and solved using a second order TVD-FV scheme on a grid of $10^4$ cells, \cite{DT11}.

\subsubsection{Solid circular explosion}
\begin{figure}[H]
	\begin{center}
		\includegraphics[clip,trim=0 200 0 0,width=0.45\textwidth]{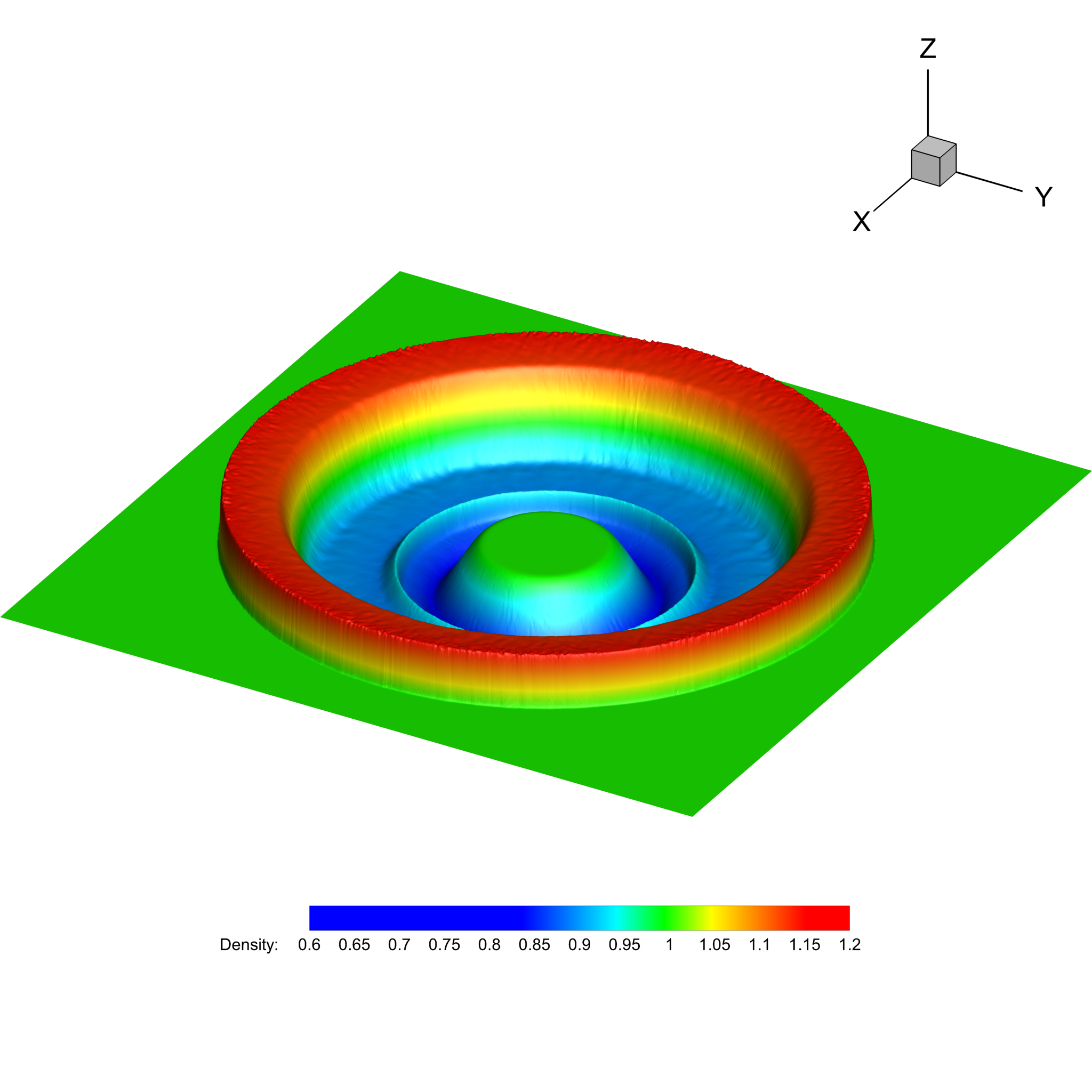}  
		\includegraphics[clip,trim=0 200 0 0,width=0.45\textwidth]{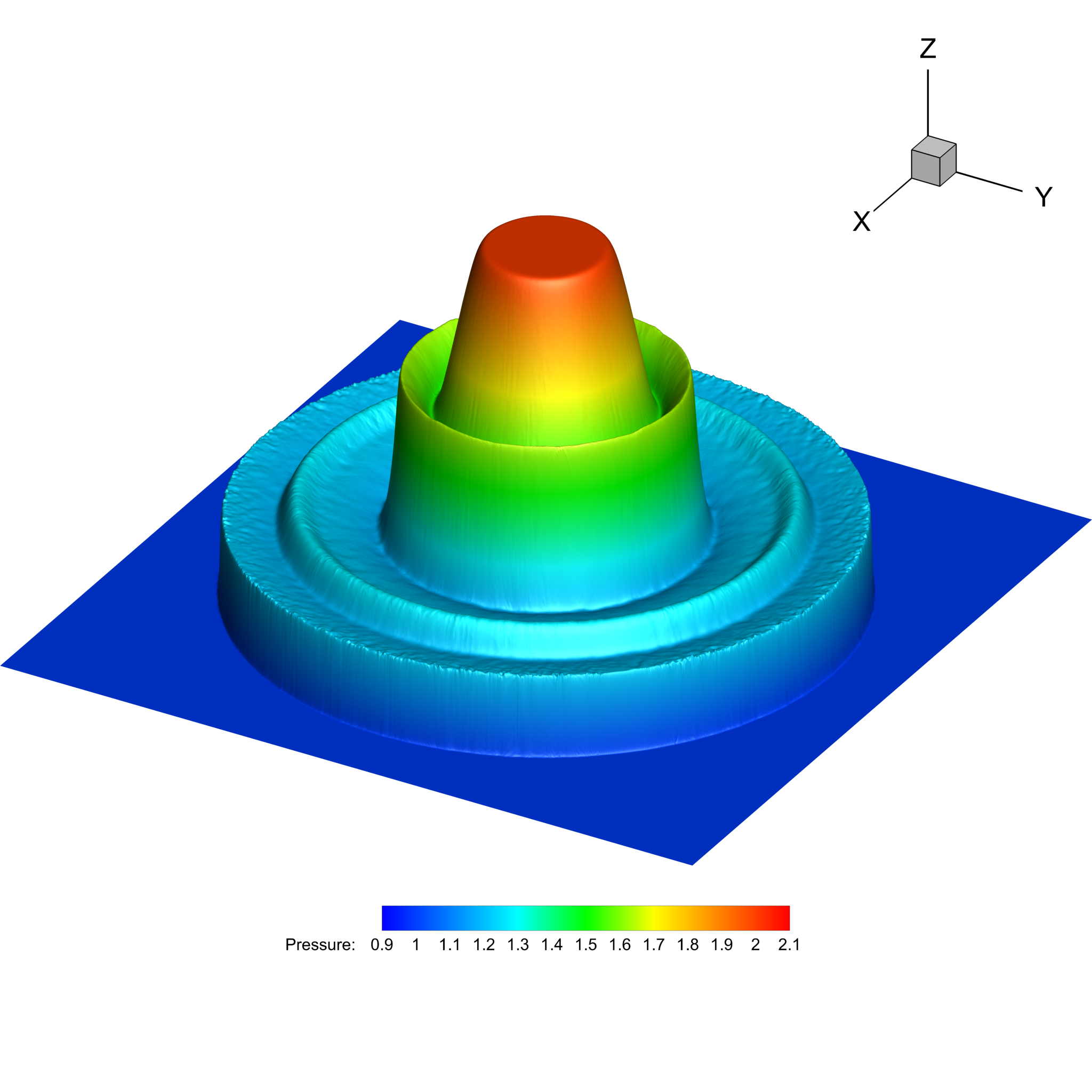}    \\ 
		\includegraphics[clip,trim=0 200 0 0,width=0.45\textwidth]{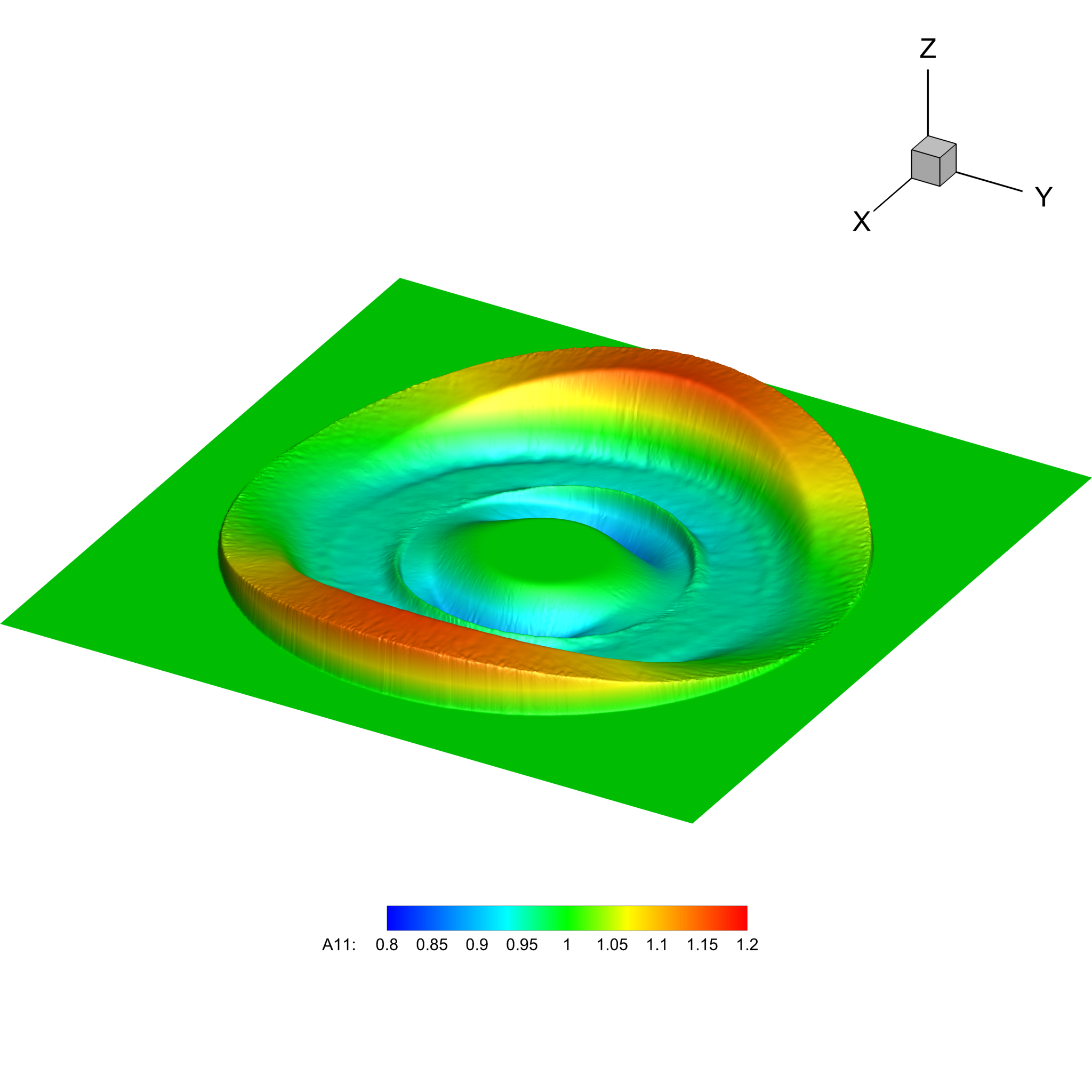}   
		\includegraphics[clip,trim=0 200 0 0,width=0.45\textwidth]{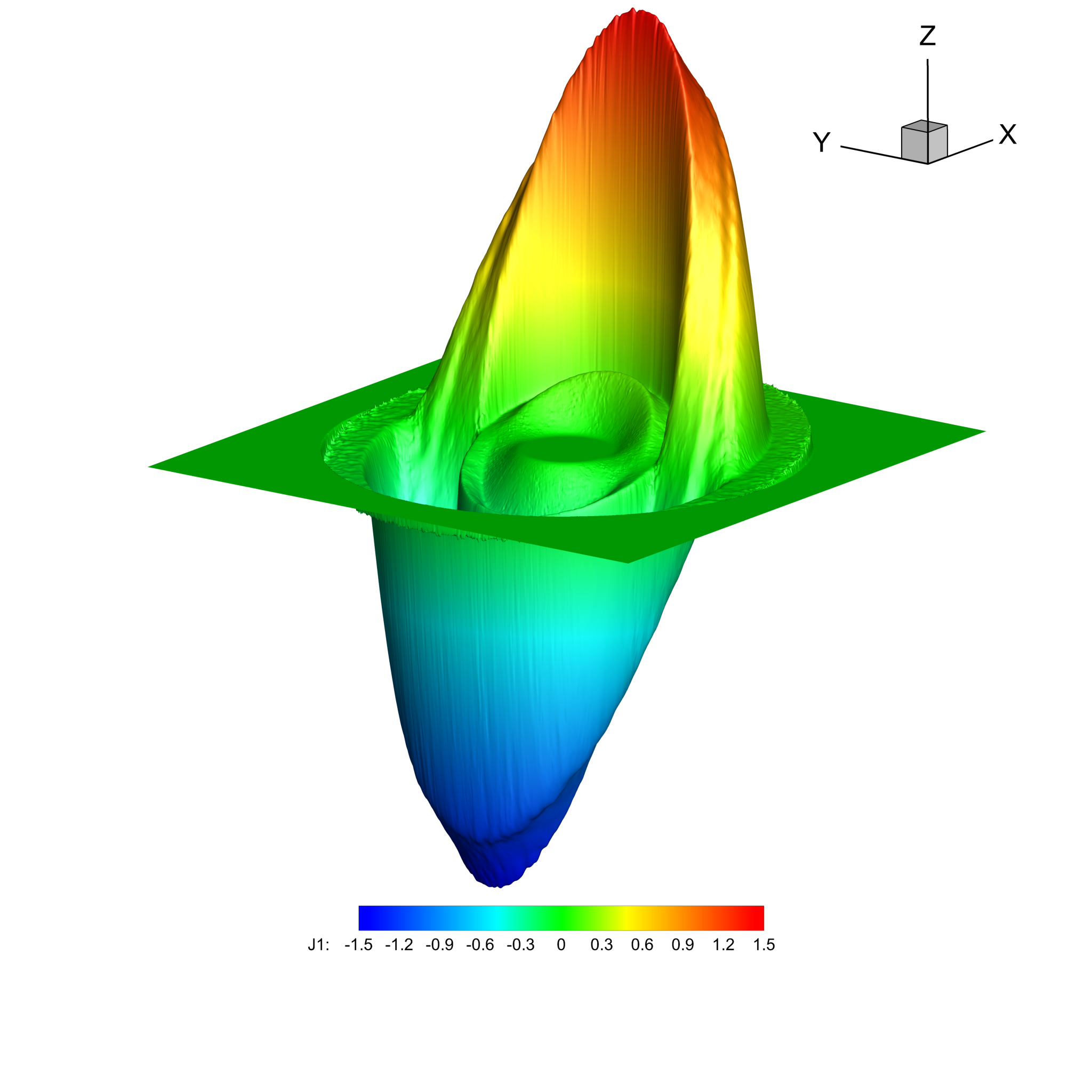}    \\ 
		\caption{Solid circular explosion. From left top to right bottom: elevated contour plots of the density, pressure, distortion component $A_{11}$ and heat flux component $J_1$.}  
		\label{fig.CEsolid3d}
	\end{center}
\end{figure}

\begin{figure}
	\begin{center}
		\includegraphics[width=0.45\textwidth]{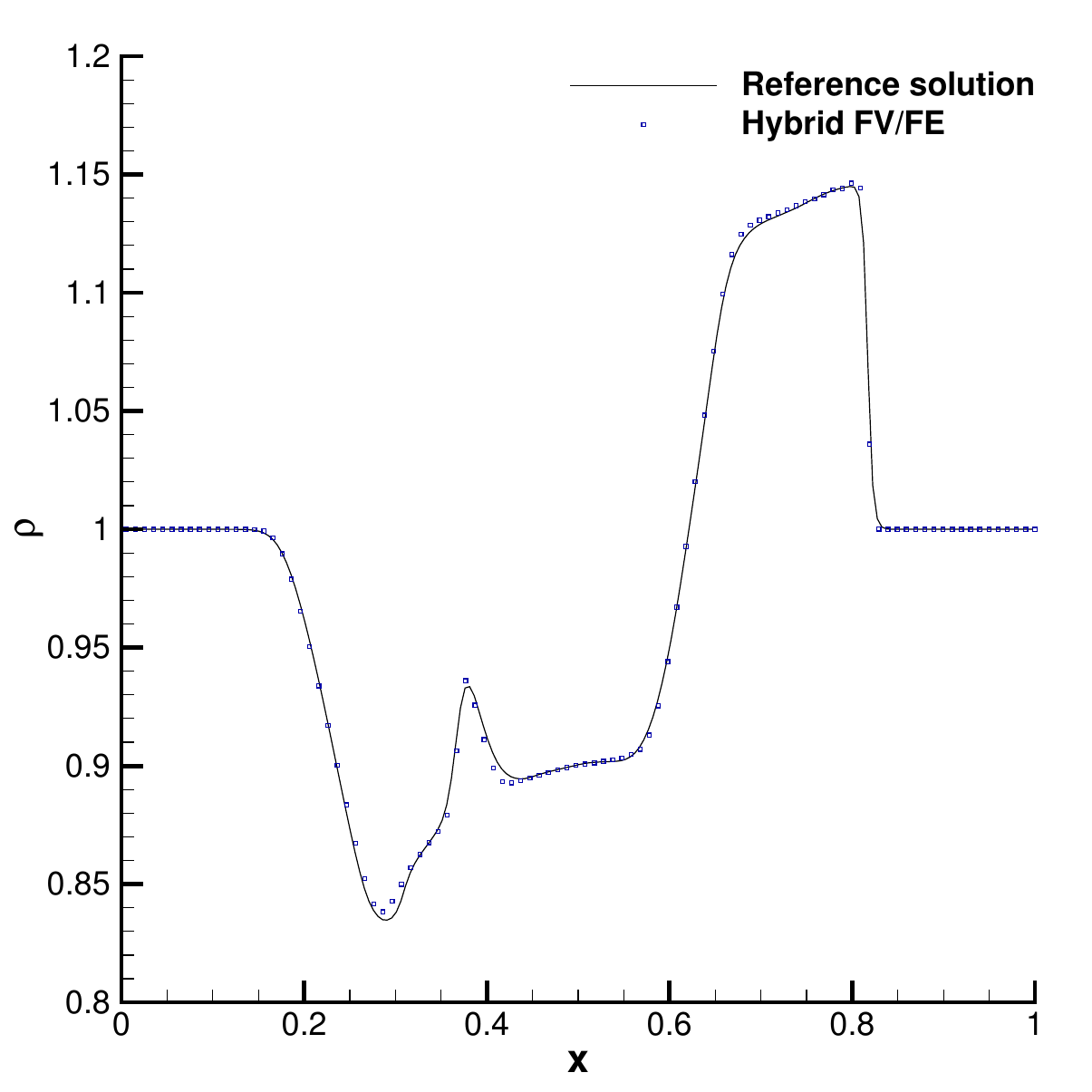}  
		\includegraphics[width=0.45\textwidth]{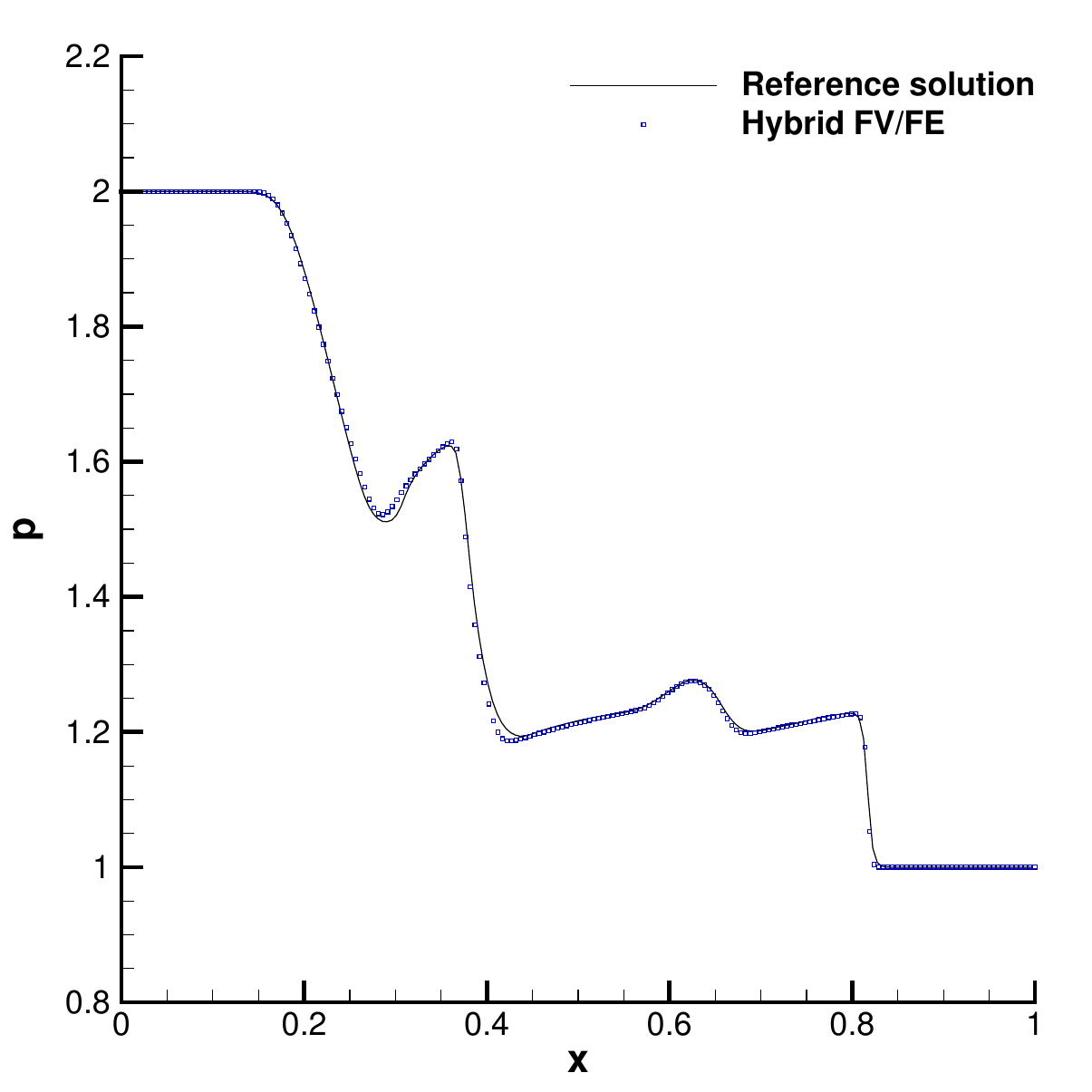}    \\ 
		\includegraphics[width=0.45\textwidth]{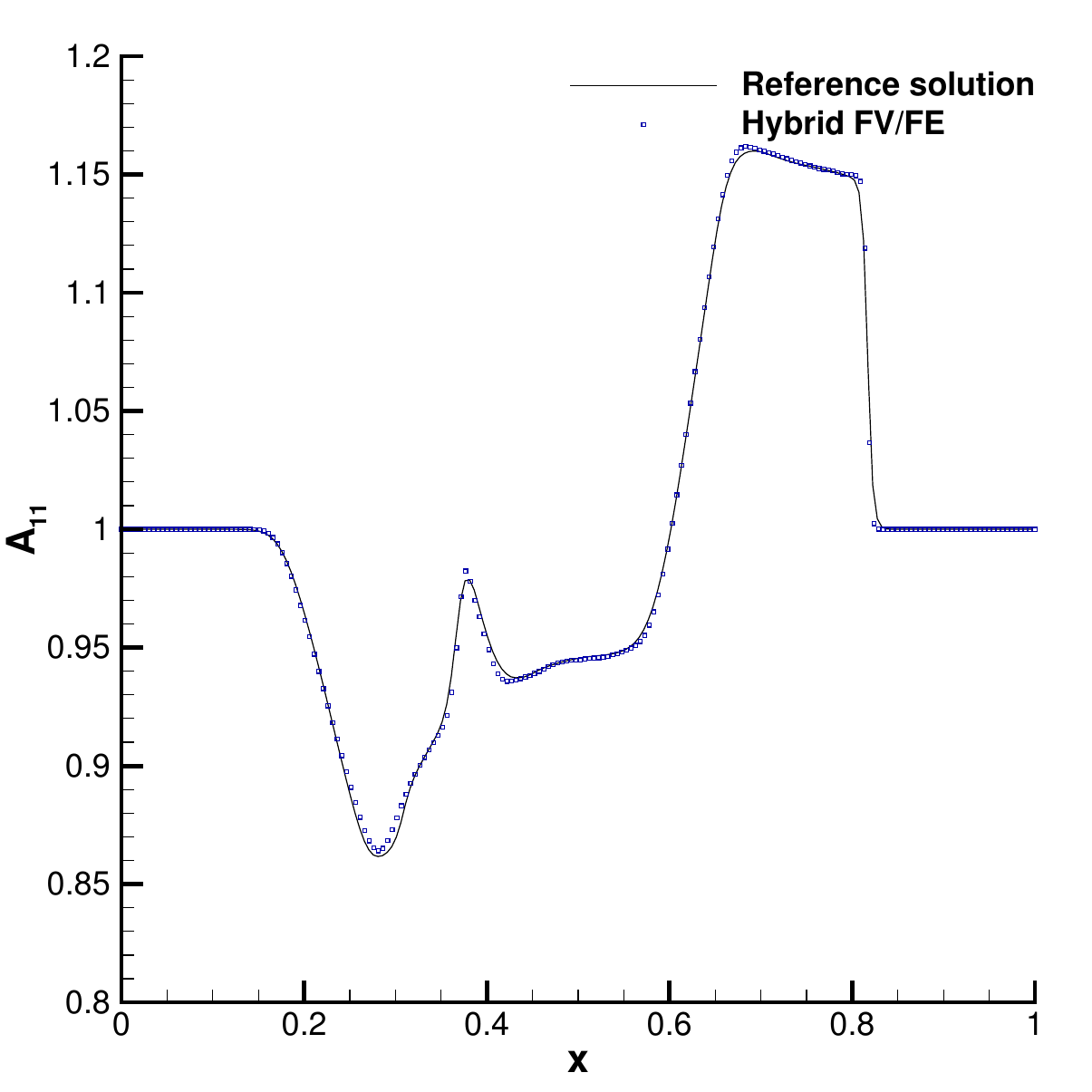}   
		\includegraphics[width=0.45\textwidth]{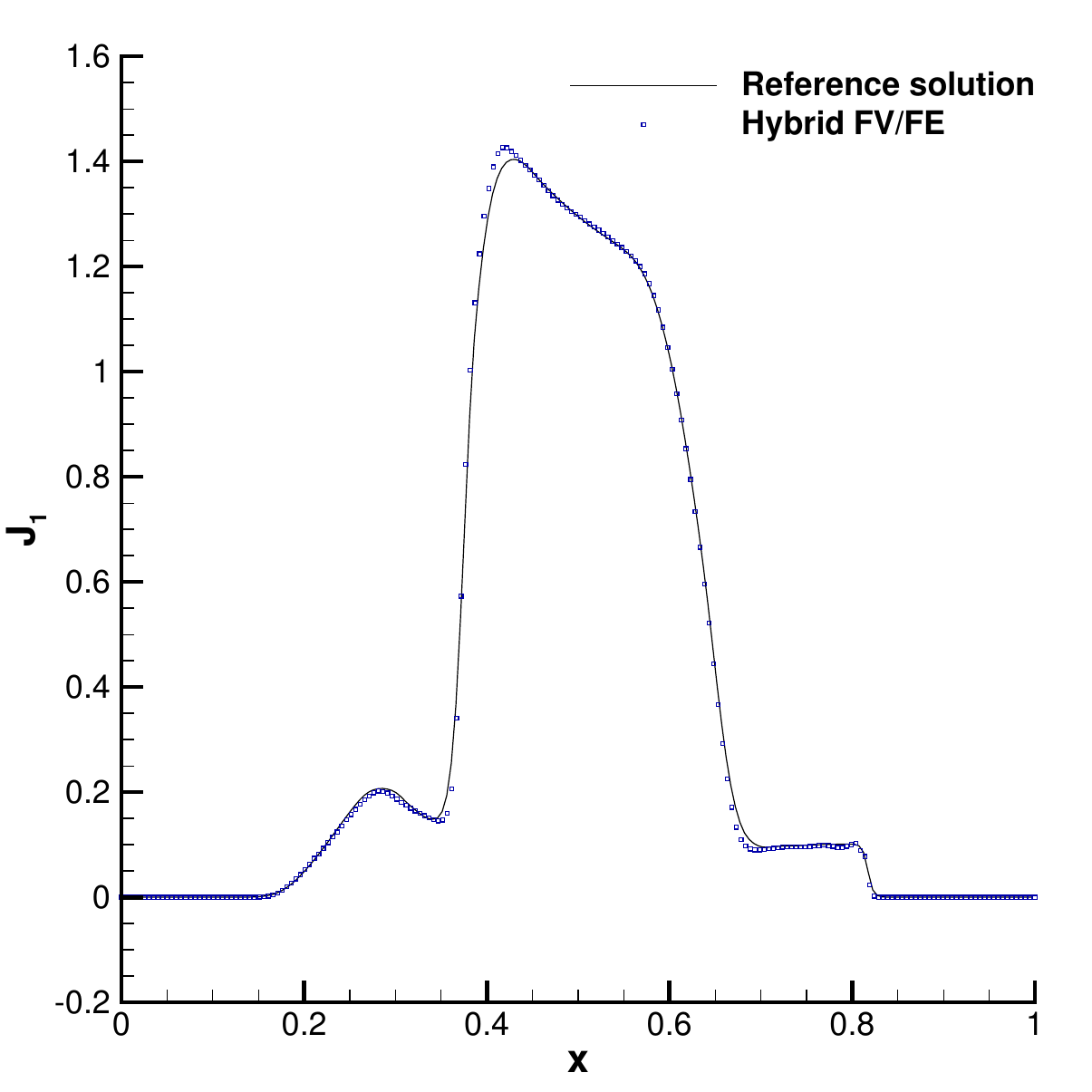}    \\ 
		\caption{Solid circular explosion. 1D cut in $x-$direction of the numerical solution obtained using the hybrid FV/FE method for the weakly compressible GPR model with the local ADER-MM approach and auxiliary artificial viscosity $c_{\alpha}=0.5$ (blue squares). Reference solution obtained with a MUSCL-Hancock FV scheme \cite{Boscheri2021SIGPR} (solid black line). From left top to right bottom: density, pressure, distortion component $A_{11}$ and heat flux component $J_1$.}  
		\label{fig.CEsolid1d}
	\end{center}
\end{figure}

To analyse the behaviour also in the solid limit, i.e., for $\tau_1\rightarrow \infty, \,\tau_2\rightarrow \infty$, we follow \cite{Boscheri2021SIGPR} and set $\tau_1=\tau_2=10^{20}$, $\rho_{0}=1$, $c_v=1.0$, $c_s=1$, $c_h=0.5$, $\gamma=1.4$ and the initial condition
\begin{equation*}
	\rho\left(\x,0\right) = 1, \;\;  \vel\left(\x,0\right)=  \boldsymbol{0},\;\; \press \left(\x,0\right) = \left\lbrace \begin{array}{lr}
		2 & \mathrm{ if } \; r \le 0.5,\\
		1 & \mathrm{ if } \; r > 0.5,
	\end{array}\right. \;\; \mathbf{A}\left(\x,0\right)=\mathbf{I}, \;\; \mathbf{J}\left(\x,0\right) = \boldsymbol{0}.
\end{equation*}
We run the simulation using the hybrid FV/FE scheme using the local ADER min-mod approach and an auxiliary artificial viscosity $c_{\alpha}=0.5$. Again, periodic boundary conditions are set in all boundaries. The results obtained at time $t=0.15$ with a fine grid of $1365504$ primal triangular elements are reported in Figure~\ref{fig.CEsolid3d}. Moreover, Figure~\ref{fig.CEsolid1d} reports the 1D cuts for the density, pressure, distortion component $A_{11}$ and heat flux component $J_1$ along a 1D cut in $x-$direction. 
An excellent agreement is observed with the reference solution computed using   
a second order MUSCL-Hancock finite volume scheme, \cite{Toro,Boscheri2021SIGPR}. 

\subsection{3D spherical explosion}
To illustrate the extension of the proposed hybrid FV/FE methodology for the weakly compressible GPR model to the three dimensional case, we study the 3D spherical explosion problem already employed in \cite{Hybrid1} to assess weakly compressible flows. The computational domain is an sphere of radius $1$ centred in $(0,0,0)$ and discretised using a primal grid of $2280182$ tetrahedra. As initial condition we set 
\begin{gather*}
	\rho\left(\x,0\right) =  \left\lbrace \begin{array}{lr}
		2 & \mathrm{ if } \; r \le 0.5,\\
		1.125 & \mathrm{ if } \; r > 0.5,
	\end{array}\right. \qquad
	\press \left(\x,0\right) = \left\lbrace \begin{array}{lr}
		2 & \mathrm{ if } \; r \le 0.5,\\
		1.1 & \mathrm{ if } \; r > 0.5,
	\end{array}\right. \notag\\
	\bvel \left(\x,0\right) = 0, 
	\qquad \mathbf{A}\left(\x,0\right)=\mathbf{I}, \qquad \mathbf{J}\left(\x,0\right) = \boldsymbol{0},
\end{gather*}
with $r$ the distance to the origin, and the model parameters are defined as $c_{s}=c_{h}=0$ and $\mu =\kappa = 0$, which correspond to the fluid limit of the model. 
\begin{figure}[H]
	\begin{center}
		\includegraphics[width=0.45\textwidth]{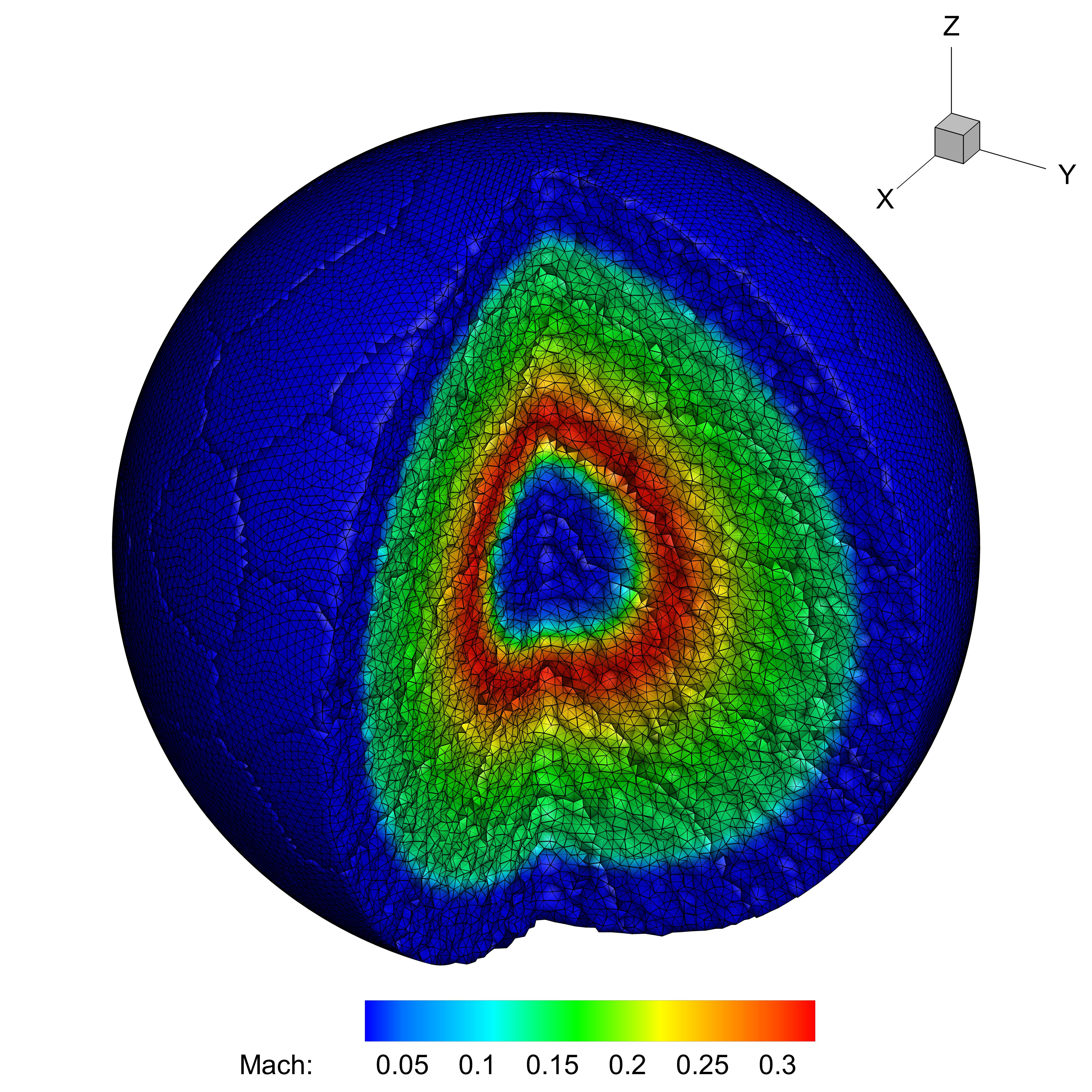}  
		\includegraphics[width=0.45\textwidth]{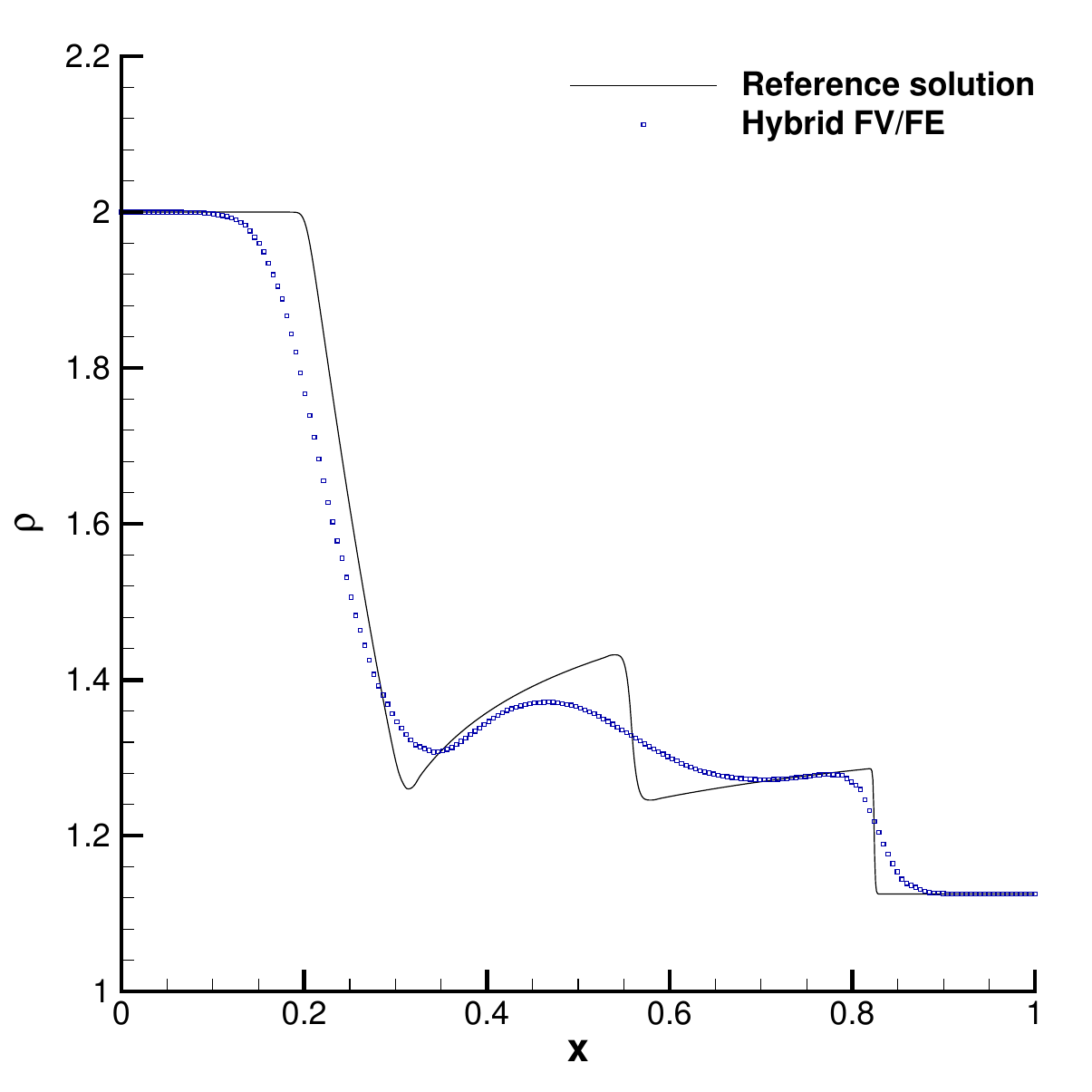}    \\ 
		\includegraphics[width=0.45\textwidth]{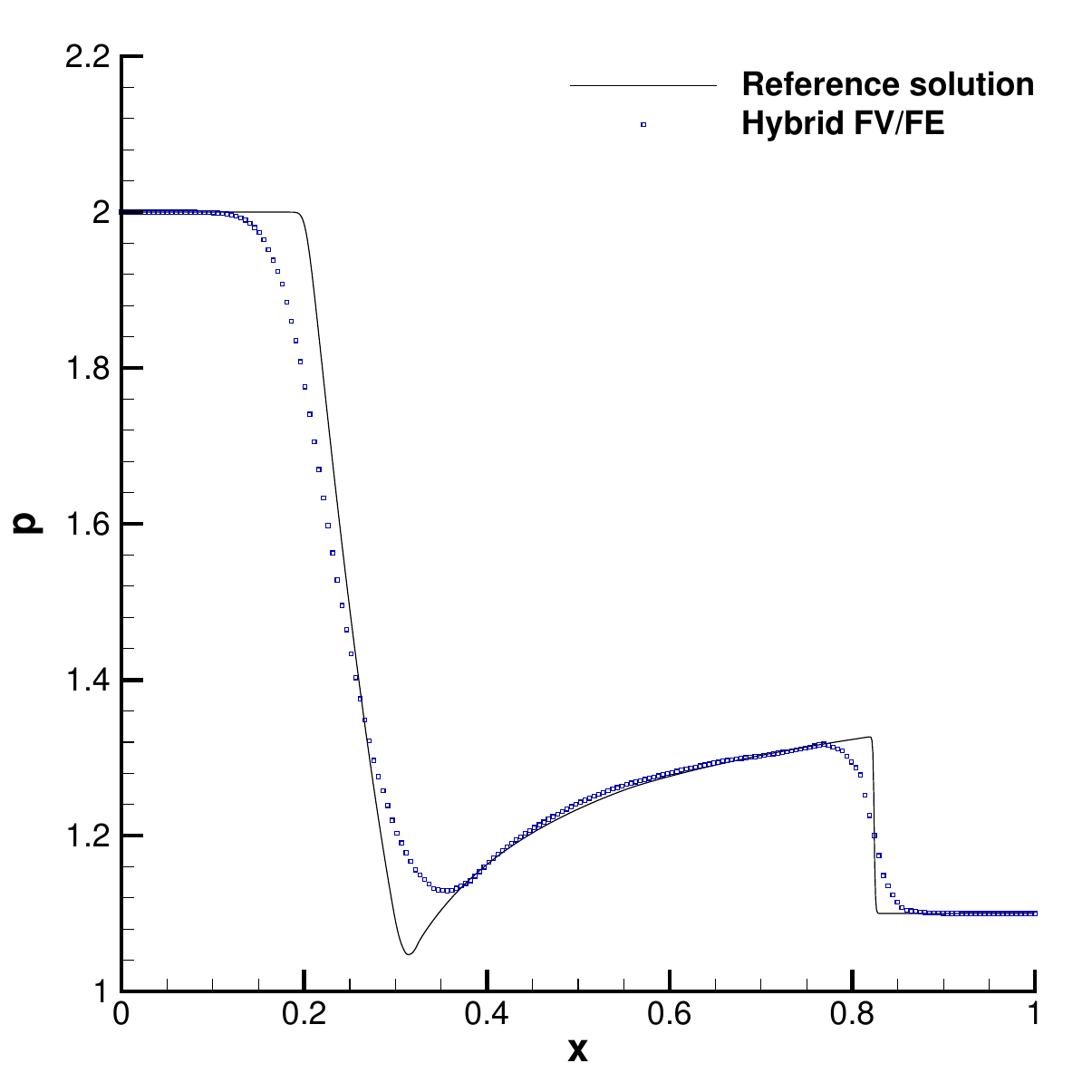}   
		\includegraphics[width=0.45\textwidth]{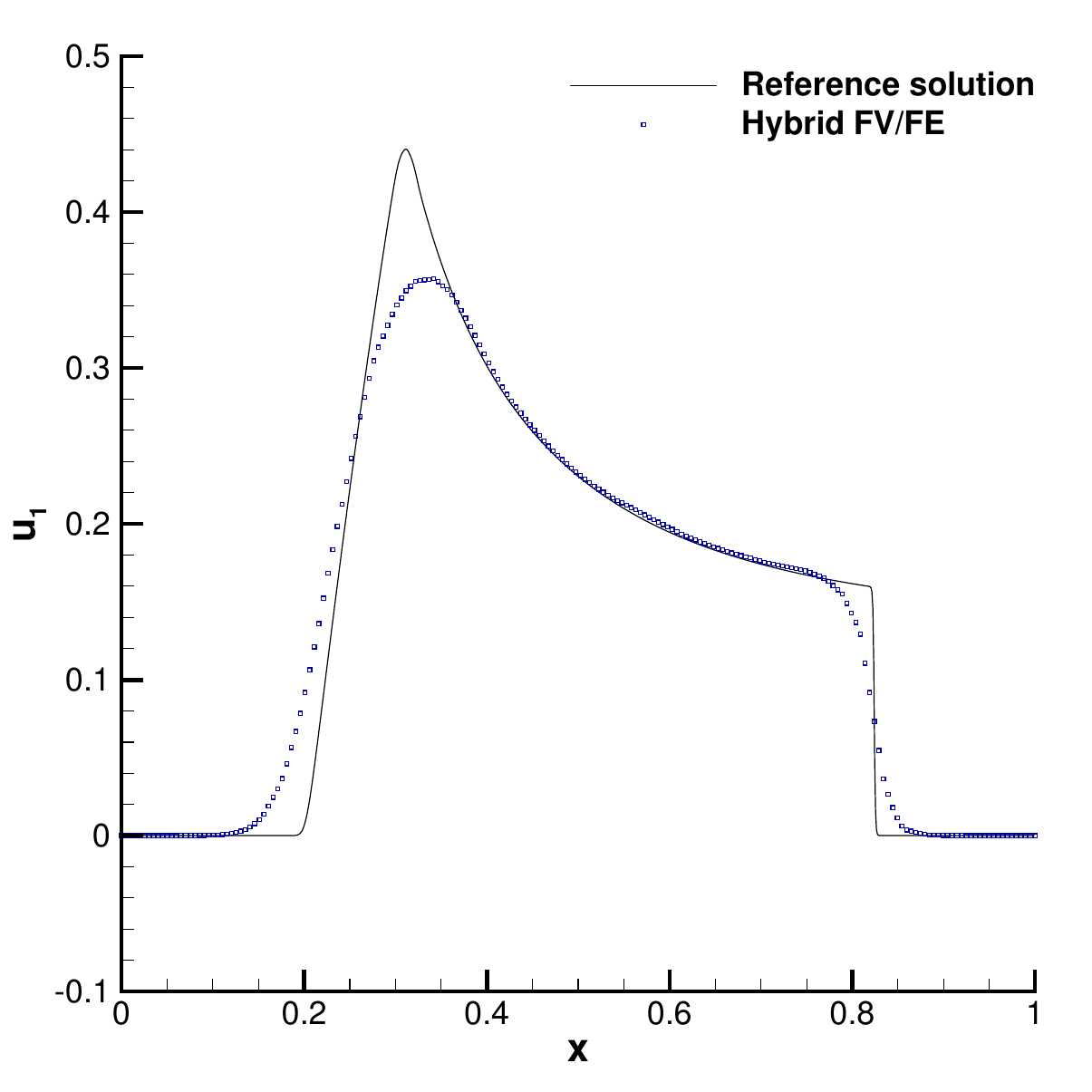}    \\ 
		\caption{3D spherical explosion. Left top: 3D mesh with MPI divisions and Mach number contours. From top right to bottom right:  1D cuts along the $x-$axis of the density, pressure, and velocity component $\vel_1$ obtained using the hybrid FV/FE method for the weakly compressible GPR model with the local ADER-ENO approach and auxiliary artificial viscosity $c_{\alpha}=3$ (blue squares). Reference solution obtained with a FV-TVD scheme solving the 1D Euler equations (solid black line).}  
		\label{fig.3DEP}
	\end{center}
\end{figure}
Dirichlet boundary conditions are imposed and the numerical simulation is carried out up to time $t=0.25$. Figure~\ref{fig.3DEP} reports the results obtained using the LADER-ENO approach for the convective terms and auxiliary artificial viscosity $c_{\alpha}=3$. We observe a good agreement with the reference solution computed using a TVD-FV scheme for the 1D compressible Euler equations with appropriate geometrical source terms, \cite{Toro}.

\subsection{Smooth acoustic wave}
\begin{figure}[h!]
	\centering
	\includegraphics[clip,trim= 200 0 10 0,width=0.45\linewidth]{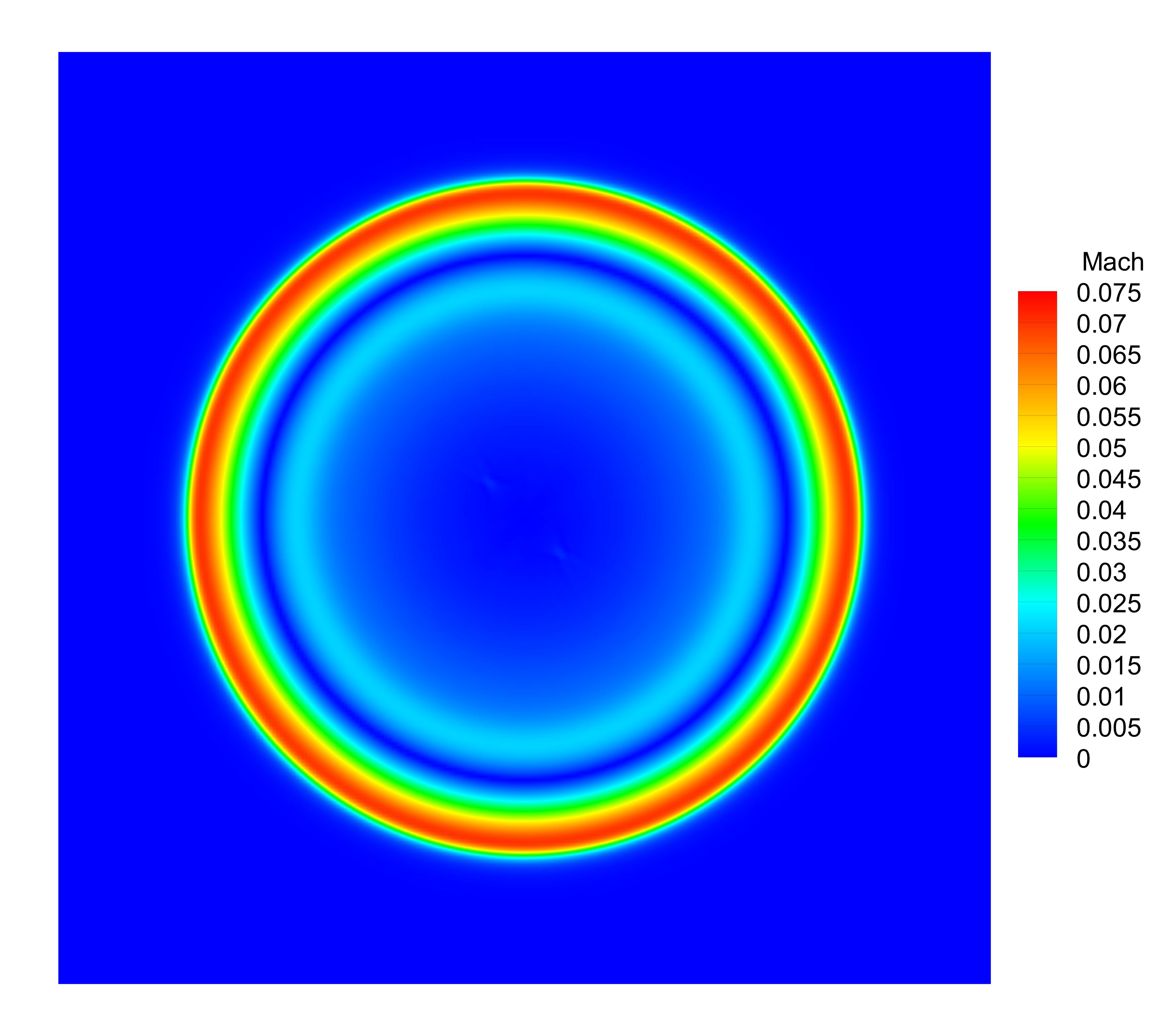}\hspace{0.05\linewidth}
	\includegraphics[width=0.45\linewidth]{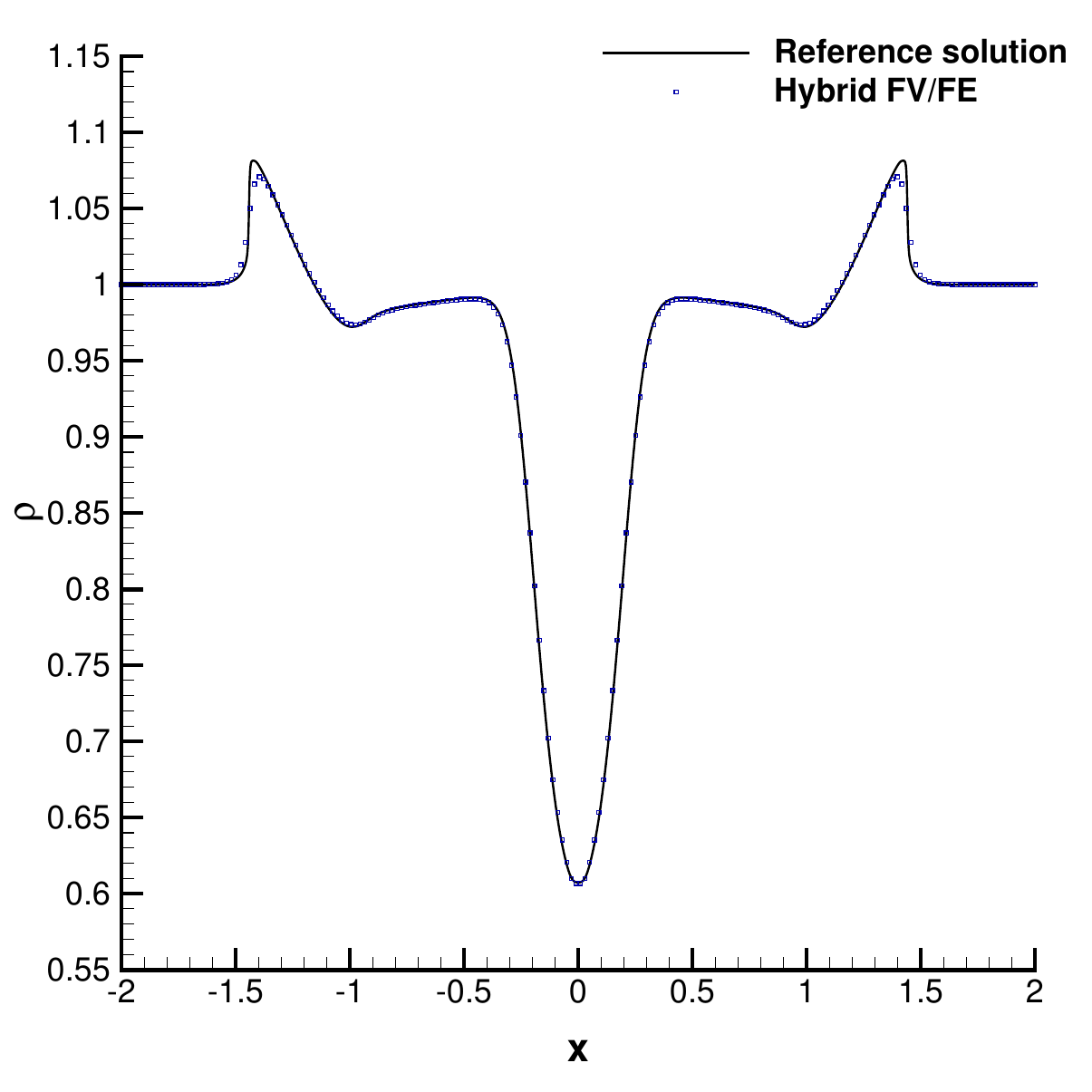}
	
	\vspace{0.05\linewidth}
	\includegraphics[width=0.45\linewidth]{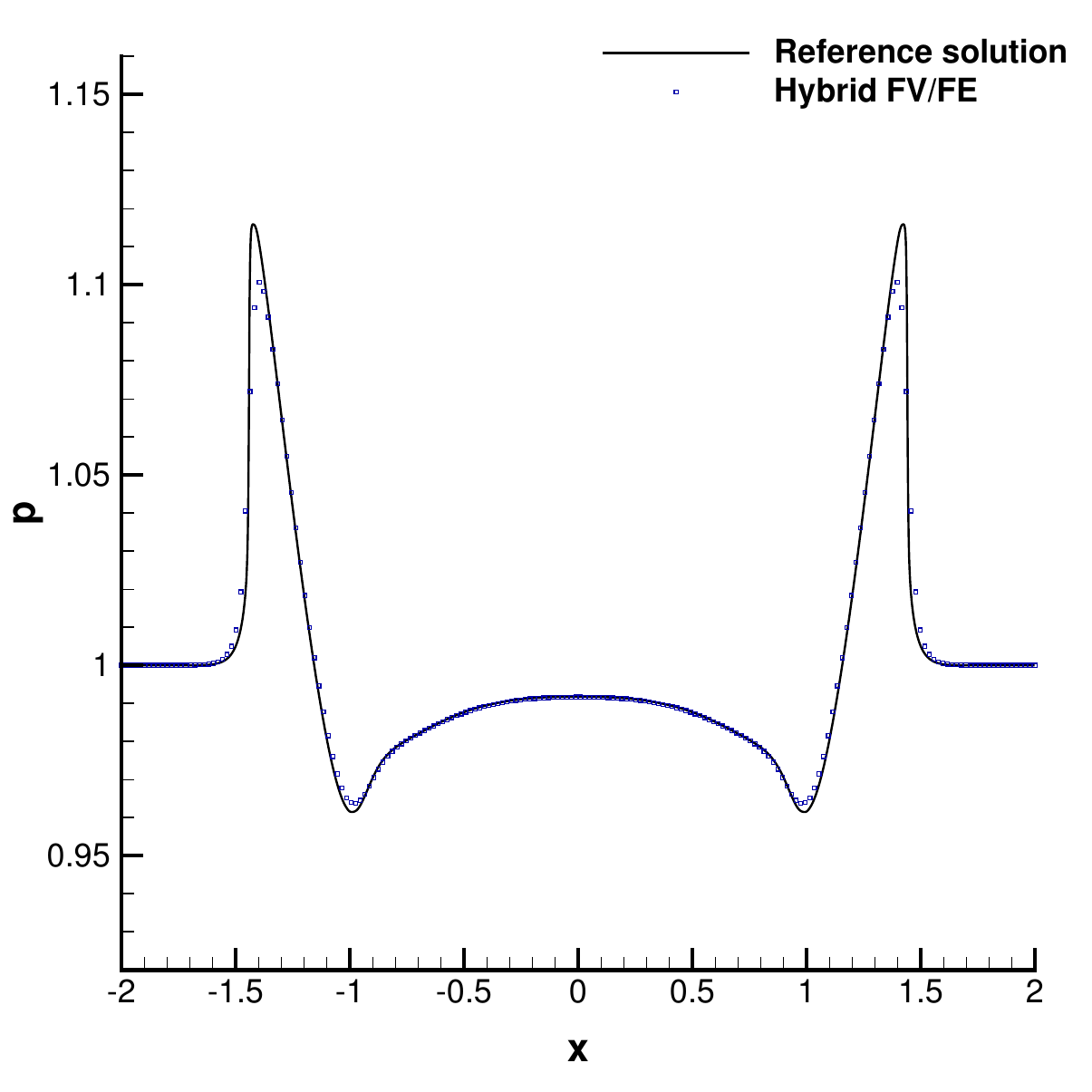}\hspace{0.05\linewidth}
	\includegraphics[width=0.45\linewidth]{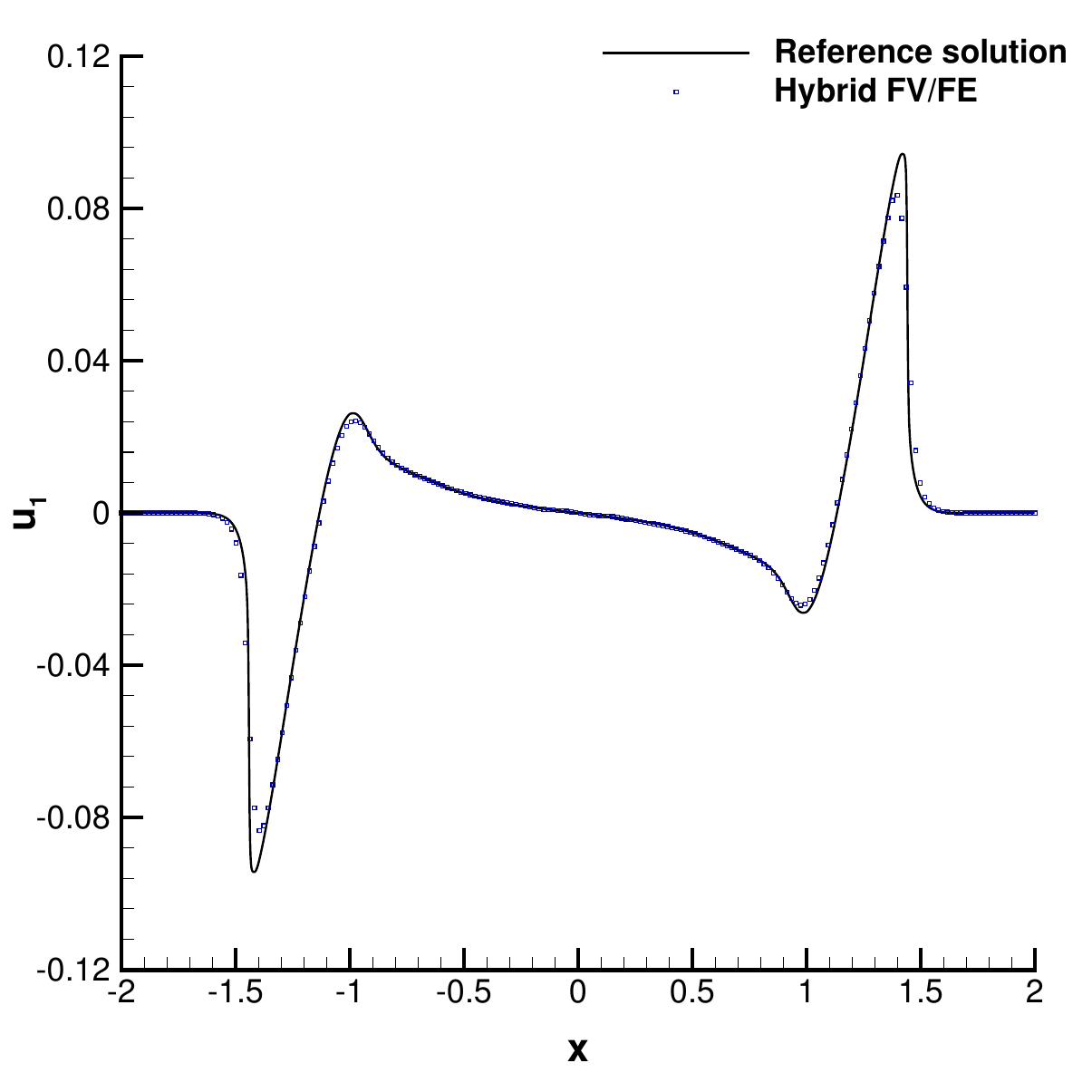}
	\caption{Smooth acoustic wave. Left top: 2D Mach number contour plot obtained at time $t=1$. From right top to right bottom: 1D cuts along $y=0$ of the density, pressure and velocity component $u_{1}$ (blue squares). Reference solution computed with the TVD-FV scheme on a grid of $10^4$ cells (black solid line).
	}
	\label{SAW256}
\end{figure}

One important difference between weakly compressible and incompressible flows is the presence of acoustic waves. To analyse the ability of the proposed semi-implicit FV/FE approach to capture acoustic waves properly we consider the smooth acoustic wave benchmark in \cite{TD17,Hybrid1}. The initial condition is defined as
\begin{gather*}
	\rho\left(\x,0\right) = 1,\;\;
	\bvel \left(\x,0\right) = 0, \;\;
	\press \left(\x,0\right) = 1+e^{-\alpha r^2}, 
	\;\; \mathbf{A} = \mathbf{I}, 
	\;\; \mathbf{J}=\boldsymbol{0},
	\;\; r = \sqrt{x^2+y^2}
\end{gather*}
and the model parameters are set to $\mu=\kappa = c_{s}=c_{h}=0$.  
The computational domain $\Omega=[0,2]^2$ is discretised employing a primal triangular grid formed by $131072$ cells, and periodic boundary conditions are defined everywhere. The second order hybrid FV/FE scheme with ENO limiters is employed to get the solution at time $t=1$. To generate a reference solution, we consider the 1D PDE in radial direction with geometrical source terms equivalent to the compressible Euler system, which is solved using a second order TVD-FV scheme on a grid of $10^4$ cells. Figure~\ref{SAW256} shows an excellent agreement between both numerical results. Let us remark that even if this test case is characterised by a low Mach number, the compressibility still plays a primal role since we observe a steep acoustic wavefront propagating in radial direction. Moreover, using a semi-implicit approach leads to a CFL number depending only on the bulk flow velocity, so we circumvent the strong time step restriction of explicit Godunov-type solvers, which is related to the sound speed. 

\section{Conclusions} \label{sec:conclusions}
We have presented a novel hybrid FV/FE methodology for the solution of the GPR model for continuum mechanics on unstructured meshes. From the proposed mathematical model and to simulate weakly compressible flows, we have derived two new formulations: the incompressible GPR model and a weakly compressible GPR model. Moreover, as for the original GPR model, able to address all Mach number flows appropriately setting the model parameters, we can simulate both solids with large deformations and fluids. To discretise these systems, a splitting of the equations is performed, leading a Poisson-type pressure system and a transport system containing convective terms and non-conservative products. This last system is solved using finite volume methods. Even if this part of the scheme is explicit, it is independent of the fast sound velocity waves, yielding a computationally efficient scheme in the low Mach regime. Moreover, to avoid the severe time-step restriction that may arise from the presence of stiff source terms in the distortion field and heat flux equations, an implicit finite volume method is employed for their discretisation. On the other hand, the pressure subsystem is solved using classical finite element methods, which are well known for their efficiency in solving Poisson-type problems. Finally, the intermediate momentum field computed within the transport stage of the algorithm is corrected to account for the new pressures, thus providing the momentum at the new time step. The final methodology has been successfully assessed employing a wide range of test cases, from solid mechanics benchmarks to the incompressible fluid limit of the equations, including the analysis of low Mach problems featuring small shocks. 

In future, we plan to extend the former hybrid methodology to deal with fluid-structure iteration problems. To this end, following the methodology in \cite{HybridALE} for the Navier-Stokes equations, the hybrid FV/FE method for the GPR model will also be developed in the ALE framework. Moreover, since the proposed scheme is at most second order accurate, we plan to extend the methodology to high order accuracy using DG schemes and IMEX methods, \cite{TCRGD20,HybridFVVEMinc}. Finally, we would like to study the development of hybrid methods verifying additional properties of the physical model, such as preserving the involution constraints at the discrete level and designing a thermodynamically compatible scheme on unstructured grids, \cite{HybridHexa1,HybridMHD,HTCAbgrall,HTCTotalAbgrall}.

\section*{Acknowledgements}
SB acknowledges support from the Spanish Ministry of Science, Innovation and Universities (MCIN), the Spanish AEI (MCIN/AEI/10.13039/501100011033) and European Social Fund Plus under the project No. RYC2022-036355-I; from FEDER and the Spanish Ministry of Science, Innovation and Universities under project No. PID2021-122625OB-I00; and from the Xunta de Galicia (Spain) under project No. GI-1563 ED431C 2021/15. LR acknowledges the support from the Italian Ministry of Education, University and Research (MIUR) in the frame of the Departments of Excellence Initiative 2018--2027 attributed to DICAM of the University of Trento (grant L. 232/2016) and in the frame of the PRIN 2022 project \textit{High order structure-preserving semi-implicit schemes for hyperbolic equations}. LR is member of the Gruppo Nazionale Calcolo Scientifico-Istituto Nazionale di Alta Matematica (GNCS-INdAM). The authors would like to acknowledge support from the CESGA, Spain, for the access to the FT3 supercomputer and to the CINECA award under the ISCRA initiative, for the availability of high performance computing resources and support (project IsB27\_NeMesiS). 

 Views and opinions expressed are however those of the author(s) only and do not necessarily reflect those of the European Union. Neither the European Union nor the granting authorities can be held responsible for them.

\bibliographystyle{abbrv}
\bibliography{./mibiblio}

\end{document}